\documentclass[12pt]{article}

\usepackage[utf8]{inputenc}
\usepackage[T1]{fontenc}

\usepackage{amsmath,
            mleftright,
            comment,
            amssymb,
            amsthm,
            nicefrac,
            xcolor,
            enumerate,
            authblk,
            color,
            etoolbox,
            mathrsfs,
            mathtools,
            bbm,
            xparse,
            stmaryrd,
            geometry,
            enumitem,
            mathtools,
            marginnote,
            tikz,
            xurl,
            graphicx,
            scalerel,
            relsize,
            stackengine
            }

\mleftright

\newcommand{\medcup}[1]{{\raisebox{0.25ex}{$\mathsmaller{\ensuremath{\bigcup}}_{#1}$}}}

\makeatletter
\newcommand{\medtimes}{\vcenter{\Large\hbox{$ \m@th\mkern-2mu\times\mkern-2mu$}}}
\makeatother

\usetikzlibrary{tikzmark,calc}

\usepackage[skins,breakable]{tcolorbox}

\usepackage[colorinlistoftodos]{todonotes}
\usepackage[english]{datetime2}

\usepackage[nocompress,sort]{cite}
                       
\usepackage[colorlinks=true]{hyperref}
            
\usepackage[sort,capitalize]{cleveref}

\crefformat{equation}{(#2#1#3)}
\crefname{enumi}{item}{items}
\crefname{equation}{}{}

\theoremstyle{plain}
\newtheorem{theorem}{Theorem}[section]

\theoremstyle{remark}

\theoremstyle{definition}
\newtheorem{definition}[theorem]{Definition}

\numberwithin{equation}{section}

\tikzset{notestyleraw/.append style={align=justify}}

\newcounter{todocounter}

\DeclareFontEncoding{LS1}{}{}
\DeclareFontSubstitution{LS1}{stix}{m}{n}
\DeclareMathAlphabet{\mathscr}{LS1}{stixscr}{m}{n}

\newcommand{\cB}{\mathcal{B}}

\newcommand{\cF}{\mathcal{F}}

\newcommand{\cU}{\mathcal{U}}

\newcommand{\fB}{\mathfrak{B}}
\newcommand{\fC}{\mathfrak{C}}

\newcommand{\fJ}{{\bf\mathfrak{J}}}

\newcommand{\fL}{\mathfrak{L}}
\newcommand{\fM}{\mathfrak{M} \cfadd{multi}}

\newcommand{\fR}{\mathfrak{R}}

\newcommand{\fa}{\mathfrak{a}}

\newcommand{\fc}{\mathfrak{c}}
\newcommand{\fd}{\mathfrak{d}}

\newcommand{\fh}{\mathfrak{h}}

\newcommand{\fp}{\mathfrak{p}}
\newcommand{\fq}{\mathfrak{q}}

\newcommand{\fu}{\mathfrak{u}}

\newcommand{\fx}{\mathfrak{x}}

\newcommand{\scrL}{\mathscr{L} \cfadd{def:lin_interp}}

\newcommand{\LL}{\mathbb{L}}

\newcommand{\induct}{\dashrightarrow}
\newcommand{\with}{\curvearrowleft}
\newcommand{\lrSpace}{\ensuremath{\mkern-1.5mu}}

\newcommand{\eps}{\varepsilon}
\newcommand{\eg}{\unskip, e.g.,\ }

\newcommand{\dpp}{\text{.}}
\newcommand{\dc}{\text{,}}
\newcommand{\dx}{\, {\rm d}}
\newcommand{\dxx}{{\rm d}}

\newcommand{\normmm}[1]{{\left\vert\kern-0.25ex\left\vert\kern-0.25ex\left\vert #1 
    \right\vert\kern-0.25ex\right\vert\kern-0.25ex\right\vert}
    \cfadd{DNN_norm}}

\DeclarePairedDelimiter{\pr}{(}{)}
\DeclarePairedDelimiter{\br}{[}{]}
\DeclarePairedDelimiter{\cu}{\{}{\}}
\DeclarePairedDelimiter{\abs}{\lvert}{\rvert}
\DeclarePairedDelimiter{\norm}{\lVert}{\rVert\cfadd{DNN_norm}}

\newcommand{\N}{\ensuremath{\mathbb N}}
\newcommand{\Z}{\ensuremath{\mathbb Z}}

\newcommand{\R}{\ensuremath{\mathbb R}}

\newcommand{\E}{\ensuremath{\mathbb E}}
\renewcommand{\P}{\ensuremath{\mathbb P}}

\newcommand{\parallelizationSpecial}{\mathbf{P} \cfadd{def:simpleParallelization}}
\newcommand{\pa}[1]{\left({#1}\right)}

\newcommand{\Ra}{\mathcal{R}_a \cfadd{def:ANNrealization}}
\newcommand{\1}{\mathbbm{1}}
\newcommand{\ii}{\mathfrak{i} \cfadd{padding}}

\newcommand{\network}{ANN \cfadd{def:ANN2}}
\newcommand{\interpolatingDNN}{\mathbf{F}}
\newcommand{\AffineANN}{\cfadd{linear}\mathbf{A}}
\newcommand{\bbigANNsum}{\cfadd{def:ANNsum:same}\mathop{\oplus}\limits}
\newcommand{\scalarMultANN}[2]{\cfadd{def:ANNscalar}#1\circledast#2}
\newcommand{\interpol}[2]{\cfadd{def:lin_interp}\mathscr{L}_{#1}^{#2}}

\newcommand{\fwpr}{W}
\newcommand{\smallU}{u}
\newcommand{\smallF}{f}

\newcommand{\LipConstF}{L}
\newcommand{\mlp}{U}

\newcommand{\indexAct}{\nu}
\newcommand{\constantAssumpMainThm}{\kappa}
\newcommand{\ANNassumpMainThm}{\mathbf{G}}
\newcommand{\ANNresultMainThm}{\mathbf{U}}
\newcommand{\fsc}{\mathscr{f}}
\newcommand{\usc}{\mathscr{u}}
\newcommand{\relu}{\mathfrak{r}}

\renewcommand{\fR}{\mathscr{n}}

\newcommand{\littleM}{m}
\newcommand{\littleMM}{m}

\newcommand{\A}{\mathbf{A} \cfadd{linear}}
\newcommand{\U}{\mathbf{U}}

\newcommand{\F}{\mathbf{F}}
\newcommand{\G}{\mathbf{G}}
\newcommand{\I}{\mathbf{I} \cfadd{def:id_net}}

\newcommand{\NN}{\mathbf{N}}
\newcommand{\lbd}{u}
\newcommand{\ubd}{v}

\newcommand{\weight}{W \cfadd{linear}}
\newcommand{\bias}{B \cfadd{linear}}

\newcommand{\andShort}{\text{ and }}

\newcommand{\smallsum}{\textstyle\sum}
\newcommand{\SmallSum}[2]{ {\textstyle\sum\limits_{#1}^{#2}}}

\newcommand{\idMatrix}{\operatorname{I} \cfadd{def:identityMatrix}}
\newcommand{\ANNs}{\mathbf{N} \cfadd{def:ANN}}

\newcommand{\activation}{a}
\newcommand{\activationDim}[1]{\mathfrak{M}_{\activation,#1} \cfadd{multi}}
\newcommand{\functionANN}[1]{\mathcal{R}_{#1} \cfadd{def:ANNrealization}}
\newcommand{\paramANN}{\mathcal{P} \cfadd{def:ANN}}

\newcommand{\lengthANN}{\mathcal{L} \cfadd{def:ANN}}
\newcommand{\inDimANN}{\mathcal{I} \cfadd{def:ANN}}
\newcommand{\compANN}[2]{{#1 \bullet \allowbreak #2} \cfadd{def:ANNcomposition}}

\newcommand{\outDimANN}{\mathcal{O} \cfadd{def:ANN}}
\newcommand{\longerANN}[1]{\mathcal{E}_{#1} \cfadd{def:ANNenlargement}}

\newcommand{\dims}{\mathcal{D} \cfadd{def:ANN}}
\newcommand{\hiddenLength}{\mathcal{H} \cfadd{def:ANN}}

\newcommand{\qandq}{\qquad\text{and}\qquad}

\newcommand{\sumANN}{\mathfrak{S} \cfadd{def:ANN:sum}}
\newcommand{\extensionANN}{\mathfrak{T} \cfadd{def:ANN:extension}}
\newcommand{\dimANNlevel}{\mathbb{D} \cfadd{def:ANN}}
\newcommand{\power}[2]{#1^{\bullet #2} \cfadd{def:iteratedANNcomposition}}
\newcommand{\scalar}[2]{ #1 \circledast #2 \cfadd{def:ANNscalar}}
\newcommand{\transpose}{* \cfadd{def:Transpose}}
\newcommand{\oSum}{\oplus \cfadd{def:ANNsum:same}}
\newcommand{\OSum}[2]{{\mathop\oplus\limits_{#1}^{#2} } \cfadd{def:ANNsum:same} }
\newcommand{\bSum}{{\mathop\boxplus} \cfadd{def:ANN:sum_diff}}
\newcommand{\BSum}[3]{{\mathop\boxplus\limits_{#1,#2}^{#3}} \cfadd{def:ANN:sum_diff}}
\newcommand{\modCont}{w \cfadd{mod_cont_def}}

\allowdisplaybreaks

\ExplSyntaxOn

\NewDocumentCommand{\enum}{ O{;} m o }
 {
  \my_enum:nnn { #1 } { #2 } { #3 }
 }

\seq_new:N \l__my_enum_seq
\tl_new:N \l__my_enum_item_tl
\prop_new:N \l__verbs
\prop_put:Nnn \l__verbs {show} {shows}
\prop_put:Nnn \l__verbs {imply} {implies}
\prop_put:Nnn \l__verbs {demonstrate} {demonstrates}
\prop_put:Nnn \l__verbs {prove} {proves}
\prop_put:Nnn \l__verbs {establish} {establishes}
\prop_put:Nnn \l__verbs {ensure} {ensures}
\prop_put:Nnn \l__verbs {assure} {assures}
\prop_put:Nnn \l__verbs {yield} {yields}

\int_new:N \l__number_of_args

\cs_new_protected:Nn \my_enum:nnn
 {
  \seq_set_split:Nnn \l__my_enum_seq { #1 } { #2 }
  \seq_remove_all:Nn \l__my_enum_seq {}
  \int_set_eq:NN \l__number_of_args { \seq_count:N \l__my_enum_seq }
  \seq_use:Nnnn \l__my_enum_seq { ~and~ } { ,~ } { ,~and~ }
  \IfNoValueTF{#3}{}{
    \space
    \int_compare:nNnTF{ \l__number_of_args } < {2}{ \prop_item:Nn \l__verbs {#3} }{ #3 }
  }
 }

\ExplSyntaxOff

\ExplSyntaxOn

\seq_new:N \g_cflist_loaded

\seq_new:N \g_cflist_pending

\NewDocumentCommand{\cfadd}{ m }

{

  \seq_if_in:NnF \g_cflist_loaded { #1 } {

    \seq_if_in:NnF \g_cflist_pending { #1 } {

      \seq_gput_right:Nn \g_cflist_pending { #1 }

    }

  }

}

\NewDocumentCommand{\cfload}{ o }

{

  \seq_if_empty:NTF \g_cflist_pending {\unskip} {

    (cf.\ \cref{\seq_use:Nn \g_cflist_pending {,}})\IfValueTF{#1}{#1~}{\unskip}

    \seq_gconcat:NNN \g_cflist_loaded \g_cflist_loaded \g_cflist_pending

    \seq_gclear:N \g_cflist_pending

  }

}

\NewDocumentCommand{\cfclear} {} {

  \seq_gclear:N \g_cflist_loaded

  \seq_gclear:N \g_cflist_pending

}

\NewDocumentCommand{\cfout}{ o }

{

  \seq_if_empty:NTF \g_cflist_pending {\unskip} {

    (cf.\ \cref{\seq_use:Nn \g_cflist_pending {,}})\IfValueTF{#1}{#1~}{\unskip}

    \seq_gclear:N \g_cflist_pending

  }

}

\NewDocumentCommand{\ifnocf} { m } {

  \seq_if_empty:NT \g_cflist_pending { #1 }

}

\NewDocumentCommand{\cfconsiderloaded}{ m }{

  \seq_gput_right:Nn \g_cflist_loaded {#1}

}

\ExplSyntaxOff
 
\ExplSyntaxOn

\NewDocumentEnvironment {athm} {m m o} {
\str_if_eq:noTF {example} {#1} {
  \bool_gset_true:N \g_example_bool
} {
  \bool_gset_false:N \g_example_bool
}
\IfNoValueTF{#3}{
\begin{#1}\label{#2}\global\def\loc{#2}
}{
\begin{#1}[#3]\label{#2}\global\def\loc{#2}
}
}{
\end{#1}
}

\NewDocumentEnvironment {adef} {m} {
\begin{definition}\label{#1}\global\def\loc{#1}
}{
\end{definition}
}

\NewDocumentEnvironment{aproof} {} {
\bool_if:NTF \g_example_bool {
  \begin{proof}[Proof~for~\cref{\loc}]
} {
  \begin{proof}[Proof~of~\cref{\loc}]
}
\bool_gset_false:N \g_finishproof_bool
}{
\bool_if:NTF \g_finishproof_bool {}
{\finishproofthus}
\end{proof}
}

\NewDocumentCommand{\finishproofthus} {} {
  \bool_gset_true:N \g_finishproof_bool 
  \bool_if:NTF \g_example_bool {
    The~proof~for~\cref{\loc}~is~thus~complete.
  } {
    The~proof~of~\cref{\loc}~is~thus~complete.
  }
}
\NewDocumentCommand{\finishproofthis} {} {
  \bool_gset_true:N \g_finishproof_bool 
  \bool_if:NTF \g_example_bool {
    This~completes~the~proof~for~\cref{\loc}.
  } {
    This~completes~the~proof~of~\cref{\loc}.
  }
}

\ExplSyntaxOff

\ExplSyntaxOn

\NewDocumentEnvironment{flexmath}{ m o }{
  \str_if_eq:noTF {a} {#1} {
    \begin{equation}
    \IfValueT{#2}{\label{eq:\loc.#2}}
    \begin{aligned}
  } {
    \catcode`&=9
    \renewcommand{\\}{}
    \str_if_eq:noTF {d} {#1} {
      \begin{equation}
      \IfValueT{#2}{\label{eq:\loc.#2}}
    } {
      \begin{math}
    }
  }
}{
  \str_if_eq:noTF {i} {#1} {
    \end{math}
    \catcode`&=4
  } {
    \str_if_eq:noTF {d} {#1} {
      \end{equation}
    } {
      \end{aligned}
      \end{equation}
    }
  }
}

\ExplSyntaxOff

\ExplSyntaxOn

\bool_new:N \g_noteobserve

\NewDocumentCommand{\setnote}{}{
  \bool_gset_true:N \g_noteobserve
}

\NewDocumentCommand{\setobserve}{}{
  \bool_gset_false:N \g_noteobserve
}

\NewDocumentCommand{\nobs}{ o }{
  \IfValueT{#1}{
    \str_if_eq:noTF {note} {#1} {
      \bool_gset_true:N \g_noteobserve
    } {
      \str_if_eq:noTF {Note} {#1} {
        \bool_gset_true:N \g_noteobserve
      } {
        \bool_gset_false:N \g_noteobserve
      }
    }
  }
  \bool_if:nTF { \g_noteobserve } {
    \bool_gset_false:N \g_noteobserve
    note
  } {
    \bool_gset_true:N \g_noteobserve
    observe
  }
  \IfValueF{#1}{~}
}

\NewDocumentCommand{\Nobs}{ o }{
  \IfValueT{#1}{
    \str_if_eq:noTF {note} {#1} {
      \bool_gset_true:N \g_noteobserve
    } {
      \str_if_eq:noTF {Note} {#1} {
        \bool_gset_true:N \g_noteobserve
      } {
        \bool_gset_false:N \g_noteobserve
      }
    }
  }
  \bool_if:nTF { \g_noteobserve } {
    \bool_gset_false:N \g_noteobserve
    Note
  } {
    \bool_gset_true:N \g_noteobserve
    Observe
  }
  \IfValueF{#1}{~}
}

\bool_new:N \g_hencetherefore

\NewDocumentCommand{\hence}{ o }{
  \IfValueT{#1}{
    \str_if_eq:noTF {hence} {#1} {
      \bool_gset_true:N \g_hencetherefore
    } {
      \str_if_eq:noTF {Hence} {#1} {
        \bool_gset_true:N \g_hencetherefore
      } {
        \bool_gset_false:N \g_hencetherefore
      }
    }
  }
  \bool_if:nTF { \g_hencetherefore } {
    \bool_gset_false:N \g_hencetherefore
    hence
  } {
    \bool_gset_true:N \g_hencetherefore
    therefore
  }
  \IfValueF{#1}{~}
}

\NewDocumentCommand{\Hence}{ o }{
  \IfValueT{#1}{
    \str_if_eq:noTF {hence} {#1} {
      \bool_gset_true:N \g_hencetherefore
    } {
      \str_if_eq:noTF {Hence} {#1} {
        \bool_gset_true:N \g_hencetherefore
      } {
        \bool_gset_false:N \g_hencetherefore
      }
    }
  }
  \bool_if:nTF { \g_hencetherefore } {
    \bool_gset_false:N \g_hencetherefore
    Hence
  } {
    \bool_gset_true:N \g_hencetherefore
    Therefore
  }
  \IfValueF{#1}{~}
}

\int_new:N \g_furthermore

\NewDocumentCommand{\Moreover}{ o o }{
  \IfValueT{#1}{
    \str_case:nn {#1} {
	  {Next} {\int_gset:Nn {\g_furthermore} {0}}      
      {Furthermore} {\int_gset:Nn {\g_furthermore} {1}}
      {Moreover} {\int_gset:Nn {\g_furthermore} {2}}
      {In~addition} {\int_gset:Nn {\g_furthermore} {3}}
      {note} {\bool_gset_true:N \g_noteobserve}
      {observe} {\bool_gset_false:N \g_noteobserve}
    }
    \IfValueT{#2}{
      \str_case:nn {#2} {
	    {Next} {\int_gset:Nn {\g_furthermore} {0}}        
        {Furthermore} {\int_gset:Nn {\g_furthermore} {1}}
        {Moreover} {\int_gset:Nn {\g_furthermore} {2}}
        {In~addition} {\int_gset:Nn {\g_furthermore} {3}}
        {note} {\bool_gset_true:N \g_noteobserve}
        {observe} {\bool_gset_false:N \g_noteobserve}
      }
    }
  }
  \int_case:nn { \int_mod:nn {\g_furthermore} {4} } {
	{ 0 } { Next,~\nobs that}    
    { 1 } { Furthermore,~\nobs that}
    { 2 } { Moreover,~\nobs that}
    { 3 } { In~addition,~\nobs that}
  }
  \int_gincr:N \g_furthermore
  \IfValueF{#1}{~}
}

\ExplSyntaxOff

\ExplSyntaxOn

\int_new:N \g_enumit

\NewDocumentCommand{\Enum}{ m o }{
  \IfValueT{#2}{
    \str_case:nn {#2} {
      {ensure} {\int_gset:Nn {\g_enumit} {0}}
      {ensures} {\int_gset:Nn {\g_enumit} {0}}       
      {assure} {\int_gset:Nn {\g_enumit} {1}}
      {assures} {\int_gset:Nn {\g_enumit} {1}}
      {prove} {\int_gset:Nn {\g_enumit} {2}}
      {proves} {\int_gset:Nn {\g_enumit} {2}}
      {show} {\int_gset:Nn {\g_enumit} {3}}
      {shows} {\int_gset:Nn {\g_enumit} {3}}
      {demonstrate} {\int_gset:Nn {\g_enumit} {4}}
      {demonstrates} {\int_gset:Nn {\g_enumit} {4}}
      {yield} {\int_gset:Nn {\g_enumit} {5}}
      {yields} {\int_gset:Nn {\g_enumit} {5}}
      {establish} {\int_gset:Nn {\g_enumit} {6}}
      {establishes} {\int_gset:Nn {\g_enumit} {6}}
      {imply} {\int_gset:Nn {\g_enumit} {7}}
      {implies} {\int_gset:Nn {\g_enumit} {7}}
    }
  }
  \int_case:nn { \int_mod:nn {\g_enumit} {8} } {
    { 0 } { \enum{ #1 }[ensure]{ } }
    { 1 } { \enum{ #1 }[assure] }  
    { 2 } { \enum{ #1 }[prove] }  
    { 3 } { \enum{ #1 }[show] }
    { 4 } { \enum{ #1 }[demonstrate] }
    { 5 } { \enum{ #1 }[yield] }
    { 6 } { \enum{ #1 }[establish] }
    { 7 } { \enum{ #1 }[imply] }
  }
  \int_gincr:N \g_enumit
} 

\ExplSyntaxOff

\ExplSyntaxOn

\seq_const_from_clist:Nn \g_prove_mru {
  establish,
  demonstrate,
  prove,
  show,
  imply,
  ensure
}

\prop_new:N \l__verbsnew
\prop_put:Nnn \l__verbsnew {show} {shows}
\prop_put:Nnn \l__verbsnew {imply} {implies}
\prop_put:Nnn \l__verbsnew {demonstrate} {demonstrates}
\prop_put:Nnn \l__verbsnew {prove} {proves}
\prop_put:Nnn \l__verbsnew {establish} {establishes}
\prop_put:Nnn \l__verbsnew {ensure} {ensures}
\prop_put:Nnn \l__verbsnew {assure} {assures}

\tl_new:N \g_wordtmp
\seq_new:N \l_mytmps

\cs_generate_variant:Nn \str_if_in:nnTF { nVTF }

\NewDocumentCommand{\prove}{ o }{
  \IfValueTF{#1}{
    \seq_clear:N \l_mytmps
    \seq_map_inline:Nn \g_prove_mru {
      \str_if_eq:nnTF {##1} {ensure} {
        \str_set:Nn \l_temps {n}
      } {
        \str_set:Nx \l_temps {\str_head_ignore_spaces:n {##1}}
      }
      \str_if_in:nVTF {#1} \l_temps {
        \seq_put_right:Nn \l_mytmps {##1}
      } { }
    }
    \seq_get_right:NN \l_mytmps \g_wordtmp
  } {
    \seq_get_right:NN \g_prove_mru \g_wordtmp
  }
  \tl_use:N \g_wordtmp
  \IfValueTF{#1}{}{~}
  \seq_gput_left:NV \g_prove_mru \g_wordtmp
  \seq_gremove_duplicates:N \g_prove_mru
}

\NewDocumentCommand{\proves}{ o }{
  \IfValueTF{#1}{
    \seq_clear:N \l_mytmps
    \seq_map_inline:Nn \g_prove_mru {
      \str_if_eq:nnTF {##1} {ensure} {
        \str_set:Nn \l_temps {n}
      } {
        \str_set:Nx \l_temps {\str_head_ignore_spaces:n {##1}}
      }
      \str_if_in:nVTF {#1} \l_temps {
        \seq_put_right:Nn \l_mytmps {##1}
      } { }
    }
    \seq_get_right:NN \l_mytmps \g_wordtmp
  } {
    \seq_get_right:NN \g_prove_mru \g_wordtmp
  }
  \str_set:NV \l_tmpa_str \g_wordtmp
  \prop_get:NVN \l__verbs \l_tmpa_str \l_tmpa_tl
  \tl_use:N \l_tmpa_tl
  \IfValueTF{#1}{}{~}
  \seq_gput_left:NV \g_prove_mru \g_wordtmp
  \seq_gremove_duplicates:N \g_prove_mru
}

\ExplSyntaxOff

\begin{document}

\title{Deep neural networks with ReLU, leaky ReLU, and softplus activation provably overcome the curse of dimensionality for Kolmogorov partial differential equations with Lipschitz nonlinearities in the $L^p$-sense}

\author{
Julia Ackermann$^{1}$,
Arnulf Jentzen$^{2,3}$,
Thomas Kruse$^{4}$,\\ 
\vspace{-0.3cm}
Benno Kuckuck$^5$,
and
Joshua Lee Padgett$^{6,7}$
\bigskip
\\
\small{$^1$ Department of Mathematics \& Informatics,}
\vspace{-0.1cm}\\
\small{University of Wuppertal, Germany, e-mail: \texttt{jackermann@uni-wuppertal.de}}
\smallskip
\\
\small{$^2$ School of Data Science and Shenzhen Research Institute of Big Data,}
\vspace{-0.1cm}\\
\small{The Chinese University of Hong Kong, Shenzhen (CUHK-Shenzhen), \vspace{-0.1cm}\\China, e-mail: \texttt{ajentzen@cuhk.edu.cn}}
\smallskip
\\
\small{$^3$ Applied Mathematics: Institute for Analysis and Numerics,}
\vspace{-0.1cm}\\
\small{University of M{\"u}nster, Germany, e-mail: \texttt{ajentzen@uni-muenster.de}}
\smallskip
\\
\small{$^4$ Department of Mathematics \& Informatics,}
\vspace{-0.1cm}\\
\small{University of Wuppertal, Germany, e-mail: \texttt{tkruse@uni-wuppertal.de}}
\smallskip
\\
\small{$^5$ Applied Mathematics: Institute for Analysis and Numerics,}
\vspace{-0.1cm}\\
\small{University of M{\"u}nster, Germany, e-mail: \texttt{bkuckuck@uni-muenster.de}}
\smallskip
\\
\small{$^6$ Department of Mathematical Sciences, University of Arkansas,}
\vspace{-0.1cm}\\
\small{Arkansas, USA, e-mail: \texttt{padgett@uark.edu}}
\smallskip
\\
\small{$^7$ Center for Astrophysics, Space Physics, and Engineering Research,}
\vspace{-0.1cm}\\
\small{Baylor University, Texas, USA, e-mail: \texttt{padgett@uark.edu}}
}

\date{\today}

\maketitle

\begin{abstract}
In recent years, several deep learning-based methods for the approximation of high-dimensio\-nal partial differential equations (PDEs) have been proposed. 
The considerable interest that these methods have generated in the scientific literature 
is in large part due to numerical simulations which appear to 
demonstrate that such deep learning-based approximation methods seem to have the capacity to overcome 
the curse of dimensionality 
(COD) 
in the numerical approximation of PDEs in the sense that 
the number of computational operations they require to achieve a certain approximation
accuracy $\eps\in(0,\infty)$ grows at most polynomially in the PDE dimension $d\in\N=\{1,2,3,\dots\}$
and the reciprocal of $\eps$.
While there is thus far no mathematical result which proves that one of these methods is indeed
capable of overcoming the 
COD 
in the numerical approximation of PDEs, 
there are now a number of rigorous mathematical results in
the scientific literature 
which show that deep neural networks (DNNs) have the
expressive power to approximate solutions of high-dimensional PDEs without the 
COD 
in the sense that the number of real parameters used to describe the approximating DNNs grows at most 
polynomially in both the PDE dimension $d\in\N$ and the reciprocal $\nicefrac1\eps$
of the prescribed approximation accuracy $\eps\in(0,\infty)$. 
More specifically, [Hutzenthaler,~M., Jentzen,~A., Kruse,~T., and Nguyen,~T.~A.,
{\it SN Part.\ Diff.\ Equ.\ Appl.\ 1}, 2 (2020)] proves that for every $T\in(0,\infty)$,
$a\in\R$, $b\in[a,\infty)$ it holds that solutions 
$u_d\colon [0,T]\times\R^d\to \R$, $d\in\N$, of semilinear heat equations with Lipschitz continuous nonlinearities can be 
approximated by DNNs 
with the rectified linear unit (ReLU) activation at the terminal time 
in the $L^2$-sense on $[a,b]^d$ 
without the 
COD 
provided that the initial value functions $\R^d\ni x\mapsto u_d(0,x)\in \R$, $d\in\N$, 
can be approximated by 
ReLU DNNs 
without the 
COD. 
It is the key contribution of 
this article to generalize this result 
by establishing this 
statement in the $L^p$-sense with $p\in(0,\infty)$ 
and by allowing the activation function to be more general 
covering the ReLU, the leaky ReLU, and the softplus activation functions 
as special cases.  
\end{abstract}

\tableofcontents

%
%
%


\section{Introduction}
\label{sec:intro}

Finding approximate solutions to high-dimensional partial differential equations (PDEs) is one
of the most challenging issues in computational mathematics. 
In recent years, several deep learning-based methods for this approximation problem have been proposed and
have received significant attention in the scientific literature. 
Some of such deep learning-based approximation methods for PDEs are based on classical or strong formulations of PDEs (cf., for example, \cite{hu2023tackling,RaissiEtAl2019PINN,Berg2018AUD,Sirignano2018}), some are based on variational or weak formulations of PDEs (cf., for example, \cite{weinan2018deep,ValsecchiOliva2022xnodewan,ZangBaoYeZhou2020weak,BaoYeZangZhou2020inverseWAN,feng2023fractional}), and some are based on suitable stochastic formulations of the Feynman--Kac type involving the associated forward stochastic differential equations (SDEs) or backward stochastic differential equations (BSDEs), respectively (cf., for example,~\cite{weinan2017deep,Han2018PNAS,HanLong2018,BeckJentzenE2019,FujiiTakahashi2019,henry2017deep,Kolmogorov,ChanMikaelWarin2019,HurePhamWarin2019,PhamWarin2019,GermainPhamWarin2022,JacquierOumgari2019,BeckBeckerCheridito2019,NueskenRichter2020,chassagneux2022deep}).

In particular, we refer, 
	for instance, to 
	\cite{weinan2017deep,Han2018PNAS,HanLong2018} 
	for certain \emph{deep BSDE approximations} for classes of semilinear parabolic PDEs. 
	We refer, for example, to  
	\cite{Sirignano2018} 
	for certain \emph{deep Galerkin approximations} 
	for general classes of PDEs. 
	We refer, for instance, to  
	\cite{BeckJentzenE2019} 
	for certain \emph{deep 2BSDE approximations} for a class of fully nonlinear 
	parabolic second-order PDEs. 
	We refer, for example, to 
	\cite{weinan2018deep} 
	for certain \emph{deep Ritz approximations} for a class of elliptic PDEs. 
	We refer, for instance, to 
	\cite{FujiiTakahashi2019} 
	for certain \emph{deep BSDE approximations involving asymptotic expansions}
	for a class of nonlinear parabolic PDEs
	and extensions of such approximation techniques to reflected BSDEs. 
	We refer, for example, to  
	\cite{henry2017deep} 
	for certain \emph{deep primal-dual approximations} for a class of nonlinear parabolic PDEs and 
	applications of such approximation methods 
	to the pricing of counterparty risks and to the computation of initial margins.  
	We refer, for instance, to  
	\cite{Berg2018AUD} 
	for certain \emph{deep artificial neural network (ANN) approximations 
		using collocation techniques} 
	for advection and diffusion type PDEs in complex geometries.
	We refer, for example, to  
	\cite{hu2023tackling,RaissiEtAl2019PINN} 
	for certain \emph{physics-informed neural network (PINN) approximations}
	for general classes of PDEs. 
	We refer, for instance, to  
	\cite{Kolmogorov} 
	for certain \emph{deep ANN approximations based on discretizations of SDEs} for a class of linear Kolmogorov PDEs on an entire region. 
	We refer, for example, to 
	\cite{ChanMikaelWarin2019} 
	for certain \emph{Monte Carlo based deep ANN approximations} 
	for a class of semilinear Kolmogorov PDEs. 
	We refer, for instance, to 
	\cite{HurePhamWarin2019,PhamWarin2019,GermainPhamWarin2022} 
	for certain \emph{deep backward dynamic programming approximations} for classes of nonlinear parabolic PDEs.
	We refer, for example, to  
	\cite{JacquierOumgari2019} 
	for certain \emph{deep BSDE based approximations}
	for a class of path-dependent PDEs arising in affine rough volatility models. 
	We refer, for instance, to  
	\cite{BeckBeckerCheridito2019} 
	for certain \emph{deep splitting approximations} for a class of nonlinear parabolic PDEs. 
	We refer, for example, to  
	\cite{NueskenRichter2020} 
	for certain \emph{iterative diffusion optimization approximations} for 
	a class of 
	Hamilton--Jacobi--Bellman PDEs. 
	We refer, for instance, to  
	\cite{chassagneux2022deep} 
	for certain \emph{deep Runge--Kutta approximations} for 
	a class of semilinear parabolic PDEs. 
For more extensive overviews 
on such and related 
deep learning-based methods for high-dimensional PDEs, 
we refer, for example, to the survey articles
\cite{BeckHutzenthalerJentzenKuckuck2023anOverview,EHanJentzen2021algorithms,BlechschmidtErnst2021,germain2021neural}. 

The considerable interest in deep learning-based approximation methods
for high-dimen\-sional PDEs is in large part due to numerical simulations which appear to 
demonstrate that some of these deep learning-based approximation methods might have the capacity to overcome 
the 
curse of dimensionality (COD) 
(cf., e.g., Bellman \cite{bellman2013dynamic} and Novak \& Wo{\'z}niakowski \cite[Chapter 1]{Novak2008}) 
in the numerical approximation of PDEs in the sense that 
the number of computational operations they require to achieve a certain approximation
accuracy $\eps\in(0,\infty)$ grows at most polynomially in the PDE dimension $d\in\N=\{1,2,3,\dots\}$
and the reciprocal of $\eps$. 
In the last 
few  
years, a number of rigorous mathematical results have appeared in
the scientific literature which show that deep 
ANNs 
have the
expressive power to approximate solutions of high-dimensional PDEs without the COD 
in the sense that the number of real parameters used to describe the approximating deep ANNs grows at most 
polynomially in both the PDE dimension $d\in\N$ and the reciprocal $\nicefrac1\eps$ 
of the prescribed approximation accuracy $\eps\in(0,\infty)$; cf., e.g., 
\cite{BernerGrohsJentzen2018,
ElbraechterSchwab2018,
GrohsWurstemberger2018,
HutzenthalerJentzenKruse2019,
JentzenSalimovaWelti2018,
KutyniokPeterseb2019,
ReisingerZhang2019,
GrohsJentzenSalimova2019,
GononGrohsEtAl2019,
HornungJentzenSalimova2020,
GrohsHerrmann2020,
GononSchwab2020,
cioicalicht2022deep,
BaggenstosSalimova2023,
grohs2021deep}.

While the articles 
\cite{BernerGrohsJentzen2018,
ElbraechterSchwab2018,
GrohsWurstemberger2018,
JentzenSalimovaWelti2018,
KutyniokPeterseb2019,
ReisingerZhang2019,
GrohsJentzenSalimova2019,
GononGrohsEtAl2019,
HornungJentzenSalimova2020,
GrohsHerrmann2020,
BaggenstosSalimova2023,
GononSchwab2020} 
prove such deep ANN approximation results for 
linear PDEs,  
the articles 
\cite{HutzenthalerJentzenKruse2019,cioicalicht2022deep,grohs2021deep}  establish deep ANN approximation results for certain nonlinear PDEs.
In the article Hutzenthaler et al.\ \cite{HutzenthalerJentzenKruse2019} it is shown that deep ANNs with the rectified linear unit (ReLU) activation function 
$\R\ni x \mapsto \max\{x,0\} \in \R$  
can approximate solutions of semilinear heat PDEs at the terminal time in the $L^2$-sense without the COD provided that the initial value functions can be approximated by deep ReLU ANNs without the COD. 
The article Cioica-Licht et al.\ \cite{cioicalicht2022deep} 
extends the findings in Hutzenthaler et al.\ \cite{HutzenthalerJentzenKruse2019} in several ways. 
Specifically, in 
Cioica-Licht et al.\ \cite{cioicalicht2022deep} 
it is shown  
that for every $T\in(0,\infty)$ solutions 
$u_d\colon [0,T]\times\R^d\to \R$, $d\in\N$, of certain semilinear Kolmogorov PDEs with Lipschitz continuous nonlinearities can be 
approximated by deep ANNs with ReLU activation 
at the terminal time 
in the $L^2$-sense 
without the COD 
provided that the initial value functions $\R^d\ni x\mapsto u_d(0,x)\in \R$, $d\in\N$, and the coefficients of the PDEs 
can be approximated by ANNs with ReLU activation without the COD. 
For $p>2$ or for the leaky ReLU or the softplus activation, 
up to our best knowledge, there is no result in the literature that shows 
that deep ANNs can overcome the COD 
in the $L^p$-approximation of nonlinear PDEs.

It is the key contribution of the present article to show that for every $p\in(0,\infty)$ we have that solutions of semilinear heat PDEs with Lipschitz continuous nonlinearities can be approximated in the $L^p$-sense 
by deep ANNs with ReLU, leaky ReLU, or softplus activation without the COD. 
More precisely, we prove that 
for any of these types of activation functions and 
for every $T\in(0,\infty)$, 
$a\in\R$, $b\in(a,\infty)$ solutions 
$u_d\colon [0,T]\times\R^d\to \R$, $d\in\N$, of semilinear heat equations with Lipschitz continuous nonlinearities can be 
approximated 
by deep ANNs in the $L^p$-sense with $p\in(0,\infty)$ on $[a,b]^d$ 
at the terminal time 
without the COD 
provided that the initial value functions $\R^d\ni x\mapsto u_d(0,x)\in \R$, $d\in\N$, 
can be approximated by ANNs without the COD (see \cref{cor:final} below for details). 
This extends the result in Hutzenthaler et al.\ \cite{HutzenthalerJentzenKruse2019} from $L^2$-approximation to $L^p$-approximation with $p\in(0,\infty)$ and from ReLU activation to ReLU, leaky ReLU, or softplus activation functions.

In order to illustrate the contribution of this article in more detail, we now present in the following result, 
\cref{thm:intro} below, a special case of \cref{theorem:final} in \cref{sec:4_1}, which is the main result of this paper. 
Below \cref{thm:intro} we add several explanatory sentences in which we aim to describe the used mathematical objects and the statement of \cref{thm:intro} in words.

\begin{samepage}
\begin{theorem}\label{thm:intro}
Let 
$T,\constantAssumpMainThm,p \in (0,\infty)$,
let $f\colon\R \to \R$ 
be Lipschitz continuous,  
for every $d\in\N$ let 
$u_d \in C^{1,2}([0,T]\times \R^d,\R)$ 
satisfy for all $t \in [0,T]$, $x = (x_1,\dots,x_d) \in \R^d$ 
that 
\begin{equation}\label{eq:1}
(\tfrac{\partial}{\partial t}u_d)(t,x) 
= 
(\Delta_x u_d)(t,x) 
+ 
f(u_d(t,x)), 
\end{equation}
let $\nu\in\{0,1\}$, 
$\alpha \in [0,\infty)\backslash\{1\}$,
$\fa_0,\fa_1\in C(\R,\R)$ 
satisfy for all 
$x \in \R$ 
that
$\fa_0(x) = \max\{x,\alpha x\}$ 
and 
$\fa_1(x) = \ln(1+\exp(x))$,
for every 
$d\in\N$, 
$x=(x_1,\allowbreak \dots, \allowbreak x_d)\in\R^d$ 
let $\mathbf A(x) \in \R^d$ satisfy 
$\mathbf A(x) = (\fa_{\nu}(x_1),\allowbreak \dots,\allowbreak\fa_{\nu}(x_d))$, 
let 
\begin{equation}
\textstyle{
\ANNs = \medcup{L\in\N}\medcup{l_0,l_1,\dots,l_L\in\N}(\bigtimes_{k=1}^{L} (\R^{l_k\times l_{k-1}}\times \R^{l_k}))
},
\end{equation}
for every  
$L\in\N$, 
$l_0,l_1,\dots,l_L\in\N$, 
$\Phi = ((W_1,B_1),\allowbreak\dots,\allowbreak(W_L,B_L)) \in (\bigtimes_{k=1}^L\allowbreak (\R^{l_k\times l_{k-1}}\times \R^{l_k}))$  
let $\mathcal{R}(\Phi) \colon \R^{l_0} \to \R^{l_L}$ 
and $\mathcal{P}(\Phi) \in \N$ 
satisfy for all 
$v_0\in \R^{l_0}$, $v_1\in\R^{l_1}$, $\dots$, $v_L\in\R^{l_L}$
with 
$\forall\, k\in \{1,2,\dots,L\}\colon \allowbreak v_k = \mathbf A(W_k v_{k-1}+B_k)$ 
that
\begin{equation}
\textstyle{
(\mathcal{R}(\Phi))(v_0) = W_L v_{L-1}+B_L
\qquad
\text{and}
\qquad
\mathcal{P}(\Phi) = \sum_{k=1}^L l_k(l_{k-1}\allowbreak+1) }
\dc
\end{equation}
and assume for all $d\in\N$, $\varepsilon\in (0,1]$ 
that there exists $\ANNassumpMainThm \in \ANNs$ such that for all 
$t\in[0,T]$, $x = (x_1,\dots,x_d) \in \R^d$ it holds that 
$\mathcal{R}(\ANNassumpMainThm) \in C(\R^d,\R)$, 
$\mathcal{P}(\ANNassumpMainThm) \le \constantAssumpMainThm d^{\constantAssumpMainThm} \varepsilon^{-\constantAssumpMainThm}$, 
and 
\begin{equation}\label{eq:thm_intro_ass_u}
\textstyle{\varepsilon \abs{u_d(t,x)} + \allowbreak \abs{u_d(0,x) \allowbreak - (\mathcal{R}(\ANNassumpMainThm))(x)} 
	\le \varepsilon \constantAssumpMainThm d^{\constantAssumpMainThm} (1 + \sum_{k=1}^d\abs{x_k}^{ \constantAssumpMainThm} )  }
.
\end{equation}
Then there exists 
$c\in \R$ 
such that for all 
$d\in\N$, 
$\varepsilon\in(0,1]$ 
there exists $\ANNresultMainThm\in\ANNs$
such that 
$\mathcal{R}(\ANNresultMainThm) \in C(\R^{d},\R)$, 
$\mathcal{P}(\ANNresultMainThm) \leq c d^c \varepsilon^{-c}$,
and
\begin{equation}\label{eq:2}
\textstyle{
\bigl[\int_{[0,1]^d} \abs{u_d(T,x) 
	- 
	(\mathcal{R}(\ANNresultMainThm))(x)}^p \dx x \bigr]^{\nicefrac{1}{p}} 
\leq
\varepsilon}
\dpp
\end{equation}
\end{theorem}
\end{samepage}

\Cref{thm:intro} is an immediate consequence of \cref{cor:final} in \cref{sec:4_2} below.
\Cref{cor:final}, in turn, follows from \cref{theorem:final}, which is the main result of this article
(see \cref{sec:4} for details).
In the following we provide some explanatory comments concerning the mathematical objects appearing in \cref{thm:intro}. 

The function $\fa_{\nu} \colon \R \to \R$ in \cref{thm:intro} serves as the activation function which we employ in the approximating ANNs in \cref{thm:intro} 
and the function 
\begin{equation}
( \medcup{d\in\N} \R^d ) \ni x \mapsto \A(x) \in ( \medcup{d\in\N} \R^d )
\end{equation}
represents a suitable multidimensional version of the activation function $\fa_{\nu} \colon \R \to \R$ in \cref{thm:intro}. 
The function $\fa_{\nu} \colon \R \to \R$ may be the ReLU activation (corresponding to the case $\nu = \alpha = 0$ in \cref{thm:intro}), the leaky ReLU activation (corresponding to the case $\nu = 0$ and $\alpha \in (0,1)$ in \cref{thm:intro}), or the softplus activation (corresponding to the case $\nu=1$ in \cref{thm:intro}).

The set 
$\mathbf N=\medcup{L\in\N} \medcup{l_0,l_1,\dots,l_L\in\N}(\medtimes_{k=1}^L (\R^{l_k\times l_{k-1}}\times \R^{l_k}))$ 
in \cref{thm:intro} 
represents the set of all ANNs which we employ 
to approximate the solutions of the PDEs under consideration. 
Observe 
that for every 
ANN 
$\Phi\in\mathbf N$ 
we have that 
\begin{equation}
\mathcal R(\Phi)\in\medcup{k,l\in\N} C(\R^k,\R^l)
\end{equation}
is the realization function of the 
ANN 
$\Phi$ with the activation function
$\fa_{\nu} \colon \R \to~\R$. 
Moreover,  
we note that for every  
ANN 
$\Phi\in\ANNs$ we have that $\mathcal P(\Phi)\in\N$ is the 
number of real numbers used to describe the ANN~$\Phi$.  
Very roughly speaking, $\mathcal P(\Phi)$ corresponds to the amount of memory that is needed 
on a computer to store the 
ANN 
$\Phi\in\ANNs$.

The real number $T\in(0,\infty)$ 
in \cref{thm:intro} 
specifies the time horizon of the PDEs (see~\cref{eq:1} above)
whose solutions we intend to approximate by deep ANNs in \cref{eq:2}
in \cref{thm:intro}.
The real number $\constantAssumpMainThm\in(0,\infty)$ 
in \cref{thm:intro}
is a constant which we employ to formulate our 
regularity and approximation hypotheses 
in \cref{thm:intro}. 
The real number $p\in(0,\infty)$ 
in \cref{thm:intro} 
is used to specify the way we measure the error between the exact solutions 
of the PDEs under consideration and their deep ANN approximations, that is, we measure the error
between the exact solutions of the PDEs under consideration and their deep ANN approximations in the $L^p$-sense
(see \cref{eq:2} for details).

In \cref{thm:intro} we assume that the initial conditions of the PDEs (see \cref{eq:1}) whose solutions we intend to 
approximate by deep ANNs without the COD  
can be approximated by ANNs without the COD 
(see~\cref{eq:thm_intro_ass_u} above).
The function $f\colon \R\to\R$ 
in \cref{thm:intro} 
specifies the nonlinearity in the PDEs (see~\cref{eq:1})
whose solutions we intend to approximate by deep ANNs in
\cref{thm:intro}. 
The functions 
\begin{equation}
u_d\colon [0,T]\times\R^d\to\R, \quad d\in\N,
\end{equation}
in \cref{thm:intro} 
describe the exact solutions of the PDEs in \cref{eq:1}. 
Observe that the hypothesis 
in 
\cref{eq:thm_intro_ass_u} 
in \cref{thm:intro} 
also ensures that for all
$d\in\N$, $\varepsilon\in(0,1]$, $t\in[0,T]$, $x\in\R^d$ 
we have that 
$\varepsilon\abs{u_d(t,x)} \le \varepsilon \constantAssumpMainThm d^{\constantAssumpMainThm} (1+\sum_{k=1}^d \abs{x_k}^{\constantAssumpMainThm})$. 
This, in turn, assures 
that for all 
$d\in\N$,  
$t\in[0,T]$, 
$x\in\R^d$  
we have that 
\begin{equation}\label{eq:introasspolynomialgrowth}
\textstyle{\abs{u_d(t,x)}\leq \constantAssumpMainThm d^{\constantAssumpMainThm} (1+\sum_{k=1}^d \abs{x_k}^{\constantAssumpMainThm} ).}
\end{equation}
Note that~\cref{eq:introasspolynomialgrowth} ensures that for all $d\in\N$ we have that  
the solution 
$u_d\colon[0,T]\times\R^d\to\R$ 
of~\cref{eq:1} grows at most polynomially.
These polynomial growth properties of the solutions assure that 
the solutions of~\cref{eq:1} with 
the fixed initial value functions 
$\R^d\ni x\mapsto u_d(0,x)\in\R$, $d\in\N$,
are unique.

\Cref{thm:intro} establishes that 
for every $d \in \N$, $\varepsilon \in (0,1]$ 
there exists 
an ANN 
$\ANNresultMainThm_{d,\varepsilon}\in\ANNs$  
such that 
the $L^p$-distance 
with respect to the Lebesgue measure
on $[0,1]^d$ 
between the exact solution $u_d\colon[0,T]\times\R^d\to\R$ 
at time~$T$ of the PDE in \cref{eq:1} and the realization 
\begin{equation}
\mathcal R(\ANNresultMainThm_{d,\eps})\colon \R^{d}\to\R
\end{equation}
of the 
ANN 
$\ANNresultMainThm_{d,\eps}$ 
is bounded by $\eps$ 
and such that the number of parameters 
of the 
ANN 
$\ANNresultMainThm_{d,\eps}\in\ANNs$ 
grows at most polynomially in both the PDE dimension $d$ and 
the reciprocal $\nicefrac 1\eps$ of the prescribed approximation accuracy $\eps$.

Although \Cref{thm:intro} is restricted to measuring the $L^p$-distance with respect to the Lebesgue measure on 
$[0,1]^d$, 
our more general deep ANN approximation results in \cref{sec:4} (see \cref{theorem:final}, \cref{cor:almostfinal}, and \cref{cor:final} in \cref{sec:4})
allow measuring the $L^p$-distance with respect to more general probability measures on $\R^d$.
In particular, for all 
$a\in\R$, 
$b\in(a,\infty)$
we have that the more general deep ANN approximation results in \cref{sec:4} 
enable 
measuring the $L^p$-distance with respect to the uniform distribution on $[a,b]^d$.
Furthermore, we note that \cref{theorem:final} in \cref{sec:4} is formulated for general activation functions provided that 
with the considered general activation function 
there exists 
a shallow ANN representation for the identity function $\R \ni x \mapsto x \in \R$ 
on the real numbers 
and provided that 
with the considered general activation function 
the Lipschitz continuous nonlinearity can be approximated 
with appropriate convergence rates 
by ANNs. 

Our proof of \cref{thm:intro} is strongly based on employing suitable nonlinear Monte Carlo approximations, so-called 
\emph{full history recursive multilevel Picard} (MLP) approximations, which provably approximate PDEs of the form \cref{eq:1} without the COD.  
MLP methods have been introduced 
in \cite{EHutzenthaler2021,HutzenthalerJentzenKruse2018,EHutzenthaler2019} 
and have been extended to more general situations in, e.g., 
\cite{BeckHornungEtAl2019,HutzenthalerPricing2019,GilesJentzenWelti2019,HutzenthalerKruse2020,HutzenthalerKruseNguyen2022,HutzenthalerNguyen2022,neufeld2022multilevel,BeckGononJentzen2020,BeckerBraunwarthEtAl2020,HutzenthalerJentzenKruse2022,HutzenthalerJentzenKruseNguyen2020}. 
In particular, the $L^p$-error with $p\in(0,\infty)$ for MLP approximations of certain semilinear PDEs has been analyzed in 
Hutzenthaler et al.~\cite{PadgettJentzen2021}
and  
in the present work we employ these findings to establish the desired ANN approximation results. 
For further references on MLP methods we refer, for example, to the survey articles~\cite{BeckHutzenthalerJentzenKuckuck2023anOverview,EHanJentzen2021algorithms}.

The remainder of this article is organized as follows:
In \cref{subsec:structured_description} we review the necessary
basic preparatory material on ANNs that we need in the later sections of this work. 
In \cref{subsubsec:ANN_repres_for_MLP} 
we study deep ANN representations for MLP approximations for PDEs of the form \cref{eq:1}. 
The ANN representations for MLP approximations from \cref{subsubsec:ANN_repres_for_MLP} are then used in \cref{sec:4} to prove \cref{thm:intro} and its above sketched generalizations.

%
%
%


%
%
%

\section{Artificial neural network (ANN) calculus}
\label{subsec:structured_description}

In this section we review the necessary basic preparatory material on ANNs that we need in the later sections of this work.
The conceptualities and the lemmas provided in this section are 
elementary or well-known 
in the literature and the specific material in this section 
mostly consists of slightly modified extracts from 
Grohs et al.\ \cite[Section 2]{GrohsHornungJentzen2019} and Grohs et al.\ \cite[Section 3]{GrohsJentzenSalimova2019}.

\subsection{Structured description of ANNs}
\label{subsubsec:structured_description_of_DNNs}

\begin{definition}[ANNs]\label{def:ANN}
We denote by $\ANNs$ the set given by 
\begin{equation}\label{eq:defANN}
\textstyle
\ANNs
\textstyle
=
\bigcup_{L \in \N}
\bigcup_{ l_0,l_1,\dots, l_L \in \N }
\pr[\big]{
\bigtimes_{k = 1}^L (\R^{l_k \times l_{k-1}} \times \R^{l_k})
}
\end{equation}
and we denote by 	
$
\paramANN \colon \ANNs \to \N
$,
$
\lengthANN \colon \ANNs \to \N
$,
$
\inDimANN \colon \ANNs \to \N
$,
$
\outDimANN \colon \ANNs \to \N
$, 
$
\hiddenLength \colon \ANNs \to 
\N_0 = 
\N\cup\{0\}
$, 
$
\dims \colon \ANNs \to \medcup{ L \in \N }\N^{L}
$,
and
$\dimANNlevel_n \colon \ANNs \to \N_0$, $n \in \N_0$,
the functions which satisfy
for all 
$L \in\N$, 
$l_0,l_1,\dots, l_L \in \N$, 
$
\Phi 
\in  \allowbreak
( \bigtimes_{k = 1}^L\allowbreak(\R^{l_k \times l_{k-1}} \times \R^{l_k}) )$,
$n \in \N_0$
that
$\paramANN(\Phi)
=
\sum_{k = 1}^L l_k(l_{k-1} + 1) 
$, 
$\lengthANN(\Phi)=L$,  
$\inDimANN(\Phi)=l_0$,  
$\outDimANN(\Phi)=l_L$, 
$\hiddenLength(\Phi)=L-1$, 
$\dims(\Phi)= (l_0,l_1,\dots, l_L)$, 
and 
\begin{align}\label{dim_operator}
\begin{split}
\dimANNlevel_n (\Phi) =
\begin{cases}
l_n &\colon n \leq L \\
0 &\colon n > L
\end{cases}
\end{split}
\dpp
\end{align}
\end{definition}

\cfclear 
\begin{definition}[ANN]\label{def:ANN2}
\cfconsiderloaded{def:ANN2}
\cfload
We say that $\Phi$ is an ANN 
if and only if it holds that $\Phi \in \ANNs$ \cfout.
\end{definition}

\begin{definition}[Standard and maximum norms]\label{DNN_norm}
We denote by 
$\norm{\cdot} \colon \medcup{d\in\N} \R^d \to \R$
and
$\normmm{\cdot} \colon \medcup{d\in\N} \R^d \to \R$ the functions which satisfy for all 
$d\in\N$, 
$x = (x_1,\dots,x_d)\in\R^d$ 
that 
$\norm{x} = [\sum_{i=1}^d \abs{x_i}^2]^{\nicefrac{1}{2}}$
and
$\normmm{x} = \max_{i\in\{1,2,\dots,d\}} \abs{x_i}$.
%
\end{definition}

\cfclear
\begin{athm}{lemma}{lem:dims_and_params}
	It holds for all $\Phi \in \ANNs$ that 
	$\lengthANN(\Phi) + \normmm{\dims(\Phi)} \le \paramANN(\Phi)$
	\cfout.
\end{athm}

\begin{aproof}
	Observe that it holds for all 
	$L \in\N$, 
	$l_0,l_1,\dots, l_L \in \N$  
	that 
	\begin{equation}
	\sum_{k=1}^L l_k l_{k-1} \ge \max\left\{ \sum_{k=1}^L l_k, \sum_{k=1}^L l_{k-1} \right\} \ge \max\{ l_0,l_1,\ldots,l_L \} .
	\end{equation}
	Hence, we obtain for all 	
	$L \in\N$, 
	$l_0,l_1,\dots, l_L \in \N$,   
	$\Phi \in  \allowbreak ( \bigtimes_{k = 1}^L\allowbreak(\R^{l_k \times l_{k-1}} \times \R^{l_k}) )$ 
	that 
	\begin{equation}
	\paramANN(\Phi) 
	= \sum_{k=1}^L l_k(l_{k-1}+1) 
	= \sum_{k=1}^L l_k + \sum_{k=1}^L l_k l_{k-1}
	\ge L +  \max\{ l_0,l_1,\ldots,l_L \} 
	= \lengthANN(\Phi) + \normmm{\dims(\Phi)} .
	\end{equation}
\end{aproof}

\begin{definition}[Multidimensional version]\label{multi}
Let $d\in\N$, $a\in C(\R,\R)$. 
Then we denote by $\mathfrak{M}_{a,d}\colon\R^d\to\R^d$ the function which satisfies for all $x=(x_1,\dots,x_d)\in\R^d$ that 
\begin{equation}\label{multidim_version:Equation}
\fM_{a,d}(x) 
= 
\pr[\big]{ a(x_1),\dots,a(x_d) }
\dpp
\end{equation}
\end{definition}

\cfclear
\begin{definition}[Realization associated to an ANN]
\label{def:ANNrealization}
\cfconsiderloaded{def:ANNrealization}
Let $a\in C(\R,\R)$ \cfload. 
Then we denote by 
$\functionANN{a} \colon \ANNs \to \pr[]{ \medcup{k,l\in\N} C(\R^k,\R^l) }$
the function which satisfies for all 
$ L\in\N$, 
$l_0,l_1,\ldots, l_L \in \N$, 
$\Phi = ((W_1, B_1),(W_2,\allowbreak B_2),\allowbreak \ldots, (W_L,\allowbreak B_L)) \in \allowbreak \bigl( \bigtimes_{k = 1}^L\allowbreak(\R^{l_k \times l_{k-1}} \times \R^{l_k})\bigr)$,
$x_0 \in \R^{l_0}$, $x_1 \in \R^{l_1}$, $\dots$, $x_{L-1} \in \R^{l_{L-1}}$
with 
$\forall \, k \in \{1,2,\ldots,L-1\} \colon x_k = \activationDim{l_k}(W_k x_{k-1} + B_k)$ 
that
\begin{equation}
\label{setting_NN:ass2}
\functionANN{a}(\Phi) \in C(\R^{l_0},\R^{l_L})
\qandq
( \functionANN{a}(\Phi) ) (x_0) = W_L x_{L-1} + B_L \ifnocf.
\end{equation}
\cfout[.]
\end{definition}

\subsection{Compositions of ANNs}
\label{subsubsec:compositions_of_dnns}

\cfclear
\begin{definition}[Composition of ANNs]
\label{def:ANNcomposition}
\cfconsiderloaded{def:ANNcomposition}
\cfload
We denote by 
$\compANN{(\cdot)}{(\cdot)}\colon\allowbreak \{(\Phi_1,\Phi_2)\allowbreak\in\ANNs\times \ANNs\colon \inDimANN(\Phi_1)=\outDimANN(\Phi_2)\}\allowbreak\to\ANNs$
the function which satisfies for all 
$ L,\mathfrak{L}\in\N$, 
$l_0,l_1,\ldots, l_L, \mathfrak{l}_0,\mathfrak{l}_1,\ldots, \mathfrak{l}_\mathfrak{L} \in \N$, 
$
\Phi_1
=
((W_1, B_1),(W_2, B_2),\allowbreak \ldots, (W_L,\allowbreak B_L))
\in  \allowbreak
\bigl( \bigtimes_{k = 1}^L\allowbreak(\R^{l_k \times l_{k-1}} \times \R^{l_k})\bigr)
$,
$
\Phi_2
=
((\mathscr{W}_1, \mathscr{B}_1),\allowbreak(\mathscr{W}_2, \mathscr{B}_2),\allowbreak \ldots, \allowbreak(\mathscr{W}_\mathfrak{L},\allowbreak \mathscr{B}_\mathfrak{L}))
\in  \allowbreak
\bigl( \bigtimes_{k = 1}^\mathfrak{L}\allowbreak(\R^{\mathfrak{l}_k \times \mathfrak{l}_{k-1}} \times \R^{\mathfrak{l}_k})\bigr)
$ 
with 
$l_0=\inDimANN(\Phi_1)=\outDimANN(\Phi_2)=\mathfrak{l}_{\mathfrak{L}}$
that
\begin{equation}\label{ANNoperations:Composition}
\begin{split}
&\compANN{\Phi_1}{\Phi_2}
=
\\
&
\begin{cases} 
\!\!\lrSpace
\begin{array}{r}
\big((\mathscr{W}_1, \mathscr{B}_1),(\mathscr{W}_2, \mathscr{B}_2),\ldots, (\mathscr{W}_{\mathfrak{L}-1},\allowbreak \mathscr{B}_{\mathfrak{L}-1}),
(W_1 \mathscr{W}_{\mathfrak{L}}, W_1 \mathscr{B}_{\mathfrak{L}}+B_{1}),\\ (W_2, B_2), (W_3, B_3),\ldots,(W_{L},\allowbreak B_{L})\big)
\end{array}
& 
\colon  
L>1<\mathfrak{L} 
\\[3ex]
\big( (W_1 \mathscr{W}_{1}, W_1 \mathscr{B}_1+B_{1}), (W_2, B_2), (W_3, B_3),\ldots,(W_{L},\allowbreak B_{L}) \big)
&
\colon 
L>1=\mathfrak{L}
\\[1ex]
\big((\mathscr{W}_1, \mathscr{B}_1),(\mathscr{W}_2, \mathscr{B}_2),\allowbreak \ldots, (\mathscr{W}_{\mathfrak{L}-1},\allowbreak \mathscr{B}_{\mathfrak{L}-1}),(W_1 \mathscr{W}_{\mathfrak{L}}, W_1 \mathscr{B}_{\mathfrak{L}}+B_{1}) \big)
&
\colon 
L=1<\mathfrak{L}  
\\[1ex]
\bigl((W_1 \mathscr{W}_{1}, W_1 \mathscr{B}_1+B_{1})\bigr)
&
\colon 
L=1=\mathfrak{L} 
\end{cases}
\end{split}
\ifnocf.
\end{equation}
\cfout[.]
\end{definition}

\subsection{Powers and extensions of ANNs}

\begin{definition}[Identity matrix]\label{def:identityMatrix}
Let $d\in\N$. 
Then we denote by $\idMatrix_{d}\in \R^{d\times d}$ the identity matrix in $\R^{d\times d}$.
\end{definition}

\cfclear
\begin{definition}[Powers of ANNs]\label{def:iteratedANNcomposition}
\cfconsiderloaded{def:iteratedANNcomposition}
\cfload
We denote by 
$\power{(\cdot)}{n} \colon \{\Phi\in \ANNs\colon \inDimANN(\Phi)=\outDimANN(\Phi)\}\allowbreak\to\ANNs$, $n\in\N_0$, 
the functions
which satisfy for all 
$n\in\N_0$, 
$\Phi\in\ANNs$ 
with 
$\inDimANN(\Phi)=\outDimANN(\Phi)$ 
that 
\begin{equation}\label{iteratedANNcomposition:equation}
\begin{split}
\power{\Phi}{n} 
=
\begin{cases} \big(\idMatrix_{\outDimANN(\Phi)},(0,0,\dots, 0)\big)\in\R^{\outDimANN(\Phi)\times \outDimANN(\Phi)}\times \R^{\outDimANN(\Phi)}
&
\colon 
n=0 
\\
\compANN{\Phi}{(\power{\Phi}{(n-1)})} 
&
\colon 
n\in\N
\end{cases}
\end{split}
\end{equation}	
\cfout.
\end{definition}

\cfclear
\begin{definition}[Extension of ANNs]\label{def:ANNenlargement}
\cfconsiderloaded{def:ANNenlargement}
\cfload
Let 
$L\in\N$, 
$\Psi\in \ANNs$ 
satisfy that 
$\inDimANN(\Psi)=\outDimANN(\Psi)$.
Then
we denote by 
$\longerANN{L,\Psi}\colon \{\Phi\in\ANNs\colon (\lengthANN(\Phi)\le L \andShort \outDimANN(\Phi)=\inDimANN(\Psi)) \}\to \ANNs$ the function which satisfies for all 
$\Phi\in\ANNs$ 
with 
$\lengthANN(\Phi)\le L$ 
and 
$\outDimANN(\Phi)=\inDimANN(\Psi)$ 
that
\begin{equation}\label{ANNenlargement:Equation}
\longerANN{L,\Psi}(\Phi)
=	 
\compANN{(\power{\Psi}{(L-\lengthANN(\Phi))})}{\Phi}
\end{equation}
\cfout.
\end{definition}

\subsection{Parallelizations of ANNs}
\label{subsubsec:parallelizations_of_dnns}

\cfclear
\begin{definition}[Parallelization of ANNs]
\label{def:simpleParallelization}
\cfconsiderloaded{def:simpleParallelization}
Let $n\in\N$. 
Then we denote by 
\begin{equation}
\parallelizationSpecial_{n}\colon \big\{(\Phi_1,\Phi_2,\dots, \Phi_n)\in\ANNs^n\colon \lengthANN(\Phi_1)= \lengthANN(\Phi_2)=\ldots =\lengthANN(\Phi_n) \big\}\to \ANNs
\end{equation}
the function which satisfies for all 
$L\in\N$,
$(l_{1,0},l_{1,1},\dots, l_{1,L}), (l_{2,0},l_{2,1},\dots, l_{2,L}),\dots,\allowbreak (l_{n,0},\allowbreak l_{n,1},\allowbreak\dots, l_{n,L})\in\N^{L+1}$, 
$\Phi_1=((W_{1,1}, B_{1,1}),(W_{1,2}, B_{1,2}),\allowbreak \ldots, (W_{1,L},\allowbreak B_{1,L}))\in \bigl( \bigtimes_{k = 1}^L\allowbreak(\R^{l_{1,k} \times l_{1,k-1}} \times \R^{l_{1,k}})\bigr) , \allowbreak \Phi_2\allowbreak=\allowbreak((W_{2,1},\allowbreak B_{2,1}),\allowbreak(W_{2,2}, B_{2,2}),\allowbreak \ldots, (W_{2,L},\allowbreak B_{2,L}))\in \bigl( \bigtimes_{k = 1}^L\allowbreak(\R^{l_{2,k} \times l_{2,k-1}} \times \R^{l_{2,k}})\bigr) , \allowbreak \ldots , \allowbreak \Phi_n= \allowbreak ((W_{n,1},\allowbreak B_{n,1}), \allowbreak (W_{n,2}, \allowbreak B_{n,2}),\allowbreak \ldots, \allowbreak(W_{n,L},\allowbreak B_{n,L}))\in \bigl( \bigtimes_{k = 1}^L\allowbreak(\R^{l_{n,k} \times l_{n,k-1}} \times \R^{l_{n,k}})\bigr)$
that
\begin{align}\label{parallelisationSameLengthDef}
\parallelizationSpecial_{n}(\Phi_1,\Phi_2,\dots,\Phi_n) =& 
\left(
\pa{\begin{pmatrix}
W_{1,1}& 0& 0& \cdots& 0\\
0& W_{2,1}& 0&\cdots& 0\\
0& 0& W_{3,1}&\cdots& 0\\
\vdots& \vdots&\vdots& \ddots& \vdots\\
0& 0& 0&\cdots& W_{n,1}
\end{pmatrix} 
,
\begin{pmatrix}
B_{1,1}\\B_{2,1}\\B_{3,1}\\\vdots\\ B_{n,1}
\end{pmatrix}}
,
\right.
\nonumber
\\
&
\quad
\pa{\begin{pmatrix}
W_{1,2}& 0& 0& \cdots& 0\\
0& W_{2,2}& 0&\cdots& 0\\
0& 0& W_{3,2}&\cdots& 0\\
\vdots& \vdots&\vdots& \ddots& \vdots\\
0& 0& 0&\cdots& W_{n,2}
\end{pmatrix} 
,
\begin{pmatrix}
B_{1,2}\\B_{2,2}\\B_{3,2}\\\vdots\\ B_{n,2}
\end{pmatrix}}
,
\ldots,
\\
&
\quad
\left.
\pa{\begin{pmatrix}
W_{1,L}& 0& 0& \cdots& 0\\
0& W_{2,L}& 0&\cdots& 0\\
0& 0& W_{3,L}&\cdots& 0\\
\vdots& \vdots&\vdots& \ddots& \vdots\\
0& 0& 0&\cdots& W_{n,L}
\end{pmatrix} 
,
\begin{pmatrix}
B_{1,L}\\B_{2,L}\\B_{3,L}\\\vdots\\ B_{n,L}
\end{pmatrix}}
\right)
\nonumber
\end{align}
\cfout.
\end{definition}

\subsection{Summations of ANNs}

\cfclear
\begin{definition}[Affine linear transformation ANN]\label{linear}
	\cfconsiderloaded{linear}
	Let 
	$m,n\in\N$, 
	$\weight\in\R^{m\times n}$, 
	$\bias\in\R^m$. 
	Then we denote by 
	$\A_{\weight,\bias} \in (\R^{m\times n}\times \R^m) \subseteq \ANNs$ the \network
	given by $\A_{\weight,\bias} = (\weight,\bias)$
	\cfout.
\end{definition}

\cfclear
\begin{definition}\label{def:ANN:sum}
\cfconsiderloaded{def:ANN:sum}
Let $m, n \in \N$. 
Then we denote by $\sumANN_{m, n} \in (\R^{m \times (nm)} \times \R^m)$ the \network
given by 
$\sumANN_{m, n} = \A_{(\idMatrix_m \,\,\,  \idMatrix_m \,\,\, \ldots \,\,\, \idMatrix_m), 0 }$
\cfout.
\end{definition}

\cfclear
\begin{definition}[Matrix transpose]\label{def:Transpose}
\cfconsiderloaded{def:Transpose}
Let 
$m, n \in \N$, 
$A \in \R^{m \times n}$. 
Then we denote by $A^\transpose \in \R^{n \times m}$ the transpose of A.
\end{definition}

\cfclear
\begin{definition}[Transpose ANN]\label{def:ANN:extension}
\cfconsiderloaded{def:ANN:extension}
Let $m, n \in \N$. 
Then we denote by $\extensionANN_{m, n} \in (\R^{(nm) \times m} \times \R^{nm})$	the \network given by 
$\extensionANN_{m, n} = \A_{(\idMatrix_m \,\,\,  \idMatrix_m \,\,\, \ldots \,\,\, \idMatrix_m)^\transpose , 0}$ 
\cfout.
\end{definition}

\cfclear
\begin{definition}[Sums of ANNs with the same length]
\label{def:ANNsum:same}
\cfconsiderloaded{def:ANNsum:same}
Let 
$\lbd \in \Z$, 
$\ubd \in \Z \cap [\lbd,\infty)$, 
$\Phi_\lbd, \Phi_{\lbd+1}, \ldots, \Phi_\ubd \in \ANNs$ 
satisfy for all  
$k \in \Z \cap [\lbd,\ubd]$ 
that
$\lengthANN(\Phi_k) = \lengthANN(\Phi_\lbd)$,
$\inDimANN(\Phi_k) = \inDimANN(\Phi_\lbd)$, 
and 
$\outDimANN(\Phi_k) = \outDimANN(\Phi_\lbd)$ \cfload.
Then we denote by $\oSum_{k=\lbd}^\ubd \Phi_k$ (we denote by $\Phi_\lbd \oSum \Phi_{\lbd+1} \oSum \ldots \oSum \Phi_\ubd$)
the \network given by
\begin{equation}
\OSum{k=\lbd}{\ubd} \Phi_k 
= 
\pr[\Big]{ \compANN{\sumANN_{\outDimANN(\Phi_\lbd), \ubd-\lbd+1}}{{\compANN{\big[\parallelizationSpecial_{\ubd-\lbd+1}(\Phi_\lbd,\Phi_{\lbd+1},\dots, \Phi_\ubd)\big]}{\extensionANN_{\inDimANN(\Phi_\lbd), \ubd-\lbd+1}}}} } 
\in 
\ANNs
\end{equation}
\cfout.
\end{definition}

\cfclear
\begin{definition}[Sums of ANNs with different lengths]\label{def:ANN:sum_diff}
\cfconsiderloaded{def:ANN:sum_diff}
Let 
$\lbd\in\Z$, 
$\ubd \in \Z \cap [\lbd,\infty)$, 
$\Phi_\lbd, \allowbreak \Phi_{\lbd+1}, \allowbreak \ldots, \allowbreak \Phi_\ubd, \Psi \in \ANNs$ 
satisfy
for all $k \in \Z \cap [\lbd,\ubd]$ that
$\inDimANN(\Phi_k) = \inDimANN(\Phi_\lbd)$, 
$\outDimANN(\Phi_k) = \inDimANN(\Psi) = \outDimANN(\Psi)$, 
and 
$\hiddenLength(\Psi) = 1$ \cfload.
Then we denote by $\bSum_{k=\lbd,\Psi}^\ubd\Phi_k$ (we denote by $\Phi_\lbd\bSum_\Psi \Phi_{\lbd+1} \bSum_\Psi \ldots \bSum_\Psi \Phi_\ubd$) the \network given by
\begin{equation}
\BSum{k=\lbd}{\Psi}{\ubd} \Phi_k 
= 
\OSum{k=\lbd}{\ubd} \longerANN{\max_{j\in\{\lbd,\lbd+1,\dots,\ubd\}}\lengthANN(\Phi_j),\Psi}(\Phi_k) 
\in 
\ANNs
\end{equation}
\cfout.
\end{definition}

\subsection{Linear combinations of ANNs}

\subsubsection{Linear combinations of ANNs with the same length}

\cfclear
\begin{definition}[Scalar multiplications of ANNs]
\label{def:ANNscalar}
\cfconsiderloaded{def:ANNscalar}
We denote by $\scalar{\left(\cdot\right)}{ \left(\cdot\right)} \colon \R \times \ANNs \to \ANNs$ the function which satisfies for all $\lambda \in \R$, $\Phi \in \ANNs$ that 
$\scalar{\lambda}{\Phi} = \compANN{\A_{\lambda \idMatrix_{\outDimANN(\Phi)},0}}{\Phi}$ 
\cfout.
\end{definition}

\cfclear
\begin{athm}{lemma}{lem:sum:ANN}
Let 
$ \lbd \in \Z $, 
$\ubd \in \Z \cap [\lbd,\infty)$, 
$n = \ubd-\lbd+1$, 
$h_\lbd, h_{\lbd+1}, \ldots, h_\ubd \in \R$, 
$t_\lbd, t_{\lbd+1}, \ldots, t_\ubd \in \R$, 
$\Phi_\lbd, \Phi_{\lbd+1}, \ldots, \allowbreak \Phi_\ubd , \allowbreak \Psi \in \ANNs$, 
$B_\lbd,B_{\lbd+1},\ldots,B_\ubd \in \R^{\inDimANN(\Phi_\lbd)}$ 
satisfy 
$\dims(\Phi_\lbd) = \dims(\Phi_{\lbd+1}) = \ldots = \dims(\Phi_\ubd)$ 
and
\begin{equation}
\Psi 
= 	
\OSum{k=\lbd}{\ubd} \pr[\Big]{ \scalar{h_k}{\pr[\big]{ \compANN{\Phi_k}{\A_{t_k\idMatrix_{\inDimANN(\Phi_k)},B_k}} }} }
\end{equation}
\cfload.
Then 
\begin{enumerate}[label=(\roman{*})]
\item\label{item:sum:ANN:1} 
it holds that
\begin{equation}
\begin{split}
\dims (\Psi) 
& 
= 
\bigl(\inDimANN(\Phi_\lbd), \smallsum_{k=\lbd}^\ubd \dimANNlevel_1(\Phi_\lbd), \smallsum_{k=\lbd}^\ubd\dimANNlevel_2(\Phi_\lbd), \ldots, \smallsum_{k=\lbd}^\ubd\dimANNlevel_{\lengthANN(\Phi_\lbd)-1}(\Phi_\lbd), \outDimANN (\Phi_\lbd) \bigr) 
\\
& 
= 
\bigl(\inDimANN(\Phi_\lbd), n\dimANNlevel_1(\Phi_\lbd), n\dimANNlevel_2(\Phi_\lbd), \ldots, n\dimANNlevel_{\lengthANN(\Phi_\lbd)-1}(\Phi_\lbd), \outDimANN (\Phi_\lbd)\bigr)
\dc
\end{split}
\end{equation}
\item\label{item:sum:ANN:3} 
it holds for all 	
$ a \in C(\R, \R)$ 
that  
$\functionANN{\activation}(\Psi) \in C(\R^{\inDimANN(\Phi_\lbd)}, \R^{\outDimANN(\Phi_\lbd)})$, 
and
\item\label{item:sum:ANN:4} 
it holds for all 
$ a \in C(\R, \R)$, 
$x \in \R^{\inDimANN(\Phi_\lbd)}$ 
that  
\begin{equation}
(\functionANN{\activation} (\Psi))(x)
= 
\SmallSum{k=\lbd}{\ubd} h_k (\functionANN{\activation} (\Phi_k))(t_k x+B_k)
\end{equation}
\end{enumerate}	
\cfout.
\end{athm}

\begin{aproof}
First, \nobs that 
\Enum{the hypothesis that 
$\dims(\Phi_\lbd) = \dims(\Phi_{\lbd+1}) = \ldots = \dims(\Phi_\ubd)$ 
}
that for all 
$k \in \{\lbd, \lbd+1, \allowbreak \ldots, \allowbreak \ubd\}$ it holds that
\begin{equation}
\dims \pr[\big]{ \A_{t_k\idMatrix_{\inDimANN(\Phi_k)},B_k} }
=
\dims \pr[\big]{ \A_{t_k\idMatrix_{\inDimANN(\Phi_\lbd)},B_k} } 
= 
\pr[\big]{ \inDimANN(\Phi_\lbd),  \inDimANN(\Phi_\lbd) }
\in 
\N^2
\dpp
\end{equation}
\Enum{
This
;
\eg Grohs et al.\ \cite[item (i) in Proposition~2.6]{GrohsHornungJentzen2019}
}
that for all 
$k \in \{\lbd, \lbd+1, \ldots, \ubd\}$ 
it holds that
\begin{equation}\label{eq:lem:resultsum:dims0}
\dims \pr[\big]{ \compANN{\Phi_k}{\A_{t_k\idMatrix_{\inDimANN(\Phi_k)},B_k}} } 
= 
\pr[\big]{ \inDimANN(\Phi_\lbd), \dimANNlevel_1(\Phi_\lbd), \dimANNlevel_2(\Phi_\lbd), \ldots, \dimANNlevel_{\lengthANN(\Phi_\lbd)}(\Phi_\lbd) }
\dpp
\end{equation}
\Nobs that
\Enum{\eg Grohs et al.\ \cite[item (i) in Lemma~3.14]{GrohsJentzenSalimova2019}
} that for all 
$k \in \{\lbd, \lbd+1, \ldots, \ubd\}$ 
it holds that
\begin{equation}
\label{eq:lem:resultsum:dims}
\dims \pr[\Big]{ h_k  \circledast \pr[\big]{ \compANN{\Phi_k}{\A_{t_k\idMatrix_{\inDimANN(\Phi_k)},B_k}} } } 
= 
\dims \pr[\big]{ \compANN{\Phi_k}{\A_{t_k\idMatrix_{\inDimANN(\Phi_k)},B_k}} } 
\dpp
\end{equation}
Combining 
\Enum{this; 
\cref{eq:lem:resultsum:dims0};
\eg
Grohs et al.\ \cite[item (ii) in Lemma~3.28]{GrohsJentzenSalimova2019}
} that
\begin{equation}
\begin{split}
\dims (\Psi) 
&
= 
\dims \pr[\Big]{ \oSum_{k=\lbd}^\ubd \pr[\big]{ h_k  \circledast ( \compANN{\Phi_k}{\A_{t_k\idMatrix_{\inDimANN(\Phi_k)},B_k}}) } }
\\
& 
= 
\bigl(\inDimANN(\Phi_\lbd), \smallsum_{k=\lbd}^\ubd \dimANNlevel_1(\Phi_\lbd), \smallsum_{k=\lbd}^\ubd\dimANNlevel_2(\Phi_\lbd), \ldots, \smallsum_{k=\lbd}^\ubd\dimANNlevel_{\lengthANN(\Phi_\lbd)-1}(\Phi_\lbd), \outDimANN (\Phi_\lbd) \bigr) 
\\
& 
=  
\pr[\big]{ \inDimANN(\Phi_\lbd), n\dimANNlevel_1(\Phi_\lbd), n\dimANNlevel_2(\Phi_\lbd), \ldots, n\dimANNlevel_{\lengthANN(\Phi_\lbd)-1}(\Phi_\lbd), \outDimANN (\Phi_\lbd) } 
\dpp
\end{split}
\end{equation}
This establishes \cref{item:sum:ANN:1}.
\Moreover
\Enum{
	\eg 
	Grohs et al.\ \cite[item (v) in Proposition~2.6]{GrohsHornungJentzen2019} 
} that
for all 
$k \in \{\lbd, \lbd+1, \ldots, \ubd\}$, 
$a \in C(\R, \R)$, 
$x \in \R^{\inDimANN(\Phi_\lbd)}$ 
it holds that 
$\functionANN{\activation} ( \compANN{\Phi_k}{\A_{t_k\idMatrix_{\inDimANN(\Phi_k)},B_k}}) \in C(\R^{ \inDimANN(\Phi_\lbd)}, \R^{\outDimANN(\Phi_\lbd)})$ 
and
\begin{equation}
\pr[\Big]{ \functionANN{\activation} \pr[\big]{ \compANN{\Phi_k}{\A_{t_k\idMatrix_{\inDimANN(\Phi_k)},B_k}} } } \lrSpace (x) 
= 
(\functionANN{\activation} (\Phi_k))(t_k x + B_k)
\end{equation}
\cfload.
Combining 
\Enum{this and, e.g., Grohs et al.\ \cite[Lemma~3.14]{GrohsJentzenSalimova2019} 
} that for all 
$k \in \{\lbd, \lbd+1, \ldots, \ubd\}$, 
$a \in C(\R, \R)$, $x \in \R^{ \inDimANN(\Phi_\lbd)}$ 
it holds that  
\begin{equation}
\functionANN{\activation} \pr[\Big]{ h_k  \circledast \pr[\big]{ \compANN{\Phi_k}{\A_{t_k\idMatrix_{\inDimANN(\Phi_k)},B_k}} } } 
\in 
C(\R^{\inDimANN(\Phi_\lbd)}, \R^{\outDimANN(\Phi_\lbd)})
\end{equation}
and
\begin{equation}
\pr[\Big]{ \functionANN{\activation} \pr[\big]{ h_k  \circledast ( \compANN{\Phi_k}{\A_{t_k\idMatrix_{\inDimANN(\Phi_k)},B_k}}) } } \lrSpace (x) 
= 
h_k (\functionANN{\activation} (\Phi_k)) (t_k x + B_k)
\dpp
\end{equation}
\Moreover
\Enum{\eg Grohs et al.\ \cite[Lemma~3.28]{GrohsJentzenSalimova2019}; 
\cref{eq:lem:resultsum:dims}}
that for all  
$a \in C(\R, \R)$, 
$x \in \R^{ \inDimANN(\Phi_\lbd)}$ 
it holds that 
$\functionANN{\activation}(\Psi) \in C(\R^{ \inDimANN(\Phi_\lbd)}, \R^{\outDimANN(\Phi_\lbd)})$  
and
\begin{equation}
\label{eq:sum:ann}
\begin{split}
(\functionANN{\activation} (\Psi))(x) 
& 
= 
\pr[\Big]{ \functionANN{\activation} \pr[\big]{\oSum_{k=\lbd}^\ubd ( h_k  \circledast ( \compANN{\Phi_k}{\A_{t_k\idMatrix_{\inDimANN(\Phi_k)},B_k}})) } } \lrSpace (x)
\\
&
=  
\SmallSum{k = \lbd}{\ubd} 
\pr[\Big]{ \functionANN{\activation} \pr[\big]{ h_k  \circledast ( \compANN{\Phi_k}{\A_{t_k\idMatrix_{\inDimANN(\Phi_k)},B_k}}) } } \lrSpace (x) 
= 
\SmallSum{k=\lbd}{\ubd} h_k (\functionANN{\activation} (\Phi_k))(t_k x+B_k)
\dpp
\end{split}
\end{equation}
This establishes \cref{item:sum:ANN:3,item:sum:ANN:4}.
\end{aproof}

\subsubsection{Linear combinations of ANNs with different lengths}

\cfclear
\begin{athm}{lemma}{lemma14a}
Let 
$ L \in \N$, 
$\lbd \in \Z$, 
$\ubd \in \Z \cap [\lbd,\infty)$, 
$h_\lbd,h_{\lbd+1},\ldots,h_\ubd\in\R$, 
$\Phi_\lbd, \allowbreak \Phi_{\lbd+1}, \allowbreak \ldots, \allowbreak \Phi_\ubd, \fJ, \Psi \in \ANNs$, 
$B_\lbd,B_{\lbd+1},\ldots,B_\ubd \in \R^{\inDimANN(\Phi_\lbd)}$, 
$a \in C(\R,\R)$ 
satisfy for all 
$j \in \Z \cap [\lbd,\ubd]$ 
that
$L = \max_{k\in\Z \cap [\lbd,\ubd]} \lengthANN(\Phi_k)$,
$\inDimANN(\Phi_j) = \inDimANN(\Phi_\lbd)$, 
$\outDimANN(\Phi_j) = \inDimANN(\fJ) = \outDimANN(\fJ)$, 
$\hiddenLength(\fJ) = 1$, 
$\Ra(\fJ) = \operatorname{id}_\R$, 
and
\begin{equation}\label{intro_2_33}
\Psi 
= 
\BSum{k=\lbd}{\fJ}{\ubd} \pr[\Big]{ \scalar{h_k}{ \pr[\big]{ \compANN{\Phi_k}{\A_{\idMatrix_{\inDimANN(\Phi_k)},B_k}} } } }
\end{equation}
\cfload. 
Then
\begin{enumerate}[label=(\roman*)]
\item\label{lemma14a_1} 
it holds that
\begin{align}
& 
\dims(\Psi) 
\\
& 
= 
\Bigl( \inDimANN(\Phi_\lbd), \SmallSum{k=\lbd}{\ubd} \dimANNlevel_1\pr[\big]{ \longerANN{L,\fJ}(\Phi_k) }, \SmallSum{k=\lbd}{\ubd} \dimANNlevel_2\pr[\big]{ \longerANN{L,\fJ}(\Phi_k) } ,  \ldots, \SmallSum{k=\lbd}{\ubd} \dimANNlevel_{L - 1} \pr[\big]{ \longerANN{L,\fJ}(\Phi_k) } , \outDimANN(\Phi_\lbd) \Bigr) 
\dc 
\nonumber
\end{align}
\item\label{lemma14a_3} 
it holds that $\Ra(\Psi) \in C(\R^{\inDimANN(\Phi_\lbd)},\R^{\outDimANN(\Phi_\lbd)})$, and
\item\label{lemma14a_4} 
it holds for all $x \in \R^{\inDimANN(\Phi_\lbd)}$ that
\begin{equation}
(\Ra(\Psi))(x) 
= 
\SmallSum{k=\lbd}{\ubd} h_k(\Ra(\Phi_k))(x+B_k)
\end{equation}
\end{enumerate}
\cfout.
\end{athm}

\begin{aproof}
\newcommand{\pp}{\Theta}
\newcommand{\ppp}{\Xi}
\newcommand{\llll}{\mathbb{L}}
\Nobs that
\cref{item:sum:ANN:1} in \cref{lem:sum:ANN}
establishes \cref{lemma14a_1}.
\Moreover
\Enum{
	\eg Grohs et al.\ \cite[item (v) in Proposition~2.6]{GrohsHornungJentzen2019}
} 
that
for all 
$k \in \Z \cap [\lbd,\ubd]$, $x \in \R^{\inDimANN(\Phi_\lbd)}$ 
it holds that 
$\functionANN{\activation} ( \compANN{\Phi_k}{\A_{\idMatrix_{\inDimANN(\Phi_k)},B_k}}) \in C(\R^{ \inDimANN(\Phi_\lbd)}, \R^{\outDimANN(\Phi_\lbd)})$ 
and
\begin{equation}
\pr[\Big]{ \functionANN{\activation} \pr[\big]{ \compANN{\Phi_k}{\A_{\idMatrix_{\inDimANN(\Phi_k)},B_k}} } } \lrSpace (x) 
= 
(\functionANN{\activation} (\Phi_k))(x + B_k) 
\dpp
\end{equation}
\Enum{
This;
e.g., Grohs et al.\ \cite[Lemma~3.14]{GrohsJentzenSalimova2019};
\eg Grohs et al.\ \cite[item~(ii) in Lemma~2.14]{GrohsHornungJentzen2019}
}
that for all 
$k \in \N \cap [\lbd,\ubd]$, 
$x \in \R^{ \inDimANN(\Phi_\lbd)}$ 
it holds that  
\begin{equation} 
\functionANN{\activation} \pr[\Big]{ \longerANN{L,\fJ}\pr[\big]{ h_k  \circledast ( \compANN{\Phi_k}{\A_{\idMatrix_{\inDimANN(\Phi_k)},B_k}}) } } 
=
\functionANN{\activation} \pr[\big]{ h_k  \circledast ( \compANN{\Phi_k}{\A_{\idMatrix_{\inDimANN(\Phi_k)},B_k}}) } 
\in 
C(\R^{\inDimANN(\Phi_\lbd)}, \R^{\outDimANN(\Phi_\lbd)})
\end{equation}
and
\begin{equation}
\begin{split}
\pr[\Big]{ \functionANN{\activation} \pr[\Big]{ \longerANN{L,\fJ}\pr[\big]{ h_k  \circledast ( \compANN{\Phi_k}{\A_{\idMatrix_{\inDimANN(\Phi_k)},B_k}}) } } } \lrSpace (x) 
& 
=
\pr[\Big]{ \functionANN{\activation} \pr[\big]{ h_k  \circledast ( \compANN{\Phi_k}{\A_{\idMatrix_{\inDimANN(\Phi_k)},B_k}}) } } \lrSpace (x) 
\\
& 
= 
h_k (\functionANN{\activation} (\Phi_k)) (x + B_k)
\end{split}
\end{equation}
\cfload.
Combining 
\Enum{
this, e.g., Grohs et al.\ \cite[Lemma~3.28]{GrohsJentzenSalimova2019}, and \cref{eq:lem:resultsum:dims} 
}
that for all $x \in \R^{ \inDimANN(\Phi_\lbd)}$ it holds that 
$\functionANN{\activation}(\Psi) \in C(\R^{ \inDimANN(\Phi_\lbd)}, \R^{\outDimANN(\Phi_\lbd)})$ 
and
\begin{align}
(\Ra(\Psi))(x)
& 
= 
\pr*{ \Ra\pr*{ \BSum{k=\lbd}{\fJ}{\ubd} \pr[\Big]{ \scalar{h_k}{\pr[\big]{ \compANN{\Phi_k}{\A_{\idMatrix_{\inDimANN(\Phi_k)},B_k}} } } } } } \lrSpace (x) 
\nonumber
\\
& 
= 
\pr*{ \Ra\pr*{ \OSum{k=\lbd}{\ubd} \longerANN{L,\fJ}\pr[\Big]{ \scalar{h_k}{ \pr[\big]{ \compANN{\Phi_k}{\A_{\idMatrix_{\inDimANN(\Phi_k)},B_k}} } } } } } \lrSpace (x) 
\\
& 
= 
\SmallSum{k = \lbd}{\ubd} \pr*{ \Ra\pr*{ \longerANN{L,\fJ}\pr[\Big]{ \scalar{h_k}{ \pr[\big]{ \compANN{\Phi_k}{\A_{\idMatrix_{\inDimANN(\Phi_k)},B_k}} } } } } } \lrSpace (x) 
= 
\SmallSum{k = \lbd}{\ubd} h_k (\functionANN{\activation} (\Phi_k)) (x + B_k) 
\nonumber
\end{align}
\cfload.
This establishes \cref{lemma14a_3,lemma14a_4}.
\end{aproof}

%
%
%


\section{ANN representations for MLP approximations}
\label{subsubsec:ANN_repres_for_MLP}

In this section we study deep ANN representations for MLP approximations for PDEs of the form~\cref{eq:1}. 
Specifically, in \cref{lemma15} below we show that realizations of suitable deep ANNs 
with a given general activation function 
coincide with appropriate MLP approximations for PDEs of the form~\cref{eq:1}
provided that with the considered general activation function there exists a shallow ANN representation for the one-dimensional identity function $\R \ni x \mapsto x \in \R$. 

In the elementary results in \cref{subsec:activANN} and \cref{subsec:ANNforId} we explicitly construct such shallow ANN representations for the one-dimensional identity function $\R \ni x \mapsto x \in \R$ in the case of several activation functions such as the ReLU (see \cref{lem:ReluLeaky:identity}), the leaky ReLU (see \cref{lem:ReluLeaky:identity}), the rectified power unit (RePU) (see \cref{lem:ReluPower:identity}), and the softplus (see \cref{lem:SoftPlus:identity}) activation functions.
In the special case of the ReLU activation results similar
 to \cref{lemma15} can, e.g., be found in   
Hutzenthaler et al.~\cite[Lemma 3.10]{HutzenthalerJentzenKruse2019} and Cioica-Licht et al.~\cite[Lemma 3.10]{cioicalicht2022deep}.

\subsection{Activation functions as ANNs}\label{subsec:activANN}

\cfclear
\begin{definition}[Activation ANN]\label{padding}
\cfconsiderloaded{padding}
Let $n\in\N$.
Then we denote by
$\ii_n \in ((\R^{n\times n}\times \R^n)\times (\R^{n\times n}\times \R^n)) \subseteq \ANNs$
the \network given by
$\ii_n = ((\idMatrix_n,0),(\idMatrix_n,0)) $
\cfout.
\end{definition}

\cfclear
\begin{athm}{lemma}{padding_lemma}
Let $n\in\N$.
Then
\begin{enumerate}[label=(\roman *)]
\item
\label{padding_item2}
it holds that $\dims(\ii_n) = (n,n,n) \in \N^3$,
\item
\label{padding_item3}
it holds for all $\activation\in C(\R,\R)$ that $\functionANN{\activation}(\ii_n) \in C(\R^n,\R^n)$, and
\item
\label{padding_item4}
it holds for all $\activation\in C(\R,\R)$ that 
$
\functionANN{\activation}(\ii_n) = \activationDim{n}
$ 
\end{enumerate}
\cfout.
\end{athm}

\begin{aproof}
\Nobs that
\Enum{
the fact that $\ii_n \in ((\R^{n\times n}\times\R^n)\times (\R^{n\times n}\times\R^n)) \subseteq \NN$
} 
that 
$\dims(\ii_n) = (n,n,n) \in \N^3$ \cfload. 
This establishes \cref{padding_item2}.
\Moreover 
\Enum{
the fact that
$\ii_n = ((\idMatrix_n,0),(\idMatrix_n,0)) \in ((\R^{n\times n}\times\R^n)\times (\R^{n\times n}\times\R^n))$
; 
\cref{setting_NN:ass2}
}
that for all 
$\activation \in C(\R,\R)$, 
$x\in\R^n$ 
it holds that 
$\functionANN{\activation}(\ii_n) \in C(\R^n,\R^n)$ 
and
\begin{equation}
(\functionANN{\activation}(\ii_n))(x) 
= 
\idMatrix_n \left(\activationDim{n}(\idMatrix_n \lrSpace x + 0) \right) 
+ 
0 
= 
\activationDim{n}(x)
\end{equation}
\cfload.
This establishes \cref{padding_item3,padding_item4}. 
\end{aproof}

\cfclear
\begin{athm}{lemma}{padding_lemma2}
Let $\Phi\in\ANNs$ \cfload.
Then
\begin{enumerate}[label=(\roman *)]
\item
\label{padding2_item1}
it holds that
\begin{equation}
\dims\pr[\big]{ \compANN{\ii_{\outDimANN(\Phi)}}{\Phi} }
= 
\pr[\big]{
\inDimANN(\Phi),\dimANNlevel_1(\Phi),\dimANNlevel_2(\Phi),\ldots,\dimANNlevel_{\lengthANN(\Phi)-1}(\Phi),\outDimANN(\Phi),\outDimANN(\Phi)
} 
\in 
\N^{\lengthANN(\Phi) + 2}
\dc
\end{equation}
\item
\label{padding2_item3}
it holds for all 
$\activation\in C(\R,\R)$ 
that 
$\functionANN{\activation}(\compANN{\ii_{\outDimANN(\Phi)}}{\Phi}) \in C(\R^{\inDimANN(\Phi)},\R^{\outDimANN(\Phi)})$,
\item
\label{padding2_item2}
it holds for all 
$\activation\in C(\R,\R)$, 
$x\in\R^{\inDimANN(\Phi)}$ 
that 
$
(\functionANN{\activation}(\compANN{\ii_{\outDimANN(\Phi)}}{\Phi}))(x) = \activationDim{\outDimANN(\Phi)}((\Ra(\Phi))(x))
$,
\item
\label{padding2_item4}
it holds that
\begin{equation}
\dims\pr[\big]{ \compANN{\Phi}{\ii_{\inDimANN(\Phi)}} } 
= 
\pr[\big]{
\inDimANN(\Phi),\inDimANN(\Phi),\dimANNlevel_1(\Phi),\dimANNlevel_2(\Phi),\ldots,\dimANNlevel_{\lengthANN(\Phi)-1}(\Phi),\outDimANN(\Phi)
} 
\in 
\N^{\lengthANN(\Phi) + 2}
\dc
\end{equation}
\item
\label{padding2_item6}
it holds for all 
$\activation \in C(\R,\R)$ 
that 
$\functionANN{\activation}(\compANN{\Phi}{\ii_{\inDimANN(\Phi)}}) \in C(\R^{\inDimANN(\Phi)},\R^{\outDimANN(\Phi)})$,
and
\item
\label{padding2_item5}
it holds for all 
$\activation\in C(\R,\R)$, 
$x\in\R^{\inDimANN(\Phi)}$ 
that
$
(\functionANN{\activation}(\compANN{\Phi}{\ii_{\inDimANN(\Phi)}}))(x) = (\Ra(\Phi))(\activationDim{\inDimANN(\Phi)}(x))
$
\end{enumerate}
\cfout.
\end{athm}

\begin{aproof}
\Nobs 
that 
\Enum{
\cref{padding_lemma}
}
that for all 
$n\in\N$, 
$\activation \in C(\R,\R)$, 
$x\in\R^n$ 
it holds that 
$\functionANN{\activation}(\ii_n) \in C(\R^n,\R^n)$ 
and
\begin{equation}
(\functionANN{\activation}(\ii_n))(x) 
= 
\activationDim{n}(x)
\end{equation}
\cfload.
Combining 
this and, 
e.g., Grohs et al.\ \cite[Proposition~2.6]{GrohsHornungJentzen2019} 
establishes
\cref{padding2_item1,%
padding2_item2,%
padding2_item3,%
padding2_item4,%
padding2_item5,%
padding2_item6}.
\end{aproof}

\subsection{ANN representations for the one-dimensional identity function}\label{subsec:ANNforId}

\cfclear
\begin{definition}[RePU monomial ANNs]\label{def:id_net}
\label{def:ReLu:identity}
\cfconsiderloaded{def:id_net}
Let $\gamma \in \N_0$.
Then we denote by $\I_\gamma \in \ANNs$ the \network given by 
\begin{equation}\label{eq:def:id:1}
\I_\gamma
= 
\left( \left(\begin{pmatrix}
1 \\
-1
\end{pmatrix}
, 
\begin{pmatrix}
0 \\
0
\end{pmatrix}
\right)
,
\Big(	
\begin{pmatrix}
1 & (-1)^\gamma
\end{pmatrix}
, 
0 \Big)
\right)  
\in 
\big((\R^{2 \times 1} \times \R^{2}) \times (\R^{1 \times 2} \times \R^1) \big)
\end{equation} 
\cfout.
\end{definition}

\cfclear
\begin{athm}{lemma}{lem:ReluLeaky:identity}[Shallow leaky ReLU ANN representation for the one-dimensional identity function]
Let 
$\alpha \in [0,\infty)$,
$\activation \in C(\R, \R)$,
$\Psi \in \ANNs$
satisfy for all $x \in \R$
that 
$\activation(x) = \max\{ x , \alpha x \}$
and
$\Psi = \scalar{ (1 + \alpha)^{-1} }{ \I_1 }$
\cfload.
Then
\begin{enumerate}[label=(\roman{*})]
\item
\label{item:lem:Relu:dims} 
it holds for all $\gamma \in \N_0$ that $\dims(\I_\gamma) = (1,2,1) \in \N^3$,
\item
\label{item:lem:ReluLeaky:identity}
it holds for all 
$x \in \R$ 
that 
$(\functionANN{\activation}(\I_1))(x) = (1+\alpha)x$,
\item
\label{item:lem:ReluLeaky:dims} 
it holds that 
$\dims(\Psi) = (1,2,1) \in \N^3$,
\item
\label{item:lem:ReluLeaky:cont} 
it holds that 
$ \functionANN{\activation}(\Psi) \in C(\R, \R)$, 
and
\item
\label{item:lem:ReluLeaky:real}
it holds for all 
$x \in \R$ 
that 
$
(\functionANN{\activation}(\Psi))(x) = x
$
\end{enumerate}
\cfout.
\end{athm}

\begin{aproof}
\Nobs that \Enum{\cref{eq:def:id:1}} that for all $\gamma \in \N_0$ it holds that
$
\dims (\I_\gamma) = (1, 2, 1) \in \N^3
$.
This establishes \cref{item:lem:Relu:dims}.
\Moreover \Enum{\cref{eq:def:id:1}} that for all $x \in \R$ it holds that
\begin{equation}\label{eq:Relu:3_6}
\begin{split}
(\functionANN{\activation}(\I_1))(x) 
= 
\activation(x) - \activation(-x) 
&
= 
\max\{x, \alpha x\} - \max\{-x, - \alpha x \} 
\\
&
=
\max\{x, \alpha x\} + \min\{x, \alpha x \} 
= 
(1 + \alpha)
x
\end{split}
\end{equation}
\cfload.
This establishes 
\cref{item:lem:ReluLeaky:identity}.
\Moreover 
\cref{item:sum:ANN:1} in \cref{lem:sum:ANN}
(applied with 
$\lbd \with 1$, 
$\ubd \with 1$, 
$h_\lbd \with (1+\alpha)^{-1}$, 
$t_\lbd \with 1$, 
$\Phi_\lbd \with \I_1$, 
$B_\lbd \with 0$, 
$\Psi \with \Psi$
in the notation of \cref{lem:sum:ANN})
establishes
\cref{item:lem:ReluLeaky:dims}.
\Enum{
Combining
\cref{eq:Relu:3_6}
and
\cref{item:sum:ANN:4} in \cref{lem:sum:ANN}
(applied with 
$\lbd \with 1$, 
$\ubd \with 1$, 
$h_\lbd \with (1+\alpha)^{-1}$, 
$t_\lbd \with 1$, 
$\Phi_\lbd \with \I_1$, 
$B_\lbd \with 0$, 
$\Psi \with \Psi$, 
$\activation \with \activation$
in the notation of \cref{lem:sum:ANN})
\hence\unskip
}[prove]
\cref{item:lem:ReluLeaky:cont,item:lem:ReluLeaky:real}.
\end{aproof}

\cfclear
\begin{athm}{lemma}{lem:Power:identity}[Shallow power ANN representation for the one-dimensional identity function]
Let 
$\gamma \in \N\backslash\{1\}$,
$b_1,b_2,\dots,b_\gamma \in \R$,
$\activation \in C(\R, \R)$
satisfy for all 
$x \in \R$
that 
$\activation(x) = x^\gamma$
and 
$b_1 < b_2 < \ldots < b_\gamma$
\cfload.
Then
\begin{enumerate}[label=(\roman{*})]
\item 
\label{item:lem:Power:constant_exist}
there exist unique 
$c_0,c_1,\dots,c_\gamma \in \R$ 
which satisfy for all 
$k \in \{0,1,\dots,\gamma\}$ 
that
$\1_{\{\gamma\}}(k) \, c_0 + \sum_{i=1}^\gamma c_i (b_i)^{k} = \1_{\{\gamma-1\}}(k)\,\gamma^{-1}$,
\item 
\label{item:lem:Power:exist}
there exists a unique $\Psi \in \ANNs$ which satisfies
\begin{equation}
\Psi
= 
\compANN{\A_{1,c_0}}{\pr[\bigg]{ \OSum{i=1}{\gamma} \pr[\Big]{ \scalar{c_i}{ \pr[\big]{ \compANN{\ii_1}{\A_{1,b_i}} } } } } }
\dc
\end{equation}
\item
\label{item:lem:Power:dims} 
it holds that 
$\dims(\Psi) = (1,\gamma,1) \in \N^3$,
\item
\label{item:lem:Power:cont} 
it holds that 
$ \functionANN{\activation}(\Psi) \in C(\R, \R)$, 
and
\item
\label{item:lem:Power:real}
it holds for all 
$x \in \R$ 
that 
$
(\functionANN{\activation}(\Psi))(x) = x
$
\end{enumerate}
\cfout.
\end{athm}

\begin{aproof}
\newcommand{\BBB}{{\bf B}}
\newcommand{\CCC}{{\bf C}}
\newcommand{\DDD}{{\bf D}}
Throughout this proof let 
$\BBB = (\BBB_{i,j})_{i,j\in\{1,2,\dots,\gamma+1\}} \in \R^{(\gamma+1)\times(\gamma+1)}$
satisfy for all 
$i,j \in \{1,2,\dots,\gamma\}$
that
$\BBB_{1,i+1} = 1$,
$\BBB_{i,1} = 0$,
$\BBB_{\gamma+1,1} = 1$,
and
$\BBB_{i+1,j+1} = (b_j)^{i}$ 
and let
$\DDD = (\DDD_1,\DDD_2,\dots,\DDD_{\gamma+1})^\transpose \in \R^{(\gamma+1)\times 1}$ 
satisfy for all 
$k \in \{1,2,\dots,\gamma+1\}$ 
that
$\DDD_k = \1_{\{\gamma\}}(k) \, \gamma^{-1}$
\cfload.
\Nobs that
\Enum{
the assumption that
$b_1 < b_2 < \ldots < b_\gamma$
;
\eg
Horn and Johnson \cite[Eq.~(0.9.11.2)]{MR832183}
}
that
\begin{equation}
\det(\BBB)
=
(-1)^{\gamma+1} \det\pr[\big]{ (\BBB_{i,j+1})_{i,j\in\{1,2,\dots,\gamma\}} }
=
(-1)^{\gamma}
\br*{
\textstyle
\prod_{\substack{i,j \in \{ 1,2,\dots,\gamma \} \\ i < j }} \pr[\big]{ b_j - b_i }
}
\neq 
0
\dpp
\end{equation}
\Enum{
This}
that there exists a unique
$\CCC = (c_0 , c_1 , \dots , c_\gamma)^\transpose \in \R^{(\gamma+1)\times 1} $
such that
$\BBB \CCC = \DDD$.
This establishes \cref{item:lem:Power:constant_exist}.
\Moreover
\cref{lem:sum:ANN,item:lem:Power:constant_exist}
establish \cref{item:lem:Power:exist}.
\Moreover
\Enum{
\eg Grohs et al.\ \cite[item (i) in Proposition~2.6]{GrohsHornungJentzen2019}
;
\cref{item:sum:ANN:1} in \cref{lem:sum:ANN}
;
\cref{padding2_item1} in \cref{padding_lemma2}
}
that
\begin{equation}
\dims(\Psi)
=
\pr[\big]{ \inDimANN(\ii_1) , \textstyle\sum_{k=1}^\gamma \dimANNlevel_1(\ii_1) , \outDimANN(\ii_1) } 
= 
(1,\gamma,1)
\end{equation}
\cfload.
This establishes \cref{item:lem:Power:dims}.
\Moreover
\Enum{
\eg  Grohs et al.\ \cite[item (v) in Proposition~2.6]{GrohsHornungJentzen2019}
;
\cref{padding_item4} in \cref{padding_lemma}
;
\cref{item:sum:ANN:4} in \cref{lem:sum:ANN}
}
that for all
$x\in\R$
it holds that
\begin{equation}
\begin{split}
\pr[]{ \functionANN{\activation}(\Psi) } (x)
&
=
c_0 
+ 
\pr[\Big]{ \functionANN{\activation}\pr[\Big]{ \oSum_{i=1}^{\gamma} \pr[\big]{ \scalar{c_i}{ \pr[]{ \compANN{\ii_1}{\A_{1,b_i}} } } } } } \lrSpace (x)
\\
&
=
c_0 
+
\SmallSum{i=1}{\gamma} c_i \pr[]{ \functionANN{\activation}( \ii_1 ) } ( x + b_i ) 
=
c_0 
+
\SmallSum{i=1}{\gamma} c_i ( x + b_i )^\gamma
\end{split}
\end{equation}
\cfload.
\Enum{
This
;
the binomial theorem
}
that for all
$x\in\R$
it holds that
\begin{equation}\label{power_3_11}
\pr[]{ \functionANN{\activation}(\Psi) } (x)
=
c_0 
+
\SmallSum{i=1}{\gamma} c_i \br*{ \SmallSum{j=0}{\gamma} \textstyle\binom{\gamma}{j} x^{\gamma-j} (b_i)^{j} }
=
c_0 
+
\SmallSum{j=0}{\gamma} \textstyle\binom{\gamma}{j} \br*{ \SmallSum{i=1}{\gamma}  c_i (b_i)^{j} } x^{\gamma-j}
\dpp
\end{equation}
Combining
\Enum{
\cref{power_3_11}
and 
\cref{item:lem:Power:constant_exist}
\hence\unskip
}
that for all $x\in \R$ it holds that
\begin{equation}
\begin{split}
\pr[]{ \functionANN{\activation}(\Psi) } (x)
&
=
c_0 
+
\SmallSum{j=0}{\gamma} \textstyle\binom{\gamma}{j} \br*{ \SmallSum{i=1}{\gamma}  c_i (b_i)^{j} } x^{\gamma-j}
\\
&
=
c_0 
+
\SmallSum{j=0}{\gamma} \textstyle\binom{\gamma}{j} \br[\big]{ \1_{\{\gamma-1\}}(j) \, \gamma^{-1} - \1_{\{\gamma\}}(j) \, c_0 } x^{\gamma-j}
=
x
\dpp
\end{split}
\end{equation}
This establishes \cref{item:lem:Power:cont,item:lem:Power:real}.
\end{aproof}

\cfclear
\begin{athm}{lemma}{lem:ReluPower:identity}[Shallow RePU ANN representation for the one-dimensional identity function]
Let 
$\gamma\in\N\backslash\{1\}$, 
$b_1,b_2,\dots,b_\gamma \in \R$,
$\activation \in C(\R, \R)$
satisfy for all 
$x \in \R$
that 
$\activation(x) = (\max\{x, 0\})^\gamma$
and 
$b_1 < b_2 < \ldots < b_\gamma$
\cfload.
Then
\begin{enumerate}[label=(\roman{*})]
\item
\label{item:lem:ReluPower:identity}
it holds for all 
$x \in \R$ 
that 
$(\functionANN{\activation}(\I_\gamma))(x) = x^\gamma$,
\item 
\label{item:lem:ReluPower:constant_exist}
there exist unique 
$c_0,c_1,\dots,c_\gamma \in \R$ 
which satisfy for all 
$k \in \{0,1,\dots,\gamma\}$ 
that
$\1_{\{\gamma\}}(k) \, c_0 + \sum_{i=1}^\gamma c_i (b_i)^{k} = \1_{\{\gamma-1\}}(k)\,\gamma^{-1}$,
\item 
\label{item:lem:ReluPower:exist}
there exists a unique $\Psi \in \ANNs$ which satisfies
\begin{equation}
\Psi
= 
\compANN{\A_{1,c_0}}{\pr[\bigg]{ \OSum{i=1}{\gamma} \pr[\Big]{ \scalar{c_i}{ \pr[\big]{ \compANN{\I_\gamma}{\A_{1,b_i}} } } } } }
\dc
\end{equation}
\item
\label{item:lem:ReluPower:dims} 
it holds that 
$\dims(\Psi) = (1,2\gamma,1) \in \N^3$,
\item
\label{item:lem:ReluPower:cont} 
it holds that 
$ \functionANN{\activation}(\Psi) \in C(\R, \R)$, 
and
\item
\label{item:lem:ReluPower:real}
it holds for all 
$x \in \R$ 
that 
$
(\functionANN{\activation}(\Psi))(x) = x
$
\end{enumerate}
\cfout.
\end{athm}

\begin{aproof}
First,
\nobs that
\Enum{\cref{eq:def:id:1}} 
that for all 
$x \in \R$ 
it holds that
\begin{equation}
\begin{split}
(\functionANN{\activation}(\I_\gamma))(x) 
= 
\activation(x) + (-1)^\gamma \activation(-x) 
&
= 
\pr[\big]{ \max\{x, 0\} }^{\!\gamma} + (-1)^\gamma \pr[\big]{ \max\{-x, 0 \} }^{\!\gamma}
\\
&
=
\pr[\big]{ \max\{x, 0\} }^{\!\gamma} + \pr[\big]{ \min\{x, 0 \} }^{\!\gamma}
= 
x^\gamma
\end{split}
\end{equation}
\cfload.
This establishes \cref{item:lem:ReluPower:identity}.
\Moreover
\cref{item:lem:Power:constant_exist} in \cref{lem:Power:identity}
establishes \cref{item:lem:ReluPower:constant_exist}.
\Moreover
\cref{lem:sum:ANN,item:lem:ReluPower:constant_exist}
establish \cref{item:lem:ReluPower:exist}.
\Moreover
\Enum{
\cref{item:sum:ANN:1} in \cref{lem:sum:ANN}
;
\cref{item:lem:Relu:dims} in \cref{lem:ReluLeaky:identity}
}
that
\begin{equation}
\dims(\Psi)
=
\pr[\big]{ \inDimANN(\I_\gamma) , \textstyle\sum_{k=1}^\gamma \dimANNlevel_1(\I_\gamma) , \outDimANN(\I_\gamma) } 
= 
(1,2\gamma,1)
\end{equation}
\cfload.
This establishes \cref{item:lem:ReluPower:dims}.
\Moreover
\Enum{
\cref{item:lem:ReluPower:identity}
;
e.g., Grohs et al.\ \cite[item (v) in Proposition~2.6]{GrohsHornungJentzen2019}
;
\cref{item:sum:ANN:4} in \cref{lem:sum:ANN}
}
that for all
$x\in\R$
it holds that
\begin{equation}
\begin{split}
\pr[]{ \functionANN{\activation}(\Psi) } (x)
&
=
c_0 
+ 
\pr[\Big]{ \functionANN{\activation}\pr[\Big]{ \oSum_{i=1}^{\gamma} \pr[\big]{ \scalar{c_i}{ \pr[]{ \compANN{\I_\gamma}{\A_{1,b_i}} } } } } } \lrSpace (x)
\\
&
=
c_0 
+
\SmallSum{i=1}{\gamma} c_i \pr[]{ \functionANN{\activation}( \I_\gamma ) } ( x + b_i ) 
=
c_0 
+
\SmallSum{i=1}{\gamma} c_i ( x + b_i )^\gamma
\end{split}
\end{equation}
\cfload.
\Enum{
This
;
the binomial theorem
}
that for all
$x\in\R$
it holds that
\begin{equation}\label{power_3_17}
\pr[]{ \functionANN{\activation}(\Psi) } (x) 
=
c_0 
+
\SmallSum{i=1}{\gamma} c_i \br*{ \SmallSum{j=0}{\gamma} \textstyle\binom{\gamma}{j} x^{\gamma-j} (b_i)^{j} }
=
c_0 
+
\SmallSum{j=0}{\gamma} \textstyle\binom{\gamma}{j} \br*{ \SmallSum{i=1}{\gamma}  c_i (b_i)^{j} } x^{\gamma-j}
\dpp
\end{equation}
Combining
\Enum{
\cref{power_3_17}
and 
\cref{item:lem:ReluPower:constant_exist}
\hence\unskip
}
that for all
$x\in\R$
it holds that
\begin{equation}
\begin{split}
\pr[]{ \functionANN{\activation}(\Psi) } (x)
&
=
c_0 
+
\SmallSum{j=0}{\gamma} \textstyle\binom{\gamma}{j} \br*{ \SmallSum{i=1}{\gamma}  c_i (b_i)^{j} } x^{\gamma-j}
\\
&
=
c_0 
+
\SmallSum{j=0}{\gamma} \textstyle\binom{\gamma}{j} \br[\big]{ \1_{\{\gamma-1\}}(j) \, \gamma^{-1} - \1_{\{\gamma\}}(j) \, c_0 } x^{\gamma-j}
=
x
\dpp
\end{split}
\end{equation}
This establishes \cref{item:lem:ReluPower:cont,item:lem:ReluPower:real}.
\end{aproof}

\cfclear
\begin{athm}{lemma}{lem:SoftPlus:identity}[Shallow softplus ANN representation for the one-dimensional identity function]
Let 
$\activation \in C(\R,\R)$ 
satisfy for all 
$x \in \R$ 
that 
$\activation(x) = \ln(1 + \exp(x))$.
Then
\begin{enumerate}[label=(\roman{*})]
\item
\label{item:lem:SoftPlus:cont} 
it holds that 
$ \functionANN{\activation}(\I_1) \in C(\R, \R)$ 
and
\item
\label{item:lem:SoftPlus:real}
it holds for all 
$x \in \R$ 
that 
$
(\functionANN{\activation}(\I_1))(x) = x
$
\end{enumerate}
\cfout.
\end{athm}

\begin{aproof}
\Nobs
that
\Enum{
\cref{eq:def:id:1}
}
that for all
$x \in \R$
it holds that
\begin{equation}
\begin{split}
(\functionANN{\activation}(\I_1))(x)
=
\activation(x) - \activation(-x)
&
=
\ln\pr[\big]{ 1 + \exp(x) } - \ln\pr[\big]{ 1 + \exp(-x) }
\\
&
=
\ln\pr*{ \frac{ 1 + \exp(x) }{ 1 + \exp(-x) } }
=
\ln\pr[\big]{ \exp(x) } 
=
x
\end{split}
\end{equation}
\cfload.
This establishes \cref{item:lem:SoftPlus:cont,item:lem:SoftPlus:real}.
\end{aproof}

\subsection{ANN representations for MLP approximations}

\cfclear
\begin{athm}{proposition}{lemma15}
Let 
$\Theta = \bigcup_{n\in\N}\! \Z^n$, 
$d,M,\fd\in\N$, 
$T\in(0,\infty)$, 
$a \in C(\R,\R)$,
$\fJ,\F,\G\in \ANNs$ 
satisfy
$\dims(\fJ) = (1,\fd,1)$,
$\Ra(\fJ) = \operatorname{id}_\R$,
$\Ra(\F) \in C(\R,\R)$,
and
$\Ra(\G) \in C(\R^d,\R)$, 
for every 
$\theta \in \Theta$
let
$\cU^\theta \colon [0,T]\to [0,T]$
and
$W^\theta\colon[0,T]\to \R^d$
be functions,
for every 
$\theta \in \Theta$, $n \in \N_0$ 
let
$\mlp_{n}^\theta \colon [0,T] \times \R^d\to \R$  
satisfy for all 
$t\in[0,T]$, 
$x\in\R^d$ 
that
\begin{align}\label{def:mlp}
& 
\mlp_{n}^\theta(t,x) 
= 
\frac{\1_\N(n)}{M^n} \left[ \sum_{k=1}^{M^n} \pr[\big]{ \Ra(\G) } \lrSpace \pr[\big]{ x + W_{T-t}^{(\theta,0,-k)} } \right]
\\
& 
+ 
\sum_{i=0}^{n-1} \frac{(T-t)}{M^{n-i}}\left[ \sum_{k=1}^{M^{n-i}} \pr[\big]{ \pr[]{ \Ra(\F) \circ \mlp_{i}^{(\theta,i,k)} } - \1_{\N}(i) \pr[]{ \Ra(\F) \circ \mlp_{\max\{i-1,0\}}^{(\theta,-i,k)} } } \lrSpace \pr[\big]{ \cU_t^{(\theta,i,k)}, x + W_{\cU_t^{(\theta,i,k)}-t}^{(\theta,i,k)} } \right]
\dc 
\nonumber
\end{align}
and 
let 
$\U_{n,t}^\theta\in \ANNs$, $t\in[0,T]$, $n\in\Z$, $\theta\in\Theta$, 
satisfy for all 
$\theta \in \Theta$, 
$n\in\N$, 
$t\in[0,T]$ 
that 
$\U_{0,t}^\theta = ((0\ 0\ \dots\ 0),0) \in \R^{1\times d}\times \R^1$ 
and
\begin{align}
\U_{n,t}^\theta 
& 
= 
\left[ \OSum{k=1}{M^n} \pr[\Big]{  \scalar{\tfrac{1}{M^n}}{\pr[\big]{ \compANN{ \G }{ \A_{\idMatrix_d,W_{T-t}^{(\theta,0,-k)}} } } } } \right]
\nonumber
\\
& 
\quad 
\bSum_{\,\fJ} 
\left[ \BSum{i=0}{\fJ}{n-1} \br[\Bigg]{ \scalar{ \pr[\Big]{ \tfrac{(T-t)}{M^{n-i}} } }{ \pr[\bigg]{ \BSum{k=1}{\fJ}{M^{n-i}} \pr[\Big]{ \compANN{ \pr[\big]{ \compANN{\F}{\U^{(\theta,i,k)}_{i,\mathcal{U}_t^{(\theta,i,k)}} } } }{ \A_{\idMatrix_d, W_{\mathcal{U}_t^{(\theta,i,k)}-t}^{(\theta,i,k)}} } } } } } \right]
\\
& 
\quad 
\bSum_{\,\fJ} 
\left[ \BSum{i=0}{\fJ}{n-1} \br[\Bigg]{ \scalar{ \pr[\Big]{ \tfrac{(t-T) \, \1_{\N}(i)}{M^{n-i}} } }{ \pr[\bigg]{ \BSum{k=1}{\fJ}{M^{n-i}} \pr[\Big]{ \compANN{ \pr[\big]{ \compANN{\F}{\U^{(\theta,-i,k)}_{\max\{i-1,0\},\mathcal{U}_t^{(\theta,i,k)}} } } }{ \A_{\idMatrix_d, W_{\mathcal{U}_t^{(\theta,i,k)}-t}^{(\theta,i,k)}} } } } } } \right] 
\nonumber
\end{align}
\cfload. 
Then
\begin{enumerate}[label=(\roman *)]
\item 
\label{lemma15_item2}
it holds for all 
$\theta_1,\theta_2\in\Theta$, 
$n\in\N_0$, 
$t_1,t_2\in [0,T]$ 
that 
$\dims(\U_{n,t_1}^{\theta_1}) 
= 
\mathcal{D}(\U_{n,t_2}^{\theta_2})$,
\item 
\label{lemma15_item2b}
it holds for all 
$\theta\in\Theta$, 
$n\in\N_0$, 
$t \in [0,T]$ 
that 
$\Ra(\U_{n,t}^\theta) \in C(\R^d,\R)$,
\item 
\label{lemma15_item5}
it holds for all 
$\theta\in\Theta$, 
$n\in\N_0$, 
$t\in[0,T]$ 
that 
$\lengthANN(\U_{n,t}^\theta) \le  \max\{\fd,\lengthANN(\G)\} + n \hiddenLength(\F)$,
\item 
\label{lemma15_item6}
it holds for all 
$\theta\in\Theta$, 
$n\in\N_0$, 
$t\in[0,T]$ 
that 
\begin{equation}\label{update1}
\normmm{\dims(\U_{n,t}^\theta)} 
\le 
\max\cu[\big]{ \fd, \normmm{\dims(\F) }, \normmm{\dims(\G)} } (3M)^n 
\dc
\end{equation}
\item 
\label{lemma15_item7}
it holds for all 
$\theta\in\Theta$, 
$n\in\N_0$, 
$t\in[0,T]$, 
$x\in\R^d$ 
that 
$U_{n}^\theta(t,x) = ((\Ra(\U_{n,t}^\theta))(x)$,
and
\item
\label{lemma15_item8}
it holds for all 
$\theta\in\Theta$, 
$n\in\N_0$, 
$t\in[0,T]$ 
that 
\begin{equation}\label{update2}
\paramANN(\U_{n,t}^\theta) 
\le 
2 \pr[\big]{ \max\{\fd,\lengthANN(\G)\} + n \hiddenLength(\F) } \pr[\Big]{ \max\cu[\big]{ \fd , \normmm{\dims(\F) }, \normmm{\dims(\G)} } }^{\!2}  (3M)^{2n}
\end{equation}
\end{enumerate}
\cfout.
\end{athm}

\begin{aproof}
Throughout this proof let 
$\Phi_{n,t}^\theta \in \ANNs$, $\theta\in\Theta$, $n\in\N$, $t\in[0,T]$, 
satisfy for all 
$\theta\in\Theta$, 
$n\in\N$, 
$t\in[0,T]$ 
that
\begin{equation}\label{eq:psi1}
\Phi_{n,t}^\theta 
= 
\OSum{k=1}{M^n} \pr[\Big]{  \scalar{\tfrac{1}{M^n}}{\pr[\big]{ \compANN{ \G }{ \A_{\idMatrix_d,W_{T-t}^{(\theta,0,-k)}} } } } }
\dc
\end{equation}
let 
$\Psi_{n,i,t}^{\theta,j} \in \ANNs$, $\theta\in\Theta$, $j\in\{0,1\}$, $n\in\N$, $i\in\{0,1,\ldots,n-1\}$, $t\in[0,T]$, 
satisfy for all 
$\theta\in\Theta$, 
$j\in\{0,1\}$, 
$n\in\N$, 
$i\in\{0,1,\ldots,n-1\}$, 
$t\in[0,T]$ 
that
\begin{equation}\label{eq:psi2}
\Psi_{n,i,t}^{\theta,j} 
= 
\BSum{k=1}{\fJ}{M^{n-i}} \pr[\Big]{ \compANN{ \pr[\big]{ \compANN{\F}{\U^{(\theta,(-1)^j i,k)}_{\max\{i-j,0\},\mathcal{U}_t^{(\theta, i,k)}} } } }{ \A_{\idMatrix_d, W_{\mathcal{U}_t^{(\theta,i,k)}-t}^{(\theta,i,k)}} } } 
\dc
\end{equation}
let 
$\Xi_{n,t}^{\theta,j} \in \ANNs$, $\theta \in \Theta$, $j\in\{0,1\}$, $n\in\N$, $t\in[0,T]$, 
satisfy for all 
$\theta \in \Theta$, 
$j\in\{0,1\}$, 
$n\in\N$, 
$t\in[0,T]$ 
that
\begin{equation}\label{eq:xi1}
\Xi_{n,t}^{\theta,j} 
= 
\BSum{i=0}{\fJ}{n-1} \br[\Big]{ \scalar{ \pr[\Big]{ \tfrac{(-1)^j(T-t) \, \1_{\N}(i-j+1)}{M^{n-i}} } }{ \Psi_{n,i,t}^{\theta,j} } } 
\dc
\end{equation}
let
$L_i \in \N$, $i\in \Z$, satisfy for all $i\in \Z$ that
\begin{equation}\label{eq:l1}
L_i=\lengthANN\pr[\Big]{ \compANN{\F}{\U^{0}_{\max\{i,0\},0} } } 
\dc
\end{equation}
let $\fL_n\in \N$, $n\in \N_0$, satisfy for all $n\in \N_0$ that 
\begin{equation}\label{eq:l1b}
\fL_n=\max_{i\in \{-1,0,\ldots, n-1\}}L_i \dc 
\end{equation}
and let 
$\LL_{n} \in \N$, $n\in\N$
satisfy for all 
$n\in\N$
that
\begin{equation}\label{eq:l2}
\LL_{n}
= 
\max\cu*{ \lengthANN(\G), \fL_{n}, \fL_{n-1}  } 
\dpp
\end{equation}
We prove 
\cref{lemma15_item2,%
lemma15_item2b,%
lemma15_item5,%
lemma15_item6,%
lemma15_item7} 
by induction on $n\in\N_0$. 
For the base case $n=0$ 
\nobs 
that 
\Enum{
the fact for all 
$\theta\in\Theta$, 
$t\in[0,T]$ 
it holds that 
$\U_{0,t}^\theta = ((0\ 0\ \dots\ 0),0) \in \R^{1\times d}\times \R^1$ 
}
that for all 
$\theta_1,\theta_2\in\Theta$, 
$t_1,t_2\in[0,T]$ 
it holds that 
\begin{equation}\label{eq:dim_equal_pre}
\dims(\U_{0,t_1}^{\theta_1}) 
= 
(d,1)
=
\dims(\U_{0,t_2}^{\theta_2}) 
\end{equation}
and that 
$\Ra(\U_{0,t_1}^{\theta_1}) \in C(\R^d,\R)$. 
\Moreover
\Enum{
the fact that \cref{def:mlp} implies that for all 
$\theta\in\Theta$, 
$t\in[0,T]$, 
$x\in\R^d$ it holds that 
$U_0^\theta(t,x) = 0$ 
;
the fact for all 
$\theta\in\Theta$, 
$t\in[0,T]$ 
it holds that 
$\U_{0,t}^\theta = ((0\ 0\ \dots\ 0),0) \in \R^{1\times d}\times \R^1$ 
}
that for all 
$\theta\in\Theta$, 
$t\in[0,T]$, 
$x\in\R^d$ 
it holds that 
\begin{equation}\label{eq:dim_equal_pre1}
\lengthANN(\U_{0,t}^\theta) 
= 
1
\dc 
\qquad
\normmm{\dims(\U_{0,t}^\theta)} 
= 
d
\dc 
\qquad \text{and} \qquad
(\Ra(\U_{0,t}^\theta))(x) 
= 
U_0^\theta(t,x) 
\end{equation}
\cfload.
Combining 
\cref{eq:dim_equal_pre}, 
\cref{eq:dim_equal_pre1},  
and the fact that the assumption that $\Ra(\G) \in C(\R^d,\R)$ implies that
$\max\cu[]{ \fd , \normmm{\dims(\F) }, \normmm{\dims(\G)} } \ge d$
hence proves 
\cref{lemma15_item2,%
lemma15_item2b,%
lemma15_item5,%
lemma15_item6,%
lemma15_item7} 
in the base case $n=0$. 
For the induction step $\N_0 \ni (n-1) \induct n \in \N$ let $n\in\N$ and assume that 
\cref{lemma15_item2,%
lemma15_item2b,%
lemma15_item5,%
lemma15_item6,%
lemma15_item7} 
hold true for all $k\in\{0,1,\dots,n-1\}$. 
\Nobs that 
\Enum{
the hypothesis that for every 
$\theta\in\Theta$, 
$t\in[0,T]$ 
it holds that 
$W_t^\theta \in \R^d$
;
\cref{lem:sum:ANN}
(applied for every 
$\theta \in \Theta$, 
$t \in [0,T]$ 
with
\begin{equation}
\begin{split}
&  
\lbd \with 1
\dc 
\quad 
\ubd \with M^n
\dc 
\quad 
(h_k)_{k\in\{\lbd,\lbd+1,\ldots,\ubd\}} \with (M^{-n})_{k\in\{1,2,\ldots,M^n\}}
\dc 
\quad 
(t_k)_{k\in\{\lbd,\lbd+1,\ldots,\ubd\}} \with (1)_{k\in\{1,2,\ldots,M^n\}}
\dc
\\
&
\Psi \with \Phi_{n,t}^\theta
\dc
\quad
(\Phi_k)_{k\in\{\lbd,\lbd+1,\ldots,\ubd\}} \with (\G)_{k\in\{1,2,\ldots,M^n\}}
\dc 
\quad 
(B_k)_{k\in\{\lbd,\lbd+1,\ldots,\ubd\}} \with (W_{T-t}^{\theta,0,-k})_{k \in \{1,2,\ldots,M^n\}}
\end{split}
\end{equation}
in the notation of \cref{lem:sum:ANN}) 
}
that for all 
$\theta\in\Theta$, 
$t\in[0,T]$, 
$x\in\R^d$ 
it holds that
\begin{equation}\label{mlp1}
\dims\pr[\big]{ \Phi_{n,t}^\theta }
= \pr[\big]{ d , M^n \dimANNlevel_1(\G), M^n \dimANNlevel_2(\G), \ldots, M^n \dimANNlevel_{\lengthANN(\G)-1}(\G), 1 } 
\in 
\N^{\lengthANN(\G)+1}
\end{equation}
and
\begin{equation}\label{eq:3_10}
\pr[\big]{ \Ra(\Phi_{n,t}^\theta) } \lrSpace (x) 
= 
\frac{1}{M^n} \br*{ \SmallSum{k=1}{M^n} \pr[\big]{ \Ra(\G) } \lrSpace \pr[\big]{ x + W_{T-t}^{(\theta,0,-k)} } }
\dpp
\end{equation}
\Moreover 
\Enum{
the induction hypothesis 
;
\eg Grohs et al.\ \cite[item~(ii) in Proposition~2.6]{GrohsHornungJentzen2019}
}
that for all 
$\theta \in \Theta$, 
$j \in \{0,1\}$, 
$i \in \{0,1,\ldots,n-1\}$,
$k \in \{1,2,\ldots,M^{n-i}\}$,
$t \in [0,T]$
it holds that
\begin{equation}
\begin{split} 
\lengthANN\pr[\Big]{ \compANN{\F}{\U^{(\theta,(-1)^j i,k)}_{\max\{i-j,0\},\mathcal{U}_t^{(\theta, i,k)}} } } 
&
= 
{ \lengthANN\pr[]{ \F } + \lengthANN\pr[\Big]{ \U^{(\theta,(-1)^j i,k)}_{\max\{i-j,0\},\mathcal{U}_t^{(\theta, i,k)} } } - 1 } 
\\
& = 
 \lengthANN\pr[]{ \F } + \lengthANN\pr[\Big]{ \U^{0}_{\max\{i-j,0\},0 } } - 1 
 \\
 &=
 \lengthANN\pr[\Big]{ \compANN{\F}{\U^{0}_{\max\{i-j,0\},0}} } 
 =
 L_{i-j}
\dpp
\end{split}
\end{equation}
\Enum{
This
;
the induction hypothesis
;
the hypothesis that for all 
$\theta\in\Theta$, 
$t\in[0,T]$ 
it holds that 
$W_t^\theta \in \R^d$
; 
the hypothesis that for all 
$\theta\in\Theta$, 
$t\in[0,T]$ 
it holds that 
$\cU^\theta_t \in [0,T]$
;
\cref{lemma14a}
(applied for every 
$j \in \{0,1\}$, 
$i \in \{0,1,\ldots,n-1\}$, 
$\theta \in \Theta$, 
$t \in [0,T]$ 
with 
\begin{equation}
\begin{split}
& 
\lbd \with 1
\dc 
\quad 
\ubd \with M^{n-i}
\dc 
\quad 
\fJ \with \fJ
\dc 
\quad 
(B_k)_{k\in\{\lbd,\lbd+1,\ldots,\ubd\}} \with \pr[\big]{ W_{\cU_t^{(\theta,i,k)}-t}^{\theta,i,k} }_{k \in \{1,2,\ldots,M^{n-i}\}}
\dc 
\\
& 
(h_k)_{k\in\{\lbd,\lbd+1,\ldots,\ubd\}} \with (1)_{k\in\{1,2,\ldots,M^{n-i}\}}
\dc 
\quad 
L \with L_{i-j}
\dc 
\quad 
\Psi \with \Psi_{n,i,t}^{\theta,j} 
\dc 
\quad a \with a 
\\
&
(\Phi_k)_{k\in\{\lbd,\lbd+1,\ldots,\ubd\}} \with \pr[\big]{ \compANN{\F}{\U^{(\theta,(-1)^j i,k)}_{\max\{i-j,0\},\mathcal{U}_t^{(\theta, i,k)}} } }_{k\in\{1,2,\ldots,M^{n-i}\}}
\end{split}
\end{equation}
in the notation of \cref{lemma14a})
;
\eg Grohs et al.\ \cite[Proposition~2.6]{GrohsHornungJentzen2019}
}
that for all 
$\theta\in\Theta$, 
$j\in\{0,1\}$, 
$i\in\{0,1,\ldots,n-1\}$, 
$t\in[0,T]$, 
$x\in\R^d$ 
it holds that
\begin{align}\label{eq:3_15a}
&
\dims\pr[\big]{ \Psi_{n,i,t}^{\theta,j} }
\nonumber
\\
& 
= 
\biggl( d , \SmallSum{k=1}{M^{n-i}} \dimANNlevel_1\lrSpace\pr[\Big]{ \longerANN{L_{i-j}, \fJ}\pr[\big]{ \compANN{\F}{\U^{(\theta,(-1)^j i,k)}_{\max\{i-j,0\},\mathcal{U}_t^{(\theta, i,k)}} } } }, \SmallSum{k=1}{M^{n-i}} \dimANNlevel_2\lrSpace\pr[\Big]{ \longerANN{L_{i-j}, \fJ}\pr[\big]{ \compANN{\F}{\U^{(\theta,(-1)^j i,k)}_{\max\{i-j,0\},\mathcal{U}_t^{(\theta, i,k)}} } } }, 
\nonumber 
\\
& 
\qquad\qquad 
\dots, 
\SmallSum{k=1}{M^{n-i}} \dimANNlevel_{L_{i-j}-1}\lrSpace\pr[\Big]{ \longerANN{L_{i-j}, \fJ}\pr[\big]{ \compANN{\F}{\U^{(\theta,(-1)^j i,k)}_{\max\{i-j,0\},\mathcal{U}_t^{(\theta, i,k)}} } } } , 1 \biggr)
\\
& 
= 
\biggl( d , M^{n-i} \dimANNlevel_1\lrSpace \pr[\big]{\compANN{\F}{\U^{0}_{\max\{i-j,0\},0} } } , M^{n-i} \dimANNlevel_2\lrSpace \pr[\big]{\compANN{\F}{\U^{0}_{\max\{i-j,0\},0} } } , 
\nonumber 
\\
& 
\qquad\qquad 
\dots, 
M^{n-i} \dimANNlevel_{L_{i-j}-1}\lrSpace \pr[\big]{\compANN{\F}{\U^{0}_{\max\{i-j,0\},0} } }  , 1 \biggr) 
\in 
\N^{L_{i-j}+1}  
\nonumber
\end{align}
and
\begin{align}\label{eq:3_15}
\pr[\big]{ \Ra(\Psi_{n,i,t}^{\theta,j}) } \lrSpace (x) 
& 
= 
\SmallSum{k=1}{M^{n-i}} \pr[\Big]{ \Ra\pr[\big]{ \compANN{\F }{ \U_{\max\{i-j,0\}, \cU_t^{(\theta,i,k)}}^{(\theta, (-1)^j i, k )} } } } \lrSpace \pr[\big]{ x + W_{\cU_t^{(\theta,i,k)}-t}^{(\theta,i,k)} } 
\nonumber 
\\
& 
= 
\SmallSum{k=1}{M^{n-i}} \pr[\Big]{ \Ra(\F) \circ \Ra\pr[\big]{ \U_{\max\{i-j,0\}, \cU_t^{(\theta,i,k)}}^{(\theta, (-1)^j i, k )} } } \lrSpace \pr[\big]{ x + W_{\cU_t^{(\theta,i,k)}-t}^{(\theta,i,k)} } 
\\
& 
= 
\SmallSum{k=1}{M^{n-i}} \pr[\Big]{ \Ra(\F) \circ \mlp_{\max\{i-j,0\}}^{(\theta, (-1)^j i, k )} } \lrSpace \pr[\big]{ \cU_t^{(\theta,i,k)}, x + W_{\cU_t^{(\theta,i,k)}-t}^{(\theta,i,k)} }
\dpp 
\nonumber
\end{align}
\Enum{
This
;
\cref{eq:l1b}
}
that for all 
$\theta \in \Theta$, 
$n\in \N$,
$j \in \{0,1\}$, 
$t \in [0,T]$
it holds that
\begin{equation}\label{eq:equal_fL}
\begin{split}
\max_{i\in \{0,1,\ldots, n-1\}}\lengthANN \pr[\Big]{ \scalar{ \pr[\Big]{ \tfrac{(-1)^j(T-t) \, \1_{\N}(i-j+1)}{M^{n-i}} } }{ \Psi_{n,i,t}^{\theta,j} } } 
&=
\max_{i\in \{0,1,\ldots, n-1\}}\lengthANN \pr[\Big]{  \Psi_{n,i,t}^{\theta,j}  }\\
&=
\max_{i\in \{0,1,\ldots, n-1\}}L_{i-j}
=\fL_{n-j}
\dpp
\end{split}
\end{equation}
Combining 
\Enum{
\cref{eq:3_15a}, 
\cref{eq:3_15}, 
\cref{eq:equal_fL}
\eg Grohs et al.\ \cite[item~(i) in Proposition~2.6]{GrohsHornungJentzen2019}, 
and 
\cref{lemma14a}
(applied for every 
$j \in \{0,1\}$, 
$\theta \in \Theta$, 
$t \in [0,T]$ 
with 
\begin{align}
& 
\lbd \with 0
\dc
\quad 
\ubd \with n-1
\dc 
\quad 
\fJ \with \fJ
\dc 
\quad 
(\Phi_k)_{k\in\{\lbd,\lbd+1,\ldots,\ubd\}} \with \pr[\big]{ \Psi_{n,i,t}^{\theta,j} }_{i\in\{0,1,\ldots,n-1\}}, \quad L \with \fL_{n-j}
\dc 
\\
& 
(h_k)_{k\in\{\lbd,\lbd+1,\ldots,\ubd\}} \with \pr[\Big]{ \tfrac{(-1)^j(T-t) \, \1_{\N}(i-j+1)}{M^{n-i}} }_{i\in\{0,1,\ldots,n-1\}}
\dc 
\quad 
(B_k)_{k\in\{\lbd,\lbd+1,\ldots,\ubd\}} \with ((0,0,\ldots,0))_{k \in \{0,1,\ldots,n-1\}} 
\nonumber 
\end{align}
in the notation of \cref{lemma14a}) 
}
that for all 
$j \in \{0,1\}$, 
$\theta \in \Theta$, 
$t \in [0,T]$, 
$x\in\R^d$ 
it holds that
\begin{equation}\label{xi_dims1}
\begin{split}
\dims\pr[\big]{ \Xi_{n,t}^{\theta,j} }
& 
= 
\biggl( d , \SmallSum{i=0}{n-1} \dimANNlevel_1\pr[\big]{  \longerANN{\fL_{n-j}, \fJ}\pr[\big]{ \Psi_{n,i,t}^{\theta,j} } }, \SmallSum{i=0}{n-1} \dimANNlevel_2\pr[\big]{ \longerANN{\fL_{n-j}, \fJ}\pr[\big]{ \Psi_{n,i,t}^{\theta,j} } }, 
\\
& 
\qquad\qquad 
\dots, 
\SmallSum{i=0}{n-1} \dimANNlevel_{\fL_{n-j}-1}\pr[\big]{ \longerANN{\fL_{n-j}, \fJ}\pr[\big]{ \Psi_{n,i,t}^{\theta,j} } } , 1 \biggr)  
\\
& = \biggl( d , \SmallSum{i=0}{n-1} \dimANNlevel_1\pr[\big]{ \longerANN{\fL_{n-j},\fJ}\pr[\big]{ \Psi_{n,i,0}^{0,j} } }, \SmallSum{i=0}{n-1} \dimANNlevel_2\pr[\big]{ \longerANN{\fL_{n-j},\fJ}\pr[\big]{ \Psi_{n,i,0}^{0,j} } }, 
\\
& 
\qquad\qquad  
\dots, \SmallSum{i=0}{n-1} \dimANNlevel_{\fL_{n-j}-1}\pr[\big]{ \longerANN{\fL_{n-j},\fJ}\pr[\big]{\Psi_{n,i,0}^{0,j} } } , 1 \biggr) 
\in 
\N^{\fL_{n-j}+1} 
\end{split}
\end{equation}
and
\begin{align}\label{eq:3_18}
\pr[\big]{ \Ra(\Xi_{n,t}^{\theta,j}) } \lrSpace (x) 
& 
= 
\SmallSum{i=0}{n-1} \pr[\Big]{ \tfrac{(-1)^j(T-t) \, \1_{\N}(i-j+1)}{M^{n-i}} } \pr[\big]{ \Ra(\Psi_{n,i,t}^{\theta,j} ) } \lrSpace (x) 
\\
& 
= 
\SmallSum{i=0}{n-1} \pr[\Big]{ \tfrac{(-1)^j(T-t) \, \1_{\N}(i-j+1)}{M^{n-i}} } \left[ \SmallSum{k=1}{M^{n-i}} \pr[\big]{ \Ra(\F) \circ \mlp_{\max\{i-j,0\}}^{(\theta, (-1)^j i, k )} } \lrSpace \pr[\big]{ \cU_t^{(\theta,i,k)}, x + W_{\cU_t^{(\theta,i,k)}-t}^{(\theta,i,k)} } \right] 
\dpp
\nonumber
\end{align}
\Moreover 
\Enum{ 
\cref{mlp1}
;
\cref{xi_dims1}
;
\cref{lemma14a_1} in \cref{lemma14a} 
(applied for every 
$\theta \in \Theta$, 
$t\in[0,T]$ 
with 
$\lbd \with 1$, 
$\ubd \with 3$, 
$L \with \LL_n$, 
$\Phi_1 \with \Phi_{n,t}^\theta$, 
$\Phi_2 \with \Xi_{n,t}^{\theta,0}$, 
$\Phi_3 \with \Xi_{n,t}^{\theta,1}$, 
$\fJ \with \fJ$, 
$h_1 \with 1$, 
$h_2 \with 1$, 
$h_3 \with 1$, 
$B_1 \with 0$, 
$B_2 \with 0$, 
$B_3 \with 0$ 
in the notation of \cref{lemma14a})
;
\eg Grohs et al.\ \cite[item~(i) in Proposition 2.6]{GrohsHornungJentzen2019}
}
that for all 
$\theta \in \Theta$, 
$t\in[0,T]$ 
it holds that
\begin{align}\label{eq:3_23}
\dims( \U_{n,t}^\theta )
& 
= 
\dims\pr[\big]{ \Phi_{n,t}^\theta \ \bSum_{\,\fJ} \ \Xi_{n,t}^{\theta,0} \ \bSum_{\,\fJ} \ \Xi_{n,t}^{\theta,1} } 
\nonumber 
\\
& 
= 
\Bigl( d , \dimANNlevel_1\pr[\big]{ \longerANN{\LL_n,\fJ}\pr[\big]{ \Phi_{n,t}^\theta } } + \dimANNlevel_1\pr[\big]{ \longerANN{\LL_n,\fJ}\pr[\big]{ \Xi_{n,t}^{\theta,0} } } + \dimANNlevel_1\pr[\big]{ \longerANN{\LL_n,\fJ}\pr[\big]{ \Xi_{n,t}^{\theta,1} } } , 
\nonumber 
\\
& 
\qquad 
\dimANNlevel_2\pr[\big]{ \longerANN{\LL_n,\fJ}\pr[\big]{ \Phi_{n,t}^\theta } } + \dimANNlevel_2\pr[\big]{ \longerANN{\LL_n,\fJ}\pr[\big]{ \Xi_{n,t}^{\theta,0} } } + \dimANNlevel_2\pr[\big]{ \longerANN{\LL_n,\fJ}\pr[\big]{ \Xi_{n,t}^{\theta,1} } } , 
\nonumber 
\\
& 
\qquad 
\dots , \dimANNlevel_{\LL_n-1}\pr[\big]{ \longerANN{\LL_n,\fJ}\pr[\big]{ \Phi_{n,t}^\theta } } + \dimANNlevel_{\LL_n-1}\pr[\big]{ \longerANN{\LL_n,\fJ}\pr[\big]{ \Xi_{n,t}^{\theta,0} } } + \dimANNlevel_{\LL_n-1}\pr[\big]{ \longerANN{\LL_n,\fJ}\pr[\big]{ \Xi_{n,t}^{\theta,1} } } , 1 \Bigr) 
\\
& 
= 
\Bigl( d , \dimANNlevel_1\pr[\big]{ \longerANN{\LL_n,\fJ}\pr[\big]{ \Phi_{n,0}^0 } } + \dimANNlevel_1\pr[\big]{ \longerANN{\LL_n,\fJ}\pr[\big]{ \Xi_{n,0}^{0,0} } } + \dimANNlevel_1\pr[\big]{ \longerANN{\LL_n,\fJ}\pr[\big]{ \Xi_{n,0}^{0,1} } } , 
\nonumber 
\\
& 
\qquad 
\dimANNlevel_2\pr[\big]{ \longerANN{\LL_n,\fJ}\pr[\big]{ \Phi_{n,0}^0 } } + \dimANNlevel_2\pr[\big]{ \longerANN{\LL_n,\fJ}\pr[\big]{ \Xi_{n,0}^{0,0} } } + \dimANNlevel_2\pr[\big]{ \longerANN{\LL_n,\fJ}\pr[\big]{ \Xi_{n,0}^{0,1} } } , 
\nonumber 
\\
& 
\qquad 
\dots , \dimANNlevel_{\LL_n-1}\pr[\big]{ \longerANN{\LL_n,\fJ}\pr[\big]{ \Phi_{n,0}^0 } } + \dimANNlevel_{\LL_n-1}\pr[\big]{ \longerANN{\LL_n,\fJ}\pr[\big]{ \Xi_{n,0}^{0,0} } } + \dimANNlevel_{\LL_n-1}\pr[\big]{ \longerANN{\LL_n,\fJ}\pr[\big]{ \Xi_{n,0}^{0,1} } } , 1 \Bigr) 
\nonumber 
\in 
\N^{\LL_n + 1}
\dpp
\nonumber 
\end{align}
Moreover, it then holds for all $\theta\in\Theta$, $t\in [0,T]$ that $\Ra(\U_{n,t}^{\theta}) \in C(\R^d,\R)$.
\Moreover[Next]
\Enum{
\cref{eq:psi1}
;
\cref{eq:psi2}
;
\cref{eq:xi1}
;
\cref{eq:3_10}
;
\cref{eq:3_18}
;
\cref{lemma14a_4} in \cref{lemma14a} 
(applied for every 
$\theta \in \Theta$, 
$t\in[0,T]$ 
with 
$\lbd \with 1$, 
$\ubd \with 3$, 
$L \with \LL_n$, 
$\Phi_1 \with \Phi_{n,t}^\theta$, 
$\Phi_2 \with \Xi_{n,t}^{\theta,0}$, 
$\Phi_3 \with \Xi_{n,t}^{\theta,1}$, 
$\fJ \with \fJ$, 
$h_1 \with 1$, 
$h_2 \with 1$, 
$h_3 \with 1$, 
$B_1 \with 0$, 
$B_2 \with 0$, 
$B_3 \with 0$ 
in the notation of \cref{lemma14a}) 
}
that for all 
$\theta \in \Theta$, 
$t\in[0,T]$, 
$x\in\R^d$ 
it holds that
\begin{equation}\label{eq:3_25}
\begin{split}
\pr[\big]{ \Ra( \U_{n,t}^\theta ) } \lrSpace (x) 
& 
= 
\pr[\big]{ \Ra\pr[]{ \Phi_{n,t}^\theta \ \bSum_{\,\fJ} \ \Xi_{n,t}^{\theta,0} \ \bSum_{\,\fJ} \ \Xi_{n,t}^{\theta,1} } } \lrSpace (x) 
\\
& 
= 
\pr[\big]{ \Ra\pr[]{ \Phi_{n,t}^\theta } } \lrSpace (x) 
+ 
\pr[\big]{ \Ra\pr[]{ \Xi_{n,t}^{\theta,0} } } \lrSpace (x) 
+ 
\pr[\big]{ \Ra\pr[]{ \Xi_{n,t}^{\theta,1} } } \lrSpace (x) 
\\
& 
= 
\tfrac{1}{M^n} \br*{ \SmallSum{k=1}{M^n} \pr[\big]{ \Ra(\G) } \lrSpace \pr[\big]{ x + W_{T-t}^{(\theta,0,-k)} } } 
\\
& 
\quad 
+ 
\SmallSum{i=0}{n-1} \tfrac{(T-t)}{M^{n-i}} \br*{ \SmallSum{k=1}{M^{n-i}} \pr[\big]{ \Ra(\F) \circ \mlp_{i}^{(\theta, i, k )} } \lrSpace \pr[\big]{ \cU_t^{(\theta,i,k)}, x + W_{\cU_t^{(\theta,i,k)}-t}^{(\theta,i,k)} } } 
\\
& 
\quad 
+ 
\SmallSum{i=0}{n-1} \tfrac{(t-T) \, \1_{\N}(i) }{M^{n-i}} \br*{ \SmallSum{k=1}{M^{n-i}} \pr[\big]{ \Ra(\F) \circ \mlp_{\max\{i-1,0\}}^{(\theta, -i, k )} } \lrSpace \pr[\big]{ \cU_t^{(\theta,i,k)}, x + W_{\cU_t^{(\theta,i,k)}-t}^{(\theta,i,k)} } }
\\
&
= 
\mlp_n^\theta(t,x) 
\dpp
\end{split}
\end{equation}
\Moreover
\Enum{ 
\cref{eq:3_23}
;
Jensen's inequality
}
that for all 
$\theta \in \Theta$, 
$t\in[0,T]$ 
it holds that
\begin{align}
& \normmm{ \dims( \U_{n,t}^\theta ) }
= \normmm{ \dims( \U_{n,0}^0 ) } 
\nonumber 
\\
& 
= 
\max\cu*{ d , \max_{k \in \{1,2,\ldots,\LL_n-1\}} \pr[\Big]{ \dimANNlevel_{k}\pr[\big]{ \longerANN{\LL_n,\fJ}\pr[]{ \Phi_{n,0}^0 } } + \dimANNlevel_{k}\pr[\big]{ \longerANN{\LL_n,\fJ}\pr[]{ \Xi_{n,0}^{0,0} } } + \dimANNlevel_{k}\pr[\big]{ \longerANN{\LL_n,\fJ}\pr[]{ \Xi_{n,0}^{0,1} } } } } 
\nonumber 
\\
& 
\le 
\max_{k \in \{1,2,\ldots,\LL_n-1\}} 
\Bigl[ 
\max\cu*{ d , \dimANNlevel_{k}\pr[\big]{ \longerANN{\LL_n,\fJ}\pr[]{ \Phi_{n,0}^0 } } } 
+ 
\max\cu*{ d , \dimANNlevel_{k}\pr[\big]{ \longerANN{\LL_n,\fJ}\pr[]{ \Xi_{n,0}^{0,0} } } } 
\\
& 
\qquad  
+ 
\max\cu*{ d , \dimANNlevel_{k}\pr[\big]{ \longerANN{\LL_n,\fJ}\pr[]{ \Xi_{n,0}^{0,1} } } } 
\Bigr]
\nonumber 
\\
& 
\le 
\normmm{ \dims\pr[\big]{ \longerANN{\LL_n,\fJ}\pr[\big]{ \Phi_{n,0}^0 } } } 
+ 
\normmm{ \dims\pr[\big]{ \longerANN{\LL_n,\fJ}\pr[]{ \Xi_{n,0}^{0,0} } } } 
+ 
\normmm{ \dims\pr[\big]{ \longerANN{\LL_n,\fJ}\pr[]{ \Xi_{n,0}^{0,1} } } } 
\dpp
\nonumber 
\end{align}
\Enum{
This
;
\cref{mlp1}
;
\cref{eq:3_15a}
;
\cref{xi_dims1}
;
Jensen's inequality
;
\eg Grohs et al.\ \cite[item~(i) in Proposition~2.6]{GrohsHornungJentzen2019}
} 
that for all 
$\theta \in \Theta$, 
$t\in[0,T]$ 
it holds that
\begin{align}
& 
\normmm{ \dims( \U_{n,t}^\theta ) } 
\nonumber 
\\
& 
\le 
\max\cu*{ \fd , M^n \normmm{\dims(\G)} } 
+ 
 \SmallSum{i=0}{n-1} \max\cu*{ \fd , \normmm{ \dims\pr[\big]{ \Psi_{n,i,0}^{0,0} } } } 
+ 
\SmallSum{i=0}{n-1} \max\cu*{ \fd , \normmm{ \dims\pr[\big]{ \Psi_{n,i,0}^{0,1} } } } 
\nonumber 
\\
& 
\le 
\max\cu*{ \fd , M^n \normmm{\dims(\G)} } 
\nonumber 
\\
&
\quad
+ 
\SmallSum{i=0}{n-1} M^{n-i} \br[\Big]{ \max\cu*{ \fd , \normmm{ \dims\pr[\big]{ \compANN{\F}{\U^{0}_{i,0} } } } } + \max\cu*{ \fd , \normmm{ \dims\pr[\big]{ \compANN{\F}{\U^{0}_{\max\{i-1,0\},0} } } } } } 
\\
& 
\le 
\max\cu*{ \fd , M^n \normmm{\dims(\G)} } 
\nonumber
\\
&
\quad
+ \SmallSum{i=0}{n-1} M^{n-i} \br[\Big]{ \max\cu*{ \fd , \normmm{ \dims\pr[]{ \F } } , \normmm{ \dims\pr[\big]{ \U^{0}_{i,0} } } } + \max\cu*{ \fd , \normmm{ \dims\pr[]{ \F } } , \normmm{ \dims\pr[\big]{ \U^{0}_{\max\{i-1,0\},0} } } } } 
\nonumber 
\dpp
\end{align}
\Enum{
The induction hypothesis \hence\unskip
}
that for all 
$\theta \in \Theta$, 
$t\in[0,T]$ 
it holds that
\begin{equation}\label{eq:3_28}
\begin{split}
\normmm{ \dims( \U_{n,t}^\theta ) } 
&\le 
\max\cu*{ \fd ,  \normmm{\dims(\G)},  \normmm{\dims(\F)}} M^n
\\
&
\quad
+ \SmallSum{i=0}{n-1} M^{n-i} 
\max\cu*{ \fd ,  \normmm{\dims(\G)},  \normmm{\dims(\F)}}
\br[\Big]{(3M)^i+(3M)^{\max\{i-1,0\}}}\\
&=
\max\cu*{ \fd ,  \normmm{\dims(\G)},  \normmm{\dims(\F)}} M^n 
\br[\Big]{1+\pr[\Big]{\SmallSum{i=0}{n-1}3^i} +1+ \tfrac{1}{M}\SmallSum{i=1}{n-1}3^{i-1}
}
\\
&\le 
\max\cu*{ \fd ,  \normmm{\dims(\G)},  \normmm{\dims(\F)}} M^n 
\br[\Big]{2+\tfrac{3^n-1}{2}+\tfrac{3^{n-1}-1}{2}
}
\\
&= 
\max\cu*{ \fd ,  \normmm{\dims(\G)},  \normmm{\dims(\F)}} (3M)^n 
\br[\Big]{\tfrac{1}{3^n}+\tfrac{2}{3}
}
\\
&\le 
\max\cu*{ \fd ,  \normmm{\dims(\G)},  \normmm{\dims(\F)}} (3M)^n 
\dpp
\end{split}
\end{equation}
\Moreover
\Enum{
the induction hypothesis
;
\cref{eq:l2}
;
\cref{eq:3_23}
;
Jensen's inequality
;
\eg Grohs et al.\ \cite[item~(ii) in Proposition~2.6]{GrohsHornungJentzen2019}
}
that for all 
$\theta \in \Theta$, 
$t\in[0,T]$ 
it holds that
\begin{align}\label{eq:3_30}
\lengthANN(\U_{n,t}^\theta) 
& 
= 
\lengthANN(\U_{n,0}^0) 
= 
\LL_n 
= 
\max\cu*{ \lengthANN(\G) , \max_{i\in \{-1,0,\ldots,n-1\} } \lengthANN\pr[\big]{ \compANN{\F}{\U^{0}_{\max\{i,0\},0} } } } 
\nonumber
\\
& 
= 
\max\cu*{ \lengthANN(\G) , \max_{i \in \{-1,0,\ldots,n-1\} } \pr[\Big]{ \lengthANN(\F) + \lengthANN\pr[\big]{ \U_{\max\{i,0\},0}^0 } - 1} } 
\nonumber
\\
& 
= 
\max\cu*{ \lengthANN(\G) , \hiddenLength(\F) +  \max_{i \in \{-1,0,\ldots,n-1\}  } \lengthANN\pr[\big]{ \U_{\max\{i,0\},0}^0 } } 
\\
& 
\le 
\max\cu*{ \lengthANN(\G) , \hiddenLength(\F) +  \max_{i \in \{0,1,\ldots,n-1\}  } \pr[\Big]{ \max\{ \fd ,\lengthANN(\G)\} + \max\{i,0\} \hiddenLength(\F) } } 
\nonumber
\\
& 
= 
\max\cu[\Big]{ \lengthANN(\G) , \hiddenLength(\F) + \br[\big]{ \max\{\fd,\lengthANN(\G)\} + (n-1) \hiddenLength(\F) } } 
\le 
\max\{\fd,\lengthANN(\G)\} + n \hiddenLength(\F) 
\dpp
\nonumber
\end{align}
Combining 
\cref{eq:3_23,eq:3_25,eq:3_28,eq:3_30} 
completes the induction step.
Induction hence establishes
\cref{lemma15_item2,%
lemma15_item2b,%
lemma15_item5,%
lemma15_item6,%
lemma15_item7}.
\Moreover
\Enum{
\cref{lemma15_item5}
;
\cref{lemma15_item6} 
}
that for all 
$\theta \in \Theta$, 
$n\in\N_0$, 
$t\in[0,T]$ 
it holds that
\begin{equation}
\begin{split}
\paramANN(\U_{n,t}^\theta) 
& 
\le 
\SmallSum{k=1}{\lengthANN(\U_{n,t}^\theta)} \normmm{ \dims( \U_{n,t}^\theta ) } \br[\Big]{ \normmm{ \dims( \U_{n,t}^\theta ) } + 1 } 
\le 
2 \lengthANN(\U_{n,t}^\theta) \normmm{ \dims( \U_{n,t}^\theta ) }^2 
\\
& 
\le 
2 \pr[\big]{ \max\{\fd,\lengthANN(\G)\} + n \hiddenLength(\F) } \pr[\Big]{ \max\cu[\big]{ \fd , \normmm{\dims(\F) }, \normmm{\dims(\G)} } }^{\!2}  (3M)^{2n}
\dpp
\end{split}
\end{equation}
This establishes \cref{lemma15_item8}.
\end{aproof}

%
%
%


\section{ANN approximations for PDEs}\label{sec:4}

In this section we use the ANN representations for MLP approximations from~\cref{subsubsec:ANN_repres_for_MLP} to 
state and prove in \cref{theorem:final} in \cref{sec:4_1} below the main ANN approximation result of this work. 
In our proof of \cref{theorem:final} we employ the error estimates for suitable MLP approximations from Hutzenthaler et al.~\cite{PadgettJentzen2021} while the arguments in our proof of \cref{theorem:final} are inspired by the arguments in Hutzenthaler et al.~\cite{HutzenthalerJentzenKruse2019}. 

In the elementary results in \cref{sec:3.2} we show that the Lipschitz continuous nonlinearity of the PDE in~\cref{eq:mainThm:heateq} in \cref{theorem:final} can be approximated with suitable convergence rates by ANNs with ReLU, leaky ReLU, or softplus activation functions. 
In the situation of the ReLU activation function we use a linear interpolation technique similar as, for example, in Hutzenthaler et al.~\cite[Section 3.4]{HutzenthalerJentzenKruse2019}. 

In \cref{sec:4_2} below we combine the ANN approximation result in \cref{theorem:final} from \cref{sec:4_1} and the 
ANN approximation results for the Lipschitz continuous PDE nonlinearities from \cref{sec:3.2} with the ANN presentation results for the one-dimensional identity function $\R \ni x \mapsto x \in\R$ from \cref{subsec:ANNforId} to establish the ANN approximation results in \cref{cor:almostfinal} and \cref{cor:final}. 
\cref{thm:intro} in the introduction, in turn, is an immediate consequence of \cref{cor:final} from \cref{sec:4_2}.

\subsection{ANN approximation results with general activation functions}\label{sec:4_1}

\cfclear
\begin{athm}{theorem}{theorem:final}
Let 
$\LipConstF, \constantAssumpMainThm, \alpha_0, \alpha_1, \beta_0, \beta_1, T \in (0,\infty)$, 
$p,r\in \N$,
$\fq \in [2,\infty)$, 
$\activation \in C(\R,\R)$,
for every $d \in \N_0$ let
$\smallF_d \in C(\R^{\max\{d,1\}},\R)$,
for every 
$d \in \N$ 
let 
$\nu_d \colon \mathcal{B}(\R^{d}) \to [0,1]$ 
be a probability measure with 
\begin{equation}
\textstyle{\int_{\R^{d}} \norm{y}^{p^2\fq} \,\nu_{d}(\dxx y) \le \constantAssumpMainThm d^{rp^2\fq}}, 
\end{equation}
let $\fJ \in \ANNs$ satisfy 
$\hiddenLength(\fJ) = 1$
and 
$\Ra(\fJ) = \operatorname{id}_\R$,
let 
$(\interpolatingDNN_{d,\varepsilon})_{(d,\varepsilon)\in\N_0\times(0,1]} \subseteq \ANNs$ 
satisfy for all 
$d \in \N_0$, 
$x \in \R^{\max\{d,1\}}$, 
$\varepsilon \in (0,1]$ 
that 
\begin{equation}
\textstyle{\Ra(\interpolatingDNN_{d,\varepsilon}) \in C(\R^{\max\{d,1\}},\R)}, \quad 
\textstyle{\varepsilon^{\alpha_{\min\{d,1\}}} \lengthANN(\interpolatingDNN_{d,\varepsilon}) 
+
\varepsilon^{\beta_{\min\{d,1\}}} \normmm{\dims(\interpolatingDNN_{d,\varepsilon})} \le \constantAssumpMainThm (\max\{d,1\})^p},
\end{equation}
\begin{equation}
\text{and} \quad  
\textstyle{\varepsilon \abs{ (\Ra(\interpolatingDNN_{d,\varepsilon}))(x) } + \abs{ \smallF_d(x) - (\Ra(\interpolatingDNN_{d,\varepsilon}))(x) } \le \varepsilon \constantAssumpMainThm (\max\{d,1\})^p(1 + \norm{x} )^{p}},  
\end{equation}
and assume for all 
$v, w \in \R$, $\varepsilon \in (0,1]$ that 
\begin{equation}
\textstyle{\max\{ \abs{ \smallF_0(v) - \smallF_0(w) } , \abs{ (\Ra(\interpolatingDNN_{0,\varepsilon}))(v) - (\Ra(\interpolatingDNN_{0,\varepsilon}))(w) } \} \le \LipConstF \abs{v-w}}
\end{equation}
\cfload.
Then
\begin{enumerate}[label=(\roman *)]
\item 
\label{final_item1}
for every
$d \in \N$ 
there exists a unique at most polynomially growing viscosity solution
$\smallU_d \in C([0,T]\times\R^d,\R)$
of
\begin{equation}\label{eq:mainThm:heateq}
\pr[]{\tfrac{\partial}{\partial t} \smallU_d}(t,x)
+
\tfrac{1}{2} \pr[]{ \Delta_x \smallU_d }(t,x)
+
\smallF_0\pr[\big]{ \smallU_d(t,x) }
=
0
\end{equation}
with $\smallU_d(T,x) = \smallF_d(x)$ for $(t,x) \in (0,T) \times \R^d$
and 
\item 
\label{final_item2}
there exist 
$(\U_{d,\varepsilon})_{(d,\varepsilon)\in\N\times(0,1]} \subseteq \ANNs$
and 
$\eta \colon (0,\infty) \to \R$ 
such that for all 
$d \in \N$, 
$\varepsilon \in (0,1]$,
$\delta \in (0,\infty)$
it holds that 
\begin{equation}
\begin{split}
\paramANN(\U_{d,\varepsilon}) 
& \le \eta(\delta) 
d^{ 3p + 2(2+\delta) (2rp + p +2) p  + (rp+p+1)p (\max\{\alpha_0,\alpha_1\} + 2 \max\{\beta_0,\beta_1\}) } \\
& \quad \cdot \varepsilon^{-(2(2+\delta) + \max\{\alpha_0,\alpha_1\} + 2 \max\{\beta_0,\beta_1\} ) },
\end{split}
\end{equation}
\begin{equation}\label{final_item2_error}
\Ra(\U_{d,\varepsilon}) \in C(\R^{d},\R), \quad 
\text{and} \quad 
	\textstyle{
	\bigl[ \int_{\R^d} \abs{ 
		\smallU_d(0,x) 
		- 
		\pr[]{ \Ra\pr[]{ \U_{d,\varepsilon} } } (x) 
	}^\fq  \nu_d(\dxx x)  \bigr]^{\nicefrac{1}{\fq}} 
	\le 
	\varepsilon }
\dpp
\end{equation}
\end{enumerate}
\end{athm}

\begin{aproof}
\newcommand{\fnl}{f}
\newcommand{\fterm}{f}
\newcommand{\fnlnn}{\interpolatingDNN}
\newcommand{\ftermnn}{\mathbf{F}}
Throughout this proof 
let 
$\fB \in[1,\infty)$, 
$(\littleMM_k)_{k\in\N} \subseteq \N$ 
satisfy for all 
$k \in \N$ 
that 
$\liminf_{j\to\infty} \littleMM_j = \infty$, 
$\limsup_{j\to\infty} \nicefrac{(\littleMM_j)^{\fq/2}}{j} < \infty$,
and
$\littleMM_{k+1} \le \fB \littleMM_k$,
let $\fd\in \N$
satisfy
$\dims(\fJ) = (1,\fd,1)$,
let $\Theta = \bigcup_{n\in\N} \! \Z^n$,
let $(\Omega, \cF, \P)$ be a probability space,  
let 
$\fu^\theta\colon \Omega \to [0,1]$, $\theta \in \Theta$, 
be i.i.d.\ random variables, 
assume for all $t\in (0,1)$ that $\P(\fu^0\le t)=t$,
let 
$\cU^\theta \colon [0,T] \times \Omega \to [0,T]$, $\theta \in \Theta$, 
satisfy for all 
$t \in [0,T]$, 
$\theta \in \Theta$ 
that 
$\cU_t^\theta = t + (T-t)\fu^\theta$, 
let 
$W^{d,\theta} \colon [0,T] \times  \Omega \to \R^d$, $d \in \N$, $\theta \in \Theta$, 
be independent standard Brownian motions,
assume for every 
$d \in \N$ 
that 
$(\cU^\theta)_{\theta\in\Theta}$ 
and 
$(W^{d,\theta})_{\theta\in\Theta}$ 
are independent,
let 
$\mlp_{n,j,\varepsilon}^{d,\theta} \colon [0,T] \times \R^d \times \Omega \to \R$, 
$d, j, n \in \Z$, 
$\theta \in \Theta$, 
$\varepsilon \in (0,1]$, 
satisfy for all 
$\varepsilon \in (0,1]$, 
$n \in \N_0$, 
$d, j \in \N$, 
$\theta \in \Theta$, 
$t \in [0,T]$, 
$x \in \R^d$ 
that 
\begin{align}
\mlp_{n,j,\varepsilon}^{d,\theta}(t,x) 
& 
= 
\frac{\1_\N(n)}{(\littleMM_{j})^n} 
\br[\Bigg]{ \SmallSum{k=1}{(\littleMM_{j})^n} 
\pr[\big]{ \Ra(\ftermnn_{d,\varepsilon})} \lrSpace \pr[\big]{ x + W_{T-t}^{d,(\theta,0,-k)} } } 
\nonumber
\\
& 
\quad 
+ 
\sum_{i=0}^{n-1} 
\frac{(T-t)}{(\littleMM_{j})^{n-i}} 
\Biggl[ \SmallSum{k=1}{(\littleMM_{j})^{n-i}} 
\Bigl[ \pr[\big]{(\Ra(\fnlnn_{0,\varepsilon})} \lrSpace \pr[\big]{ \mlp_{i,j,\varepsilon}^{d,(\theta,i,k)} \pr[\big]{ \cU_t^{(\theta,i,k)}, x + W_{\cU_t^{(\theta,i,k)}-t}^{d,(\theta,i,k)} } } 
\\
& 
\quad 
- 
\1_{\N}(i) 
\pr[\big]{(\Ra(\fnlnn_{0,\varepsilon})} \lrSpace \pr[\big]{ \mlp_{\max\{i-1,0\},j,\varepsilon}^{d,(\theta,-i,k)} \pr[\big]{ \cU_t^{(\theta,i,k)}, x + W_{\cU_t^{(\theta,i,k)}-t}^{d,(\theta,i,k)} } } \Bigr] \Biggr]
\dc
\nonumber
\end{align}
let 
$\U_{n,j,t}^{d,\theta,\varepsilon}\colon \Omega \to \ANNs$, 
$d, j, n \in \Z$, 
$\theta \in \Theta$, 
$t \in [0,T]$, 
$\varepsilon \in (0,1]$, 
satisfy for all 
$\varepsilon \in (0,1]$, 
$\theta \in \Theta$, 
$d, j, n \in \N$, 
$t \in [0,T]$, 
$\omega\in \Omega$ 
that 
$\U_{0,j,t}^{d,\theta,\varepsilon}(\omega) = ((0\ 0\ \dots\ 0),0) \in \R^{1\times d}\times \R^1$ 
and
\begin{align}\label{mlp_final_th_rep}
\U_{n,j,t}^{d,\theta,\varepsilon} 
& 
= 
\left[ 
\OSum{k=1}{(\littleMM_j)^n} 
\pr[\Big]{  \scalar{\tfrac{1}{(\littleMM_{j})^n}}{\pr[\big]{ \compANN{ \ftermnn_{d,\varepsilon} }{ \A_{\idMatrix_d,W_{T-t}^{d,(\theta,0,-k)}} } } } } 
\right]
\nonumber
\\
& 
\quad 
\bSum_{\,\fJ} 
\left[ 
\BSum{i=0}{\fJ}{n-1} 
\br[\Bigg]{ \scalar{ \pr[\Big]{ \tfrac{(T-t)}{(\littleMM_{j})^{n-i}} } }{ \pr[\bigg]{ \BSum{k=1}{\fJ}{(\littleMM_{j})^{n-i}} \pr[\Big]{ \compANN{ \pr[\big]{ \compANN{\fnlnn_{0,\varepsilon}}{\U^{d,(\theta,i,k),\varepsilon}_{i,j,\mathcal{U}_t^{(\theta,i,k)}} } } }{ \A_{\idMatrix_d, W_{\mathcal{U}_t^{(\theta,i,k)}-t}^{d,(\theta,i,k)}} } } } } } 
\right]
\\
& 
\quad 
\bSum_{\,\fJ} 
\left[ 
\BSum{i=0}{\fJ}{n-1} 
\br[\Bigg]{ \scalar{ \pr[\Big]{ \tfrac{(t-T) \, \1_{\N}(i)}{(\littleMM_{j})^{n-i}} } }{ \pr[\bigg]{ \BSum{k=1}{\fJ}{(\littleMM_{j})^{n-i}} \pr[\Big]{ \compANN{ \pr[\big]{ \compANN{\fnlnn_{0,\varepsilon}}{\U^{d,(\theta,-i,k),\varepsilon}_{\max\{i-1,0\},j,\mathcal{U}_t^{(\theta,i,k)}} } } }{ \A_{\idMatrix_d, W_{\mathcal{U}_t^{(\theta,i,k)}-t}^{d,(\theta,i,k)}} } } } } } 
\right] 
\nonumber
\end{align}
(cf.\ \cref{lemma15}),
assume without loss of generality that
$\max\{ \abs{\fnl_0(0)} , \fd, 1 \}
\le \constantAssumpMainThm$, 
let 
$c_d, \fc_d \in [1,\infty)$, $d \in \N$, 
satisfy for all 
$d \in \N$ 
that
\begin{equation}
\begin{split}
c_d 
& 
=  
\constantAssumpMainThm 2^{p-1}d^p
\pr[\big]{ e^{\LipConstF T}(T+1) }^{p+1}
((2\constantAssumpMainThm d^p)^{p}+1)
\\
& 
\quad 
\cdot
\left( 1+ 
\left(
\int_{\R^d} 
\norm{x}^{p^2 \fq}
 \, \nu_d(\dxx x) \right)^{\!\!\nicefrac{1}{p^2\fq}}
 +\pr[\big]{ \E\br[\big]{ \norm{ \fwpr_T^{d,0} }^{p^2} } }^{\!\!\nicefrac{1}{p^2}}
\right)^{p^2} 
\end{split}
\end{equation}
and 
\begin{equation}\label{eq:deffcd}
\begin{split}
\fc_d 
& 
=  
\constantAssumpMainThm 2^{2(p+1)}d^p
e^{\LipConstF T}(T+1) 
(\sqrt{\fq-1})
\left( 1+ 
\left(\int_{\R^d} 
\norm{x}^{p \fq}
\, \nu_d(\dxx x) \right)^\frac{1}{\fq}
+ \left(\sup_{s\in [0,T]} \E\br[\big]{ \norm{ \fwpr_s^{d,0} }^{p\fq} } 
\right)^{\frac{1}{\fq}} \right) 
\dc
\end{split}
\end{equation}
let $\fR\colon \N \times (0,1] \to [1,\infty]$ satisfy for all $d\in\N$, $\delta\in(0,1]$ that
\begin{equation}\label{eq:deffR}
\begin{split}
\fR(d,\delta) 
& = 
\inf\left( \left\{ n \in \N \colon \fc_d \left( \frac{(1+2LT)\exp\left(\frac{(\littleMM_n)^\frac{\fq}{2}}{n} \right)}{(\littleMM_n)^\frac{1}{2}} \right)^n \le \delta \right\} \cup \{\infty\} \right)
\dc
\end{split}
\end{equation}
let $\kappa_{\delta}$, $\delta \in (0,\infty)$, satisfy for all $\delta \in (0,\infty)$ 
that 
$\kappa_{\delta} = 2(2+\delta)+\max\{\alpha_0,\alpha_1\} + 2 \max\{\beta_0,\beta_1\}$, 
and 
let 
$\delta_{d,\varepsilon} \in (0,1]$, $d \in \N$, $\varepsilon \in (0,1]$, 
satisfy for all 
$d \in \N$, 
$\varepsilon \in (0,1]$ 
that 
$\delta_{d,\varepsilon} 
= 
\nicefrac{\varepsilon}{(c_d + 1)}$
\cfload.
\Nobs that 
\Enum{
the assumption that for all 
$w, z \in \R$, $\varepsilon \in (0,1]$
it holds that 
$\abs{ (\Ra(\fnlnn_{0,\varepsilon}))(z) - (\Ra(\interpolatingDNN_{0,\varepsilon}))(w) } \le \LipConstF \abs{z-w}$
;
the assumption that for all 
$d\in \N_0$, $x \in \R^{\max\{d,1\}}$, $\varepsilon \in (0,1]$  
it holds that 
$ \abs{ (\Ra(\ftermnn_{d,\varepsilon}))(x) }\le \constantAssumpMainThm (\max\{d,1\})^p(1 + \norm{x} )^{p}
$
;
Beck et al.\ \cite[Corollary~3.10]{MR4259658}
(applied for every 
$d \in \N$, $\varepsilon\in (0,1]$
with
$d \with d$,
$m \with d$,
$L \with \LipConstF$,
$T \with T$,
$\mu \with (\R^d \ni x \mapsto (0,0,\dots,0) \in \R^d)$,
$\sigma \with \idMatrix_d$,
$f \with ([0,T]\times\R^d\times\R \ni (t,x,w) \mapsto (\Ra(\fnlnn_{0,\varepsilon}))(w) \in \R)$,
$g \with \Ra(\ftermnn_{d,\varepsilon})$,
$W \with \fwpr^{d,0}$
in the notation of 
Beck et al.\ \cite[Corollary~3.10]{MR4259658})
}
that for every $d \in \N$, $\varepsilon \in (0,1]$ there exists a unique at most polynomially growing
$v_{d,\varepsilon} \in C([0,T]\times \R^d, \R)$ such that
for all
$t \in [0,T]$, 
$x \in \R^d$ 
it holds that
\begin{equation}
v_{d,\varepsilon}(t,x) 
= 
\E\br[\Big]{ \pr[\big]{\Ra(\ftermnn_{d,\varepsilon})} \lrSpace (x+W^{d,0}_{T-t}) } 
+ 
\int_t^T \E\br[\Big]{ \pr[\big]{\Ra(\interpolatingDNN_{0,\varepsilon})} \lrSpace (v_{d,\varepsilon}(s,x+W^{d,0}_{s-t})) } \dx s
\dpp
\end{equation}
\Moreover
\Enum{
the triangle inequality
;
the assumption that for all 
$d \in \N$, 
$x \in \R^d$, 
$\varepsilon \in (0,1]$ 
it holds that 
$\varepsilon \abs{(\Ra(\ftermnn_{d,\varepsilon}))(x)} + \abs{\fterm_d(x) - (\Ra(\ftermnn_{d,\varepsilon}))(x)} \le \varepsilon \constantAssumpMainThm d^p (1 + \norm{x})^{p}$
}
that for all 
$d \in \N$, 
$x \in \R^d$ 
it holds that
\begin{equation}\label{135}
\abs{\fterm_d(x)} 
\le 
\abs{\fterm_d(x) - (\Ra(\ftermnn_{d,1}))(x)} + \abs{(\Ra(\ftermnn_{d,1}))(x)} 
\le 
\constantAssumpMainThm d^p(1 + \norm{x})^{p}
\dpp
\end{equation}
Combining 
this, 
the assumption that for all 
$w, z \in \R$ 
it holds that 
$\abs{\fnl_0(w) - \fnl_0(z)} \le \LipConstF \abs{w-z}$,
and
Beck et al.\ \cite[Theorem~1.1]{beck2021nonlinear}
(applied for every 
$d \in \N$
with
$d \with d$,
$m \with d$,
$L \with \LipConstF$,
$T \with T$,
$\mu \with (\R^d \ni x \mapsto (0,0,\dots,0) \in \R^d)$,
$\sigma \with \idMatrix_d$,
$f \with (\R^d\times\R \ni (x,w) \mapsto \fnl_0(w) \in \R)$,
$g \with \fterm_d$,
$W \with \fwpr^{d,0}$
in the notation of 
Beck et al.\ \cite[Theorem~1.1]{beck2021nonlinear}) 
establishes 
\cref{final_item1}.
\Moreover
\Enum{
the fact that for all 
$d \in \N$, $s\in(0,T]$ 
the random variable 
$\norm{ \nicefrac{\fwpr_s^{d,0}}{\sqrt{s}} }^2$ 
is a chi-squared distributed random variable with $d$-degrees of freedom
;
Jensen's inequality
;
\eg
Simon \cite[Eq.~(2.35)]{simon2007probability}
}
that for all 
$\gamma, d \in \N$, $s\in[0,T]$
it holds that
\begin{equation}\label{eq:estimateBrownianExpect}
\pr[\Big]{ \E\br[\big]{ \norm{ \fwpr_s^{d,0} }^{\gamma} } }^{\!2} 
\le 
\E\br[\big]{ \norm{ \fwpr_s^{d,0} }^{2\gamma} } 
= 
(2s)^{\gamma} \br*{ \frac{\Gamma(\frac{d}{2} + \gamma)}{\Gamma(\frac{d}{2})} } 
= 
(2s)^{\gamma} \br*{ \prod_{k=0}^{\gamma-1} \left(\tfrac{d}{2} + k\right) }
\dpp
\end{equation}
\Enum{
This
}
that for all 
$d \in \N$ 
it holds that
\begin{equation}
\pr[\Big]{ \E\br[\big]{ \norm{ \fwpr_T^{d,0} }^{p^2} } }^{\!\!\nicefrac{1}{(p^2)}} 
\le 
\sqrt{2T} \br*{ \prod_{k=0}^{p^2-1} \left(\tfrac{d}{2} + k\right) }^{\!\nicefrac{1}{(2p^2)}} 
\le 
\sqrt{2T \pr*{ \tfrac{d}{2} + p^2 - 1 } }
\dpp
\end{equation}
Combining 
\Enum{
this 
and 
the assumption that for all 
$d \in \N$ 
it holds that 
$\int_{\R^{d}} \norm{x}^{p^2\fq} \, \nu_d(\dxx x) \le \constantAssumpMainThm d^{rp^2\fq}$ 
}
that there exists
$\overline{C} \in [1,\infty)$ 
such that for all 
$d \in \N$ 
it holds that
\begin{equation}\label{eq:boundforcd}
c_d 
\le 
\overline{C} d^{p+(r+1)p^2}
\dpp
\end{equation}
\Moreover
\Enum{
the triangle inequality
}
that for all 
$n \in \N_0$, 
$d \in \N$, 
$\delta \in (0,1]$ 
it holds that
\begin{align}\label{big_item1}
& 
\left( \int_{\R^d} \E\br[\Big]{ \abs[\big]{ \smallU_d(0,x) - \mlp_{n,j,\delta}^{d,0}(0,x) }^\fq } \, \nu_d(\dxx x) \right)^{\!\!\nicefrac{1}{\fq}} 
\\
& 
\le 
\left( \int_{\R^d} \abs[\big]{ \smallU_d(0,x) - v_{d,\delta}(0,x) }^\fq \, \nu_d(\dxx x) \right)^{\!\!\nicefrac{1}{\fq}} 
+ 
\left( \int_{\R^d} \E\br[\Big]{ \abs[\big]{ v_{d,\delta}(0,x) - \mlp_{n,j,\delta}^{d,0}(0,x) }^\fq } \, \nu_d(\dxx x) \right)^{\!\!\nicefrac{1}{\fq}} 
\dpp
\nonumber
\end{align}
\Moreover
\Enum{
the assumption that for all
$x \in \R$,
$\delta \in (0,1]$
it holds that
$
 \abs{ \fnl_0(x) - (\Ra(\fnlnn_{0,\delta}))(x) } \le \delta \constantAssumpMainThm (1 + \abs{x} )^{p}
$
}
that for all 
$\delta \in (0,1]$ 
it holds that
\begin{equation}\label{137}
\abs{ (\Ra(\fnlnn_{0,\delta}))(0) } 
\le 
\abs{ (\Ra(\fnlnn_{0,\delta}))(0) - \fnl_0(0) } 
+ 
\abs{\fnl_0(0)} 
\le 
\delta \constantAssumpMainThm
+ 
\abs{\fnl_0(0)} 
\le 
\constantAssumpMainThm
+ 
\abs{\fnl_0(0)} 
\le 2 \constantAssumpMainThm
\dpp
\end{equation}
\Moreover
\Enum{
the assumption that for all 
$d \in \N$,
$x \in \R^d$,
$w \in \R$, 
$\delta \in (0,1]$ 
it holds that 
$
 \abs{ \fnl_0(w) - (\Ra(\fnlnn_{0,\delta}))(w) } \le \delta \constantAssumpMainThm (1 + \abs{w} )^{p}
$
and
$\abs{ \fterm_d(x) - (\Ra(\ftermnn_{d,\delta}))(x) } \le \delta \constantAssumpMainThm d^p (1 + \norm{x} )^{p}$
}
that for all 
$d \in \N$, 
$w \in \R$, 
$x \in \R^d$, 
$\delta \in (0,1]$ 
it holds that
\begin{align}\label{zero_cond_1}
& 
\max\cu[\big]{ 
\abs{ \fnl_0(w) - (\Ra(\fnlnn_{0,\delta}))(w) } , 
\abs{ \fterm_d(x) - (\Ra(\ftermnn_{d,\delta}))(x) } 
} 
\\
& 
\le 
\max\cu[\big]{ 
\delta \constantAssumpMainThm (1 + \abs{w} )^{p}, 
\delta \constantAssumpMainThm d^p (1 + \norm{x} )^{p} 
} 
\le \delta \constantAssumpMainThm d^p
\max\cu[\big]{ 
 (1 + \norm{x}+\abs{w} )^{p}, 
(1 + \norm{x} )^{p} 
}
\nonumber
\\
&
\le \delta \constantAssumpMainThm d^p
 (1 + \norm{x}+\abs{w} )^{p}
 \le \delta \constantAssumpMainThm d^p2^{p-1}
 ((1 + \norm{x})^{p}+\abs{w}^{p})
 \le \delta \constantAssumpMainThm d^p2^{p-1}
 ((1 + \norm{x})^{p^2}+\abs{w}^{p})
 \nonumber
\dpp
\end{align} 
Combining 
\Enum{
this
;
\cref{135}
;
\cref{137}
;
the assumption that for all 
$d \in \N$, 
$w, z \in \R$, 
$x \in \R^d$, 
$\delta \in (0,1]$ 
it holds that
$\max\{\abs{ \fnl_0(w)-\fnl_0(z) } , \abs{ (\Ra(\fnlnn_{0,\delta}))(w) - (\Ra(\fnlnn_{0,\delta}))(z) } \} \le \LipConstF \abs{ w-z }$ 
and 
$\delta \abs{ (\Ra(\ftermnn_{d,\delta}))(x) } + \abs{ \fterm_d(x) - (\Ra(\ftermnn_{d,\delta}))(x) } \le \delta \constantAssumpMainThm d^p (1 + \norm{x} )^{p}$
;
Hutzenthaler et al.\ \cite[Lemma 2.3]{HutzenthalerJentzenKruse2019}
(applied for every 
$d \in \N$, 
$\delta \in (0,1]$ 
with 
$f_1 \with \fnl_0$, 
$f_2 \with \Ra(\fnlnn_{0,\delta})$, 
$g_1 \with \fterm_d$, 
$g_2 \with \Ra(\ftermnn_{d,\delta})$, 
$T \with T$, 
$L \with \LipConstF$, 
$B \with 2 \constantAssumpMainThm d^p$, 
$\delta \with \delta \constantAssumpMainThm 2^{p-1}d^p$, 
${\bf W} \with \fwpr^{d,0}$, 
$u_1 \with \smallU_d$, 
$u_2 \with v_{d,\delta}$, 
$p \with p$, 
$q \with p$ 
in the notation of Hutzenthaler et al.\ \cite[Lemma 2.3]{HutzenthalerJentzenKruse2019})
}
that for all 
$d \in \N$, 
$\delta \in (0,1]$ 
it holds that
\begin{align}\label{big_item2}
\left( 
\int_{\R^d} \abs{ \smallU_d(0,x) - v_{d,\delta}(0,x) }^\fq \, \nu_d(\dxx x)
\right)^{\!\!\nicefrac{1}{\fq}} 
& 
\le 
\delta \constantAssumpMainThm 2^{p-1}d^p
\pr[\big]{ e^{\LipConstF T}(T+1) }^{p+1}
((2\constantAssumpMainThm d^p)^{p}+1)
\\
& 
\quad 
\cdot
\left( 
\int_{\R^d} 
\pr[\Big]{ 1 + \norm{x} + \pr[\big]{ \E\br[\big]{ \norm{ \fwpr_T^{d,0} }^{p^2} } }^{\!\!\nicefrac{1}{p^2}} }^{p^2 \fq}
 \, \nu_d(\dxx x)
\right)^{\!\!\nicefrac{1}{\fq}} 
\dpp 
\nonumber
\end{align}
Combining 
\Enum{
\cref{big_item2} 
and 
the triangle inequality
\hence\unskip
}
that for all 
$d \in \N$, 
$\delta \in (0,1]$ 
it holds that
\begin{equation}\label{big_item2a}
\left( \int_{\R^d} \abs{ \smallU_d(0,x) - v_{d,\delta}(0,x) }^\fq \, \nu_d(\dxx x) \right)^{\!\!\nicefrac{1}{\fq}} 
\le
c_d \delta
\dpp
\end{equation}
\Moreover
\Enum{
	\cref{137}
	;
	the assumption that for all 
	$d \in \N$, 
	$x \in \R^d$, 
	$\delta \in (0,1]$ 
	it holds that
	$\delta \abs{ (\Ra(\ftermnn_{d,\delta}))(x) } + \abs{ \fterm_d(x) - (\Ra(\ftermnn_{d,\delta}))(x) } \le \delta \constantAssumpMainThm d^p (1 + \norm{x} )^{p}$
}
that for all 
$d \in \N$, 
$x \in \R^d$, 
$\delta \in (0,1]$ 
it holds that
\begin{equation}
\begin{split}
\max\cu[\big]{ 
	\abs{ (\Ra(\fnlnn_{0,\delta}))(0) } , 
	\abs{ (\Ra(\ftermnn_{d,\delta}))(x) } 
}
& 
\le 
\max\cu[\big]{ 
	2\constantAssumpMainThm , 
	\constantAssumpMainThm d^p (1 + \norm{x} )^{p} 
} 
\le 
2\constantAssumpMainThm d^p (1 + \norm{x} )^{p} 
\dpp
\end{split} 
\end{equation}
This,
the assumption that for all $v,w\in \R$, $\delta\in (0,1]$ it holds that
	$
	\abs{ (\Ra(\interpolatingDNN_{0,\delta}))(v) - (\Ra(\interpolatingDNN_{0,\delta}))(w) } \le \LipConstF \abs{v-w}
	$,
and Hutzenthaler et al.\ \cite[Corollary 3.15]{PadgettJentzen2021}
	(applied 
	for every 
	$\delta \in (0,1]$, $d,j\in \N$
	with 
	$T \with T$, 
	$L\with \LipConstF$,
	$\mathfrak{L}\with 2^p \constantAssumpMainThm d^p$,
	$p \with p$, 
	$\mathfrak{p} \with \fq$,
	$m\with m_j$,
	$f \with ( [0,T]\times \R^d \times \R \ni (t,x,v)\mapsto (\Ra(\fnlnn_{0,\delta}))(v)\in \R)$,
	$g \with \Ra(\ftermnn_{d,\delta})$,
	$\Theta \with \Theta$, 
	$(\fu^\theta)_{\theta\in\Theta} \with (\fu^\theta)_{\theta\in\Theta}$, 
	$(\cU^\theta)_{\theta\in\Theta} \with (\cU^\theta)_{\theta\in\Theta}$, 
	$(W^{\theta})_{\theta\in \Theta} \with (W^{d,\theta})_{\theta\in\Theta}$, 
	$u \with v_{d,\delta}$, 
	$(U_{n}^{\theta})_{(n,\theta) \in \Z\times\Theta} \with (\mlp_{n,j,\delta}^{d,\theta})_{(n,\theta) \in \Z\times\Theta}$
	in the notation of Hutzenthaler et al.\ \cite[Corollary 3.15]{PadgettJentzen2021})
assure 
that for all $\delta\in (0,1]$, $d,j\in\N$, $n\in\N_0$, $x \in \R^d$ it holds that 
\begin{equation}\label{eq:1408}
\begin{split}
&\E\br[\Big]{ \abs[\big]{ v_{d,\delta}(0,x) - \mlp_{n,j,\delta}^{d,0}(0,x) }^\fq } \\
&\le 
\left( \frac{2(\sqrt{\fq-1})2^p \constantAssumpMainThm d^p(T+1)\exp(\LipConstF T)(1+2\LipConstF T)^n}{(m_j)^{n/2}\exp(-(m_j)^{\fq/2}/\fq)}
	\sup_{s\in [0,T]}
		\pr[\big]{
			\E\br[\big]{
				(1+\norm{x+W^{d,0}_s}^p)^\fq
			}
		}^{1/\fq}
\right)^{\fq} .
\end{split}
\end{equation}
It follows from \cref{eq:1408} for all $\delta\in (0,1]$, $d,j\in\N$, $n\in\N_0$ that 
\begin{equation}\label{eq:1905}
\begin{split}
&
\left( \int_{\R^d} 
\E\br[\Big]{ \abs[\big]{ v_{d,\delta}(0,x) - \mlp_{n,j,\delta}^{d,0}(0,x) }^\fq } 
 \, \nu_d(\dxx x) \right)^{\!\!\nicefrac{1}{\fq}} \\
&\le 
\frac{(\sqrt{\fq-1})2^{p+1} \constantAssumpMainThm d^p(T+1)\exp(\LipConstF T)(1+2\LipConstF T)^n}{(m_j)^{n/2}\exp(-(m_j)^{\fq/2}/\fq)}
\left(
\int_{\R^d} 
\sup_{s\in [0,T]}
		\E\br[\big]{
			(1+\norm{x+W^{d,0}_s}^p)^\fq
}
\, \nu_d(\dxx x) 
\right)^{\!\!\nicefrac{1}{\fq}}
.
\end{split}
\end{equation}
Moreover, note that Jensen's inequality and the triangle inequality ensure for all
$x \in \R$, $s \in [0,T]$, $d\in\N$ that 
$(1+\norm{x+W^{d,0}_s}^p)^{\fq} \le 2^{\fq -2 + p\fq } \left( 1 +  \norm{x}^{p\fq} + \norm{W_s^{d,0}}^{p\fq} \right)$. 
This and the assumption that for all $d \in \N$ it holds that $\nu_d$ is a probability measure imply for all $d\in\N$ that 
\begin{equation}
	\begin{split}
		\int_{\R^d} 
		\sup_{s\in [0,T]}
		\E\br[\big]{
			(1+\norm{x+W^{d,0}_s}^p)^\fq
		}
		\, \nu_d(\dxx x)  
		& \le 
		2^{(p+1)\fq} 
		\left( \!
		1\!
		+ \! \int_{\R^d} 
		\norm{x}^{p\fq} 
		\, \nu_d(\dxx x)
		+ \!
		\sup_{s\in [0,T]} 
		\E\br[\big]{\norm{W_s^{d,0}}^{p\fq} } 
		\right) 
		.
	\end{split}
\end{equation}
This, \cref{eq:1905}, and \cref{eq:deffcd} 
ensure   
for all $\delta\in (0,1]$, $d,j\in\N$, $n\in\N_0$ that 
\begin{equation}\label{eq:1553}
\begin{split}
\left( \int_{\R^d} 
\E\br[\Big]{ \abs[\big]{ v_{d,\delta}(0,x) - \mlp_{n,j,\delta}^{d,0}(0,x) }^\fq } 
\, \nu_d(\dxx x) \right)^{\!\!\nicefrac{1}{\fq}} 
&\le 
\fc_d \left(\frac{(1+2\LipConstF T)\exp\left(\frac{(m_j)^{\frac{\fq}{2}}}{n}\right)}{(m_j)^{1/2}} \right)^n
.
\end{split}
\end{equation}
Moreover, note that \cref{eq:estimateBrownianExpect} assures for all $s\in [0,T]$, $d\in\N$ that 
\begin{equation}
\E\br[\big]{
	\norm{W^{d,0}_s}^{p\fq} }
\le 1 + (2T+1)^{\frac{p(\fq+1)}{2}} \left( \frac{d}{2} + p(\fq+1) -1 \right)^{\frac{p(\fq+1)}{2}} 
.
\end{equation}
It follows for all $d\in\N$ that 
\begin{equation}\label{eq:estimateSupEW}
\left( \sup_{s\in[0,T]} \E\br[\big]{
	\norm{W^{d,0}_s}^{p\fq} } \right)^{\frac{1}{\fq}}
\le 1 + (2T+1)^p \left( \frac{d}{2} + p(\fq+1) -1 \right)^p .
\end{equation}
Furthermore, 
the assumptions that for all $d\in\N$ it holds that $(\int_{\R^{d}} \norm{y}^{p^2\fq} \,\nu_{d}(\dxx y)) \le \constantAssumpMainThm d^{rp^2\fq}$ and that $\nu_d$ is a probability measure 
imply for all $d\in\N$ that 
\begin{equation}
	\left(\int_{\R^d} 
	\norm{x}^{p \fq}
	\, \nu_d(\dxx x) \right)^\frac{1}{\fq} 
	\le \left( 1 + \int_{\R^d} 
	\norm{x}^{p^2 \fq}
	\, \nu_d(\dxx x)  \right)^\frac{1}{\fq}
	\le 1 + \constantAssumpMainThm d^{rp^2} .
\end{equation}
Combining this, \cref{eq:estimateSupEW}, and \cref{eq:deffcd} 
proves that there exists $\overline{\fC} \in [1,\infty)$ such that 
for all $d \in \N$ it holds that 
\begin{equation}\label{eq:boundforfcd}
\fc_d \le \overline{\fC} d^{rp^2 + p} .
\end{equation} 
The fact that 
$\limsup_{n\to\infty} \nicefrac{(\littleM_n)^{\fq/2}}{n} < \infty$ and the fact that $\liminf_{n\to\infty} \littleM_n = \infty$
imply that 
\begin{equation}
\limsup_{n\to\infty} 
\, \left(\frac{(1+2\LipConstF T)\exp\left(\frac{(m_n)^{\frac{\fq}{2}}}{n}\right)}{(m_n)^{1/2}} \right)^{\! n}
= 0 .
\end{equation}
This, \cref{eq:boundforfcd}, and \cref{eq:deffR} show that it holds for all $d\in\N$, $\delta \in (0,1]$ that $\fR(d,\delta)<\infty$. 
It thus follows from \cref{eq:deffR} and \cref{eq:1553} 
for all $\delta \in (0,1]$, $d\in\N$ that 
\begin{equation}\label{big_item3}
\begin{split}
\left( \int_{\R^d} 
\!
\E\br[\Big]{ \abs[\big]{ v_{d,\delta}(0,x) - \mlp_{\fR(d,\delta),\fR(d,\delta),\delta}^{d,0}(0,x) }^\fq } 
\, \nu_d(\dxx x) \right)^{\!\!\nicefrac{1}{\fq}} 
& \! \le 
\fc_d \left(\frac{(1+2\LipConstF T)\exp\left(\frac{(m_{\fR(d,\delta)})^{\frac{\fq}{2}}}{\fR(d,\delta)}\right)}{(m_{\fR(d,\delta)})^{1/2}} \right)^{\!\!\fR(d,\delta)}
\!\!\!\!\! \le \delta .
\end{split}
\end{equation}
Combining 
\Enum{
	this, 
	\cref{big_item1}, and 
	\cref{big_item2a}
	\hence\unskip
}
that for all 
$d \in \N$, 
$\delta \in (0,1]$ 
it holds that
\begin{equation}
\left( \int_{\R^d} \E\br[\Big]{ \abs[\big]{ u_d(0,x) - \mlp_{\fR(d,\delta),\fR(d,\delta),\delta}^{d,0}(0,x) }^\fq } \, \nu_d(\dxx x) \right)^{\!\!\nicefrac{1}{\fq}} 
\le 
c_d \delta + \delta 
\dpp
\end{equation}
\Enum{
	This
	;
	Fubini's theorem
}
that 
for all 
$d \in \N$, 
$\varepsilon \in (0,1]$ 
it holds that
\begin{equation}\label{error_bd_final1}
\begin{split}
& 
\E \left[ \int_{\R^d} \abs[\Big]{ \smallU_d(0,x) - \mlp_{\fR(d,\delta_{d,\varepsilon}),\fR(d,\delta_{d,\varepsilon}),\delta_{d,\varepsilon}}^{d,0}(0,x) }^\fq \, \nu_d(\dxx x) \right] 
\\
& 
= 
\int_{\R^d} \E\br[\bigg]{ \abs[\Big]{ \smallU_d(0,x) - \mlp_{\fR(d,\delta_{d,\varepsilon}),\fR(d,\delta_{d,\varepsilon}),\delta_{d,\varepsilon}}^{d,0}(0,x) }^\fq } \, \nu_d(\dxx x) 
\le 
\pr[\big]{ c_d \delta_{d,\varepsilon} + \delta_{d,\varepsilon} }^{\!\fq}
= 
\varepsilon^\fq 
\dpp
\end{split}
\end{equation}
Hence, there exists 
$(\omega_{d,\varepsilon})_{(d,\varepsilon)\in \N\times (0,1]} \subseteq \Omega$, which is assumed to be fixed for the remainder of this proof,  
such that for all $d \in \N$, $\varepsilon \in (0,1]$ it holds that 
\begin{equation}\label{int_bd_final1}
\int_{\R^d} \abs[\big]{ \smallU_d(0,x) - \mlp_{\fR(d,\delta_{d,\varepsilon}),\fR(d,\delta_{d,\varepsilon}),\delta_{d,\varepsilon}}^{d,0} (0,x,\omega_{d,\varepsilon}) }^\fq \, \nu_d(\dxx x) 
\le 
\varepsilon^\fq 
\dpp
\end{equation}
Furthermore, note that 
\Enum{
	\cref{mlp_final_th_rep}
	;
	\cref{lemma15_item2b,lemma15_item7} in \cref{lemma15} 
	(applied for every 
	$d, j \in \N$, 
	$\varepsilon\in(0,1]$, 
	$\omega \in \Omega$ 
	with 
	$\Theta \with \Theta$, 
	$d \with d$, 
	$M \with \littleMM_j$, 
	$\F \with \fnlnn_{0,\varepsilon}$, 
	$\G \with \ftermnn_{d,\varepsilon}$, 
	$(\cU^\theta)_{\theta\in\Theta} \with (\cU^\theta(\omega))_{\theta\in\Theta}$, 
	$(W^\theta)_{\theta\in\Theta} \with (W^\theta(\omega))_{\theta\in\Theta}$, 
	$(\mlp_n^\theta)_{(n,\theta) \in \Z \times \Theta} \with (\mlp_{n,j,\varepsilon}^{d,\theta}(\omega))_{(n,\theta)\in \Z \times \Theta}$, \\
	$(\U_{n,t}^\theta)_{(n,t,\theta) \in \Z \times [0,T] \times \Theta} \with (\U_{n,j,t}^{d,\theta,\varepsilon}(\omega))_{(n,t,\theta) \in \Z \times [0,T] \times \Theta}$ 
	in the notation of \cref{lemma15})
}
that for all 
$d \in \N$, 
$\varepsilon \in (0,1]$,
$x \in \R^d$ 
it holds that 
$(\Ra( \U_{\fR(d,\delta_{d,\varepsilon}),\fR(d,\delta_{d,\varepsilon}),0}^{d,0,\delta_{d,\varepsilon}} (\omega_{d,\varepsilon}) )) \in C(\R^d,\R)$ and 
\begin{equation}
\left(\Ra\left( \U_{\fR(d,\delta_{d,\varepsilon}),\fR(d,\delta_{d,\varepsilon}),0}^{d,0,\delta_{d,\varepsilon}} (\omega_{d,\varepsilon}) \right)\right)(x) = \mlp_{\fR(d,\delta_{d,\varepsilon}),\fR(d,\delta_{d,\varepsilon}),\delta_{d,\varepsilon}}^{d,0} (0,x,\omega_{d,\varepsilon}) 
.
\end{equation} 
Combining 
\Enum{
	this 
	and 
	\cref{int_bd_final1}
}[establish]
\cref{final_item2_error}. 
\Moreover
\Enum{
	\cref{lemma15_item8} in \cref{lemma15} 
}
that 
for all 
$d \in \N$, 
$\varepsilon \in (0,1]$ 
it holds that
\begin{align}\label{big_param_est}
& 
\paramANN\pr[\big]{
\U_{\fR(d,\delta_{d,\varepsilon}),\fR(d,\delta_{d,\varepsilon}),0}^{d,0,\delta_{d,\varepsilon}} (\omega_{d,\varepsilon}) 
} 
\nonumber 
\\
& 
\le
2
\pr[\big]{ \max\cu[\big]{\fd,\lengthANN(\ftermnn_{d,\delta_{d,\varepsilon}})} + \fR(d,\delta_{d,\varepsilon}) \hiddenLength(\fnlnn_{0,\delta_{d,\varepsilon}}) }
\left( \max\cu[\big]{ \fd , \normmm{\dims(\fnlnn_{0,\delta_{d,\varepsilon}}) }, \normmm{\dims(\ftermnn_{d,\delta_{d,\varepsilon}})} }  \right)^2 
\nonumber
\\
& \quad \cdot  (3\littleM_{\fR(d,\delta_{d,\varepsilon})})^{2\fR(d,\delta_{d,\varepsilon})}
\nonumber 
\\
& 
\le
2
\pr[\big]{ \max\cu[\big]{\fd,\lengthANN(\ftermnn_{d,\delta_{d,\varepsilon}})} + \hiddenLength(\fnlnn_{0,\delta_{d,\varepsilon}}) } 
\left( \max\cu[\big]{ \fd , \normmm{\dims(\fnlnn_{0,\delta_{d,\varepsilon}}) }, \normmm{\dims(\ftermnn_{d,\delta_{d,\varepsilon}})} }  \right)^2  \\
& \quad \cdot \br[\Big]{ \pr[\big]{ \fR(d,\delta_{d,\varepsilon}) }^{\nicefrac{1}{2}}  \pr[\big]{3 \littleMM_{\fR(d,\delta_{d,\varepsilon})} }^{\fR(d,\delta_{d,\varepsilon})} }^2 
\dpp
\nonumber
\end{align}
\Moreover
\Enum{
	the assumption that for all 
	$d \in \N$, 
	$\varepsilon \in (0,1]$ 
	it holds that
	$\lengthANN(\fnlnn_{0,\varepsilon}) \le \constantAssumpMainThm \varepsilon^{-\alpha_0}$,
	$\normmm{\dims(\fnlnn_{0,\varepsilon})} \le \constantAssumpMainThm \varepsilon^{-\beta_0}$,
	$\lengthANN(\ftermnn_{d,\varepsilon}) \le  \constantAssumpMainThm d^p\varepsilon^{-\alpha_1}$,
	and 
	$\normmm{\dims(\ftermnn_{d,\varepsilon})} \le \constantAssumpMainThm d^p\varepsilon^{-\beta_1}$
}
that for all 
$d \in \N$, 
$\varepsilon \in (0,1]$ 
it holds that
\begin{equation}\label{big_param_est_first_part}
\begin{split}
& 
\pr[\big]{ \max\cu[\big]{\fd,\lengthANN(\ftermnn_{d,\delta_{d,\varepsilon}})} + \hiddenLength(\fnlnn_{0,\delta_{d,\varepsilon}}) } 
\left( \max\cu[\big]{ \fd , \normmm{\dims(\fnlnn_{0,\delta_{d,\varepsilon}}) }, \normmm{\dims(\ftermnn_{d,\delta_{d,\varepsilon}})} }  \right)^2 
\\
& 
\le 
\pr[\big]{ \constantAssumpMainThm d^p \delta_{d,\varepsilon}^{-\alpha_1} + \constantAssumpMainThm \delta_{d,\varepsilon}^{-\alpha_0} } 
\pr[\Big]{ \max\cu[\big]{ \constantAssumpMainThm , \constantAssumpMainThm \delta_{d,\varepsilon}^{-\beta_0}, \constantAssumpMainThm d^p \delta_{d,\varepsilon}^{-\beta_1} } }^{\!2} 
\\
& 
\le 
2 \constantAssumpMainThm^3 d^{3p} \max\cu[\big]{ \delta_{d,\varepsilon}^{-\alpha_1} , \delta_{d,\varepsilon}^{-\alpha_0} }
\pr[\Big]{ \max\cu[\big]{  \delta_{d,\varepsilon}^{-\beta_0},  \delta_{d,\varepsilon}^{-\beta_1} } }^{\!2} 
\\
& 
\le 
2\constantAssumpMainThm^3 d^{3p} 
\delta_{d,\varepsilon}^{-(\max\{\alpha_0,\alpha_1\} + 2 \max\{\beta_0,\beta_1\})}
\dpp
\end{split}
\end{equation}
It follows from \cref{eq:deffR} that it holds for all $d\in\N$, $\varepsilon \in (0,1]$ with $\fR(d,\delta_{d,\varepsilon}) \in \N\cap[2,\infty)$ that 
\begin{equation}
\fc_d  \left(\frac{(1+2\LipConstF T)\exp\left(\frac{(m_{\fR(d,\delta_{d,\varepsilon})-1})^{\frac{\fq}{2}}}{\fR(d,\delta_{d,\varepsilon})-1}\right)}{(m_{\fR(d,\delta_{d,\varepsilon})-1})^{1/2}} \right)^{\!(\fR(d,\delta_{d,\varepsilon})-1)}
> \delta_{d,\varepsilon} .
\end{equation}
This implies for all $d\in\N$, $\varepsilon \in (0,1]$, $\delta\in (0,\infty)$ with $\fR(d,\delta_{d,\varepsilon}) \in \N\cap[2,\infty)$
that 
\begin{equation}\label{eq:1023}
\begin{split}
& \pr[\big]{ \fR(d,\delta_{d,\varepsilon}) }^{\nicefrac{1}{2}}  \pr[\big]{3 \littleMM_{\fR(d,\delta_{d,\varepsilon})} }^{\fR(d,\delta_{d,\varepsilon})} \\ 
& 
\le 
\pr[\big]{ \fR(d,\delta_{d,\varepsilon}) }^{\nicefrac{1}{2}}  \pr[\big]{3 \littleMM_{\fR(d,\delta_{d,\varepsilon})} }^{\fR(d,\delta_{d,\varepsilon})} 
\left( \fc_d \delta_{d,\varepsilon}^{-1}  \left(\frac{(1+2\LipConstF T)\exp\left(\frac{(m_{\fR(d,\delta_{d,\varepsilon})-1})^{\frac{\fq}{2}}}{\fR(d,\delta_{d,\varepsilon})-1}\right)}{(m_{\fR(d,\delta_{d,\varepsilon})-1})^{1/2}} \right)^{\!(\fR(d,\delta_{d,\varepsilon})-1)} \right)^{2+\delta} \\
& 
\le \frac{\fc_d^{2+\delta}}{\delta_{d,\varepsilon}^{2+\delta}}
\, \sup_{n\in\N} 
\left( \frac{\littleM_{n+1}^{n+1}}{\littleM_n^{n(1+\frac{\delta}{2})}}
(n+1)^{\frac{1}{2}} 
\left( 3 (1+2\LipConstF T) \exp\left( \frac{\littleM_n^{\frac{\fq}{2}}}{n} \right) \right)^{n(2+\delta)}
\right)
.
\end{split}
\end{equation}
Observe that the fact that it holds for all $n\in\N$ that $\littleM_{n+1}\le \fB \littleM_n$ ensures for all $n\in\N$ that $(\nicefrac{\littleM_{n+1}}{\littleM_n})^n \leq \fB^n$ and $\littleM_{n+1} \le \fB^n \littleM_1$. Therefore, it holds for all $n\in\N$ that $\nicefrac{\littleM_{n+1}^{n+1}}{\littleM_n^{n(1+\frac{\delta}{2})}} \le \nicefrac{\fB^{2n} \littleM_1}{\littleM_n^{\nicefrac{n\delta}{2}}}$. 
This and \cref{eq:1023} yield for all $d\in\N$, $\varepsilon \in (0,1]$, $\delta\in (0,\infty)$ with $\fR(d,\delta_{d,\varepsilon}) \in \N\cap[2,\infty)$
that  
\begin{equation}\label{eq:1101}
\begin{split}
\pr[\big]{ \fR(d,\delta_{d,\varepsilon}) }^{\nicefrac{1}{2}}  \pr[\big]{3 \littleMM_{\fR(d,\delta_{d,\varepsilon})} }^{\fR(d,\delta_{d,\varepsilon})}
& \le 
\frac{\fc_d^{2+\delta} \littleM_1}{\delta_{d,\varepsilon}^{2+\delta}}
\, \sup_{n\in\N} 
\left( 
\frac{(n+1)^{\frac{1}{2n}} \left( 3 (1+2\LipConstF T)\fB \exp\left( \frac{\littleM_n^{\frac{\fq}{2}}}{n} \right) \right)^{2+\delta} }{\littleM_n^{\frac{\delta}{2}}}
\right)^{\! n}
.
\end{split}
\end{equation}
Observe that the fact that $\littleM_1\in \N$, the fact that $\fB\in [1,\infty)$, the fact that $\fq \in [2,\infty)$, the fact that for all $d\in\N$ it holds that $\fc_d\in [1,\infty)$, and the fact that for all $d\in\N$, $\varepsilon\in (0,1]$ it holds that $\delta_{d,\varepsilon}\in(0,1]$ 
show for all $d\in\N$, $\varepsilon \in (0,1]$, $\delta\in (0,\infty)$ 
that 
\begin{equation}
\begin{split}
& \frac{\fc_d^{2+\delta} \littleM_1}{\delta_{d,\varepsilon}^{2+\delta}}
\, \sup_{n\in\N} 
\left( 
\frac{(n+1)^{\frac{1}{2n}} \left( 3 (1+2\LipConstF T)\fB \exp\left( \frac{\littleM_n^{\frac{\fq}{2}}}{n} \right) \right)^{2+\delta} }{\littleM_n^{\frac{\delta}{2}}}
\right)^{\! n} \\
& \ge 
\littleM_1 \cdot 
\frac{2^{\frac{1}{2}} \left( 3 (1+2\LipConstF T)\fB \exp\left( \littleM_1^{\frac{\fq}{2}} \right) \right)^{2+\delta} }{\littleM_1^{\frac{\delta}{2}}}
\ge 3 \littleM_1 \littleM_1^{\fq + \frac{\delta}{2} \fq - \frac{\delta}{2}} 
\ge 3 \littleM_1 .
\end{split}
\end{equation}
Combining this and the fact that it holds for all $d\in\N$, $\varepsilon \in (0,1]$ that $\fR(d,\delta_{d,\varepsilon})<\infty$
assures 
that \cref{eq:1101} holds  
for all $d\in\N$, $\varepsilon \in (0,1]$, $\delta\in (0,\infty)$.
This, \cref{big_param_est}, and \cref{big_param_est_first_part} 
ensure for all $d\in\N$, $\varepsilon \in (0,1]$, $\delta\in (0,\infty)$ 
that
\begin{equation}\label{eq:1541a}
\begin{split}
\paramANN\pr[\big]{
	\U_{\fR(d,\delta_{d,\varepsilon}),\fR(d,\delta_{d,\varepsilon}),0}^{d,0,\delta_{d,\varepsilon}} (\omega_{d,\varepsilon}) }
& \le 
4 \constantAssumpMainThm^3 d^{3p} 
\delta_{d,\varepsilon}^{-(2(2+\delta)+\max\{\alpha_0,\alpha_1\} + 2 \max\{\beta_0,\beta_1\})} 
\fc_d^{2(2+\delta)} \littleM_1^2 \\
& \quad \cdot \sup_{n\in\N} 
\left( 
\frac{(n+1)^{\frac{1}{2n}} \left( 3 (1+2\LipConstF T)\fB \exp\left( \frac{\littleM_n^{\frac{\fq}{2}}}{n} \right) \right)^{2+\delta} }{\littleM_n^{\frac{\delta}{2}}}
\right)^{\!2 n}
.
\end{split}
\end{equation}
Note that \cref{eq:boundforcd} establishes for all $d\in\N$, $\varepsilon\in(0,1]$, $\delta \in (0,\infty)$ that 
\begin{equation}\label{eq:1541b}
\begin{split}
& \delta_{d,\varepsilon}^{-(2(2+\delta)+\max\{\alpha_0,\alpha_1\} + 2 \max\{\beta_0,\beta_1\})} 
= \left(\frac{c_d+1}{\varepsilon}\right)^{\kappa_{\delta}} 
\le (2 \overline{C})^{\kappa_{\delta}} 
\varepsilon^{-\kappa_{\delta}}
d^{(p+(r+1)p^2)\kappa_{\delta}} .
\end{split}
\end{equation}
Moreover, \cref{eq:boundforfcd} proves for all $d\in\N$, $\delta \in (0,\infty)$ that 
\begin{equation}\label{eq:1541c}
\begin{split}
\fc_d^{2(2+\delta)} 
& \le \overline{\fC}^{2(2+\delta)} d^{2(2+\delta)(rp^2 + p)} .
\end{split}
\end{equation} 
Furthermore, 
the fact that 
$\limsup_{n\to\infty} \nicefrac{(\littleM_n)^{\fq/2}}{n} < \infty$, the fact that $\liminf_{n\to\infty} \littleM_n = \infty$, and the fact that  $\lim_{n\to\infty} (n+1)^{\frac{1}{2n}}=1$ 
yield for all $\delta \in (0,\infty)$ that 
\begin{equation}\label{eq:1541d}
\sup_{n\in\N} 
\left( 
\frac{(n+1)^{\frac{1}{2n}} \left( 3 (1+2\LipConstF T)\fB \exp\left( \frac{\littleM_n^{\frac{\fq}{2}}}{n} \right) \right)^{2+\delta} }{\littleM_n^{\frac{\delta}{2}}}
\right)^{\! 2 n} 
< \infty .
\end{equation} 
Combining \cref{eq:1541a}, \cref{eq:1541b}, \cref{eq:1541c}, and \cref{eq:1541d} shows that there exists $\eta\colon (0,\infty) \to \R$ such that 
for all $d\in\N$, $\varepsilon \in (0,1]$, $\delta\in (0,\infty)$ 
it holds that 
\begin{equation}
\paramANN\pr[\big]{
	\U_{\fR(d,\delta_{d,\varepsilon}),\fR(d,\delta_{d,\varepsilon}),0}^{d,0,\delta_{d,\varepsilon}} (\omega_{d,\varepsilon}) } 
\le 
\eta(\delta)  
d^{3p + 2(2+\delta)(rp^2 + p) + (p+(r+1)p^2)\kappa_{\delta}}
\varepsilon^{-\kappa_{\delta}}
.
\end{equation}
This 
establishes
\cref{final_item2}.
\end{aproof}

%
%
%


\subsection{One-dimensional ANN approximation results}\label{sec:3.2}

\subsubsection{The modulus of continuity}

\begin{definition}[Modulus of continuity]\label{mod_cont_def}
Let $A \subseteq \R$ be a set and let $f \colon A \to \R$ be a function.
Then we denote by 
$
\modCont_f \colon [0,\infty] \to [0,\infty]
$ 
the function which satisfies for all $h\in [0,\infty]$ that
\begin{equation}\label{eq:mod_cont}
\modCont_f(h) 
= 
\sup\pr[\Big]{ \cu[\big]{ \abs{ f(x)-f(y) } \in [0,\infty) \colon \pr[\big]{ x,y\in A \text{ with } \abs{x-y} \le h } } \cup \{0\} }
\end{equation}
and we call $\modCont_f$ the modulus of continuity of $f$.
\end{definition}

\cfclear
\begin{athm}{lemma}{lem:mod_continuity}
Let 
$b_1 \in [-\infty,\infty]$, 
$b_2 \in[b_1,\infty]$ 
and let 
$ f \colon ( [b_1,b_2] \cap \R ) \to \R $ 
be a function. 
Then
\begin{enumerate}[label=(\roman *)]
\item
\label{mod_item1}
it holds that $ \modCont_f $ is non-decreasing,
\item
\label{mod_item2}
it holds that
$ f $ is uniformly continuous if and only if 
$\lim_{ h \searrow 0 } \modCont_f( h ) = 0$,
\item
\label{mod_item3}
it holds that
$ f $ is globally bounded if and only if
$ \modCont_f( \infty ) < \infty $,
\item
\label{mod_item4}
it holds
for all 
$ x, y \in [b_1,b_2] \cap \R $ 
that
$
\abs{f(x) - f(y) } 
\leq
\modCont_f (\abs{x- y})
$,
and
\item
\label{mod_item5}
it holds for all 
$h, \fh \in[0,\infty]$ 
that 
$\modCont_f(h + \fh) \le \modCont_f(h) + \modCont_f(\fh)$
\end{enumerate}
\cfout.
\end{athm}

\begin{aproof}
First, \nobs that 
\Enum{
\cref{eq:mod_cont}
}[imply]
\cref{mod_item1,mod_item2,mod_item3,mod_item4}.
\Moreover
\Enum{
\cref{eq:mod_cont}
}
that for all $h, \fh \in [0,\infty]$ and for all $x,y \in [b_1,b_2] \cap \R$ that satisfy $\abs{ x-y } \le \max\{h,\fh\}$ 
it holds that 
$\abs[]{ f(x)-f(y) } \le \modCont_f(h) + \modCont_f(\fh)$. 
\Moreover
\Enum{
\cref{eq:mod_cont}
}
that for all $h, \fh \in [0,\infty]$ and for all $x,y \in [b_1,b_2] \cap \R$ that satisfy $\max\{h,\fh\} < \abs{ x-y }  \le (h+\fh)$ 
it holds that 
$x-h\tfrac{x-y}{\abs{x-y}} \in [b_1,b_2] \cap \R$, 
that 
$\abs{ x -  \pr[\big]{ x-h\tfrac{x-y}{\abs{x-y}} }  } =h$, 
that 
$\abs{ \pr[\big]{ x-h\tfrac{x-y}{\abs{x-y}} }  - y } =\abs{x-y}-h\le \fh$, 
and thus 
\begin{align}
\abs[]{ f(x)-f(y) } 
& \le \abs[\big]{ f(x) - f\pr[\big]{ x-h\tfrac{x-y}{\abs{x-y}} } } 
+ \abs[\big]{ f\pr[\big]{ x-h\tfrac{x-y}{\abs{x-y}} } - f(y) }
\nonumber 
\\
& \le \modCont_f(h) + \modCont_f(\fh) .
\end{align}
Combining both cases, 
\Enum{
	\cref{eq:mod_cont}
}
that for all $h, \fh \in [0,\infty]$ it holds that
$\modCont_f(h+\fh) \le \modCont_f(h) + \modCont_f(\fh)$ 
\cfload.
This establishes \cref{mod_item5}.
\end{aproof}

\cfclear
\begin{athm}{lemma}{lem:mod_lip}
Let 
$A \subseteq \R$, 
$L \in [0,\infty)$, 
and let 
$ f \colon A \to \R $ 
satisfy for all 
$x,y\in A$ 
that 
$\abs{f(x)-f(y)} \le L \abs{x-y}$.
Then it holds for all 
$h \in [0,\infty)$ 
that 
$\modCont_f(h) \le Lh$
\cfout.
\end{athm}

\begin{aproof}
\Nobs that
\Enum{
the 
assumption that for all 
$x,y\in A$ 
it holds that 
$\abs{f(x)-f(y)} \le L\abs{x-y}$
;
\cref{eq:mod_cont}
}
that for all 
$h \in [0,\infty)$ 
it holds that 
\begin{align}
\modCont_f(h) 
& 
=  
\sup\pr[\Big]{ \cu[\big]{ \abs{ f(x)-f(y) } \in [0,\infty) \colon \pr[\big]{ x,y\in A \text{ with } \abs{x-y} \le h } } \cup \{0\} } 
\\
& 
\le  
\sup\pr[\Big]{ \cu[\big]{ L \abs{ x-y } \in [0,\infty) \colon \pr[\big]{ x,y\in A \text{ with } \abs{x-y} \le h } } \cup \{0\} } 
\le 
\sup\pr[\big]{ \cu[]{ Lh,0 } } 
= 
Lh 
\nonumber
\end{align}
\cfout.
\end{aproof}

\subsubsection{Linear interpolation of one-dimensional functions}

\begin{definition}[Linear interpolation function]
\label{def:lin_interp}
Let 
$K\in\N$, 
$\fx_0,\fx_1,\dots,\fx_K, f_0, f_1, \dots, f_K \in\R$ 
satisfy 
$\fx_0 < \fx_1 < \dots < \fx_K$. 
Then we denote by 
$\scrL_{\fx_0,\fx_1,\dots,\fx_K}^{f_0,f_1,\dots,f_K} \colon \R \to \R$ 
the function which satisfies for all 
$k\in\{1,2,\dots,K\}$, 
$x\in(-\infty,\fx_0)$, 
$y\in[\fx_{k-1},\fx_k)$, 
$z\in[\fx_K,\infty)$ 
that 
$(\scrL_{\fx_0,\fx_1,\dots,\fx_K}^{f_0,f_1,\dots,f_K})(x) = f_0$, 
$(\scrL_{\fx_0,\fx_1,\dots,\fx_K}^{f_0,f_1,\dots,f_K})(z) = f_K$, 
and
\begin{equation}\label{eq:lin_interp}
(\scrL_{\fx_0,\fx_1,\dots,\fx_K}^{f_0,f_1,\dots,f_K})(y) 
= 
f_{k-1} + \pr[\big]{ \tfrac{y - \fx_{k-1}}{\fx_k - \fx_{k-1}} } (f_k - f_{k-1})
\dpp
\end{equation}
\end{definition}

\cfclear
\begin{athm}{lemma}{lin_interp1_new}
Let 
$K\in\N$, 
$\fx_0,\fx_1,\dots,\fx_K, f_0, f_1, \dots, f_K \in\R$ 
satisfy that 
$\fx_0 < \fx_1 < \dots < \fx_K$. 
Then
\begin{enumerate}[label=(\roman *)]
\item
\label{lin_interp1_item1_new}
it holds for all 
$k\in\{0,1,\dots,K\}$ 
that 
$(\scrL_{\fx_0,\fx_1,\dots,\fx_K}^{f_0,f_1,\dots,f_K})(\fx_k) = f_k$,
\item
\label{lin_interp1_item2_new}
it holds for all 
$k\in\{1,2,\dots,K\}$, 
$x\in[\fx_{k-1},\fx_k]$ 
that 
\begin{equation}
(\scrL_{\fx_0,\fx_1,\dots,\fx_K}^{f_0,f_1,\dots,f_K})(x) 
= 
f_{k-1} + \pr[\big]{ \tfrac{ x-\fx_{k-1} }{ \fx_k - \fx_{k-1} } } ( f_k - f_{k-1} )
\dc
\end{equation}
and
\item
\label{lin_interp1_item3_new}
it holds for all 
$k \in \{1,2,\dots,K\}$, 
$x \in [\fx_{k-1},\fx_k]$ 
that
\begin{equation}
(\scrL_{\fx_0,\fx_1,\dots,\fx_K}^{f_0,f_1,\dots,f_K})(x) 
= 
\pr[\big]{ \tfrac{ \fx_k - x }{ \fx_k - \fx_{k-1} } } f_{k-1} + \pr[\big]{ \tfrac{ x-\fx_{k-1} }{ \fx_k - \fx_{k-1} } } f_k
\end{equation}
\end{enumerate}
\cfout.
\end{athm}

\begin{aproof}
\Nobs that 
\Enum{
\cref{eq:lin_interp}
}
\cref{lin_interp1_item1_new,lin_interp1_item2_new}.
\Moreover
\Enum{
\cref{lin_interp1_item2_new}
}
that for all 
$k \in \{1,2,\dots,K\}$, 
$x\in[\fx_{k-1},\fx_k]$ 
it holds that
\begin{equation}
(\scrL_{\fx_0,\fx_1,\dots,\fx_K}^{f_0,f_1,\dots,f_K})(x) 
= 
\br[\Big]{ 1 - \pr[\big]{ \tfrac{ x-\fx_{k-1} }{ \fx_k - \fx_{k-1} } } } f_{k-1}
+ \pr[\big]{ \tfrac{ x-\fx_{k-1} }{ \fx_k - \fx_{k-1} } } f_k
= 
\pr[\big]{ \tfrac{ \fx_k - x }{ \fx_k - \fx_{k-1} } } f_{k-1} + \pr[\big]{ \tfrac{ x-\fx_{k-1} }{ \fx_k - \fx_{k-1} } } f_k 
\end{equation}
\cfload.
This proves \cref{lin_interp1_item3_new}.
\end{aproof}

\cfclear
\begin{athm}{lemma}{lin_interp1_new1}
Let 
$K\in\N$, 
$\fx_0,\fx_1,\dots,\fx_K\in\R$ 
satisfy 
$\fx_0 < \fx_1 < \dots < \fx_K$ 
and let 
$ f \colon [\fx_0,\fx_K] \to \R$ 
be a function. 
Then
\begin{enumerate}[label=(\roman *)]
\item
\label{lin_interp1_item1_new1}
it holds for all 
$x,y\in\R$ 
that 
\begin{equation}
\begin{split}
&
\abs[\big]{ (\scrL_{\fx_0,\fx_1,\dots,\fx_K}^{f(\fx_0),f(\fx_1),\dots,f(\fx_K)})(x) - (\scrL_{\fx_0,\fx_1,\dots,\fx_K}^{f(\fx_0),f(\fx_1),\dots,f(\fx_K)})(y) } 
\\
&
\le 
\br[\bigg]{ \max_{k\in\{1,2,\dots,K\}} \pr[\Big]{ \tfrac{\modCont_f(\abs{\fx_k-\fx_{k-1}})}{\abs{\fx_k - \fx_{k-1}}} } } \abs{ x-y }
\end{split}
\end{equation}
and
\item
\label{lin_interp1_item2_new1}
it holds that 
$
\sup_{x\in[\fx_0,\fx_K]} \abs[\big]{ (\scrL_{\fx_0,\fx_1,\dots,\fx_K}^{f(\fx_0),f(\fx_1),\dots,f(\fx_K)})(x) - f(x) } 
\le 
\modCont_f( \max_{k\in\{1,2,\dots,K\}} \abs{ \fx_k - \fx_{k-1} } ) $
\end{enumerate}
\cfout.
\end{athm}

\begin{aproof}
Throughout this proof let 
$l \colon \R \to \R$ 
satisfy for all 
$x\in\R$ 
that 
$l(x) = (\scrL_{\fx_0,\fx_1,\dots,\fx_K}^{f(\fx_0),f(\fx_1),\dots,f(\fx_K)})(x)$ 
and let
$L \in [0,\infty]$ 
satisfy
\begin{equation}\label{eq:3.12}
L 
= 
\max_{k\in\{1,2,\dots,K\}} \pr[\Big]{ \tfrac{\modCont_f(\abs{\fx_k-\fx_{k-1}})}{\abs{\fx_k - \fx_{k-1}}} }
\end{equation}
\cfload.
\Nobs that 
\Enum{
\cref{mod_item4} in \cref{lem:mod_continuity}
;
\cref{eq:3.12}
;
\cref{lin_interp1_item2_new} in \cref{lin_interp1_new} 
}
that for all 
$k\in\{1,2,\dots,K\}$, 
$x,y\in[\fx_{k-1},\fx_k]$
it holds that
\begin{equation}\label{eq:3.13}
\begin{split}
\abs{ l(x) - l(y) } 
& 
= 
\abs[\Big]{ \pr[\big]{ \tfrac{ x-\fx_{k-1} }{ \fx_k - \fx_{k-1} } } \pr[\big]{ f(\fx_k) - f(\fx_{k-1}) } - \pr[\big]{ \tfrac{ y-\fx_{k-1} }{ \fx_k - \fx_{k-1} } } \pr[\big]{ f(\fx_k) - f(\fx_{k-1}) } } 
\\
& 
= 
\abs[\bigg]{ \pr[\bigg]{ \frac{ f(\fx_k) - f(\fx_{k-1})}{\fx_k - \fx_{k-1}} } ( x-y ) } 
\le \pr[\bigg]{ \frac{ \modCont_f( \abs{ \fx_k - \fx_{k-1} })}{ \abs{ \fx_k - \fx_{k-1} }} } \abs{ x-y }
\le 
L \abs{x-y} 
\dpp
\end{split}
\end{equation}
\Enum{
This
;
\cref{mod_item4} in \cref{lem:mod_continuity}
;
\cref{lin_interp1_new}
;
\cref{eq:3.12}
;
the triangle inequality
}
that for all 
$j,k \in \{1,2,\dots,K\}$, 
$x \in [\fx_{j-1},\fx_j]$, 
$y \in [\fx_{k-1},\fx_k]$ 
with $j < k$ 
it holds that
\begin{align}
\abs{ l(x) - l(y) } 
& 
\le 
\abs{ l(x) - l(\fx_j) } + \abs{ l(\fx_j) - l(\fx_{k-1}) } + \abs{ l(\fx_{k-1}) - l(y) } \nonumber 
\\
& 
= 
\abs{ l(x) - l(\fx_j) } 
+ 
\abs{ f(\fx_j) - f(\fx_{k-1}) } 
+ 
\abs{ l(\fx_{k-1}) - l(y) } 
\nonumber 
\\
& 
\le 
\abs{ l(x) - l(\fx_j) } 
+ 
\left[ \SmallSum{i=j+1}{k-1} \abs{f(\fx_i) - f(\fx_{i-1}) } \right] 
+ 
\abs{ l(\fx_{k-1}) - l(y) } 
\\
& 
\le 
\abs{ l(x) - l(\fx_j) } 
+ 
\left[ \SmallSum{i=j+1}{k-1} \modCont_f\pr[\big]{ \abs{\fx_i - \fx_{i-1} } } \right] 
+ 
\abs{ l(\fx_{k-1}) - l(y) } 
\nonumber 
\\
& 
\le 
L \left( ( \fx_j - x ) + \left[ \SmallSum{i=j+1}{k-1} ( \fx_i - \fx_{i-1} ) \right] + (y - \fx_{k-1}) \right) 
= 
L \abs{x-y}
\dpp
\nonumber
\end{align}
Combining 
\Enum{
this and \cref{eq:3.13}
}
that for all 
$x,y\in [\fx_0,\fx_K]$ 
it holds that 
$\abs{l(x)-l(y)} \le L\abs{x-y}$. 
\Enum{
This
;
the fact that for all 
$x,y\in(-\infty,\fx_0]$ 
it holds that 
$\abs{l(x)-l(y)} = 0 \le L\abs{x-y}$
;
the fact that for all 
$x,y\in[\fx_K,\infty)$ 
it holds that 
$\abs{l(x)-l(y)} = 0 \le L\abs{x-y}$
;
the triangle inequality
}
that for all 
$x,y\in\R$ 
it holds that 
$\abs{l(x)-l(y)} \le L\abs{x-y}$. 
This proves \cref{lin_interp1_item1_new1}.
\Moreover
\Enum{
\cref{eq:mod_cont}
;
\cref{lem:mod_continuity}
;
\cref{lin_interp1_item3_new} in \cref{lin_interp1_new}
;
the triangle inequality
}
that for all 
$k\in\{1,2,\dots,K\}$, 
$x\in[\fx_{k-1},\fx_k]$ 
it holds that
\begin{align}
\abs{l(x) - f(x) }
& 
= 
\abs[\bigg]{ \pr[\bigg]{ \frac{ \fx_k - x }{ \fx_k - \fx_{k-1} } } f(\fx_{k-1}) + \pr[\bigg]{ \frac{ x - \fx_{k-1} }{ \fx_k - \fx_{k-1} } } f(\fx_{k}) - f(x) } 
\nonumber 
\\
& 
= 
\abs[\bigg]{ \pr[\bigg]{ \frac{ \fx_k - x }{ \fx_k - \fx_{k-1} } } ( f(\fx_{k-1}) - f(x) ) + \pr[\bigg]{ \frac{ x - \fx_{k-1} }{ \fx_k - \fx_{k-1} } } ( f(\fx_{k}) - f(x) ) } 
\nonumber 
\\
& 
\le 
\pr[\bigg]{ \frac{ \fx_k - x }{ \fx_k - \fx_{k-1} } } \abs{ f(\fx_{k-1}) - f(x) } 
+ 
\pr[\bigg]{ \frac{ x - \fx_{k-1} }{ \fx_k - \fx_{k-1} } } \abs{ f(\fx_{k}) - f(x) } 
\\
& 
\le 
\modCont_f\pr[\big]{ \abs{\fx_k - \fx_{k-1}} } \pr[\bigg]{ \frac{ \fx_k - x }{ \fx_k - \fx_{k-1} } + \frac{ x - \fx_{k-1} }{ \fx_k - \fx_{k-1} } } 
\nonumber 
\\
& 
= 
\modCont_f\pr[\big]{ \abs{\fx_k - \fx_{k-1} } } 
\le 
\modCont_f\pr[\big]{ \textstyle\max_{j \in \{1,2,\ldots,K\}} \abs{\fx_j - \fx_{j-1}} }
\dpp 
\nonumber
\end{align}
This establishes \cref{lin_interp1_item2_new1}.
\end{aproof}

\cfclear
\begin{athm}{lemma}{lin_interp1_new2}
Let 
$K\in\N$, 
$L, \fx_0,\fx_1,\dots,\fx_K\in\R$ 
satisfy 
$\fx_0 < \fx_1 < \dots < \fx_K$ 
and let 
$ f \colon [\fx_0,\fx_K] \to \R$ 
satisfy for all 
$x,y\in[\fx_0,\fx_K]$ 
that 
$\abs{f(x)-f(y)} \le L\abs{x-y}$. 
Then
\begin{enumerate}[label=(\roman *)]
\item
\label{lin_interp1_item1_new2}
it holds for all $x,y\in\R$ that 
\begin{equation}
\abs[\big]{ (\scrL_{\fx_0,\fx_1,\dots,\fx_K}^{f(\fx_0),f(\fx_1),\dots,f(\fx_K)})(x) - (\scrL_{\fx_0,\fx_1,\dots,\fx_K}^{f(\fx_0),f(\fx_1),\dots,f(\fx_K)})(y) } 
\le 
L \abs{ x-y }
\end{equation}
and
\item
\label{lin_interp1_item2_new2}
it holds that 
$
\sup_{x\in[\fx_0,\fx_K]} \abs[\big]{ (\scrL_{\fx_0,\fx_1,\dots,\fx_K}^{f(\fx_0),f(\fx_1),\dots,f(\fx_K)})(x) - f(x) } 
\le 
L ( \max_{k\in\{1,2,\dots,K\}} \abs{ \fx_k - \fx_{k-1} } ) 
$
\end{enumerate}
\cfout.
\end{athm}

\begin{aproof}
First, \nobs that 
\Enum{
the assumption that for all 
$x,y\in[\fx_0,\fx_K]$ 
it holds that 
$\abs{f(x)-f(y)} \le L \abs{x-y}$
;
\cref{lem:mod_lip}
;
\cref{lin_interp1_item1_new1} in \cref{lin_interp1_new1} 
}
that for all 
$x,y\in\R$ 
it holds that 
\begin{equation}
\begin{split}
& 
\abs[\big]{ (\scrL_{\fx_0,\fx_1,\dots,\fx_K}^{f(\fx_0),f(\fx_1),\dots,f(\fx_K)})(x) - (\scrL_{\fx_0,\fx_1,\dots,\fx_K}^{f(\fx_0),f(\fx_1),\dots,f(\fx_K)})(y) } 
\\
& 
\le 
\br[\bigg]{ \max_{k\in\{1,2,\dots,K\}} \pr[\Big]{ \tfrac{ L \abs{\fx_k-\fx_{k-1}} }{\abs{\fx_k - \fx_{k-1}}} } } \abs{ x-y } 
= 
L \abs{x-y} 
\end{split}
\end{equation}
\cfload.
This proves \cref{lin_interp1_item1_new2}.
\Moreover
\Enum{
the assumption that for all 
$x,y\in[\fx_0,\fx_K]$ 
it holds that 
$\abs{f(x)-f(y)} \le L\abs{x-y}$
; 
\cref{lem:mod_lip}
;
\cref{lin_interp1_item2_new1} in \cref{lin_interp1_new1}
}
that
\begin{equation}
\sup_{x\in[\fx_0,\fx_K]} \abs[\big]{ (\scrL_{\fx_0,\fx_1,\dots,\fx_K}^{f(\fx_0),f(\fx_1),\dots,f(\fx_K)})(x) - f(x) }
\le 
L \left( \max_{k\in\{1,2,\dots,K\}} \abs{ \fx_k - \fx_{k-1} } \right)
\dpp
\end{equation}
This establishes \cref{lin_interp1_item2_new2}.
\end{aproof}

\subsubsection{Linear interpolation with ANNs}

\cfclear
\begin{athm}{lemma}{non_rep2}
	Let $K\in\N_0$, 
	$f_0, c_0,c_1,\ldots, c_K, \beta_0, \beta_1,\ldots, \beta_K, \alpha_0,\alpha_1,\ldots,\alpha_K\in \R$, $a\in C(\R,\R)$, 
	$\F \in \ANNs$ 
	satisfy
	\begin{equation}\label{non_rep2_ANN}
	\F 
	= 
	\compANN{ \A_{1,f_0} }{ \left( \OSum{k=0}{K} \left( \scalar{ c_k }{ ( \compANN{ \ii_{1} }{ \A_{\alpha_k,\beta_k} } )\right) } \right) }
	\end{equation}
	\cfload.
	Then 
	\begin{enumerate}[label=(\roman *)]
		\item 
		\label{non_rep2_item1}
		it holds that 
		$\dims(\F) = (1,K+1,1) \in \N^3$,
		\item 
		\label{non_rep2_item2}
		it holds that 
		$\functionANN{a}(\F) \in C(\R,\R)$,
		\item 
		\label{non_rep2_item3}
		it holds for all $x \in \R$ that 
		$(\functionANN{a}(\F))(x) = f_0 + \sum_{k=0}^K c_k (a(\alpha_k x + \beta_k))$, 
		and
		\item 
		\label{non_rep2_item4}
		it holds that 
		$\paramANN(\F) = 3K+4 = 3(\dimANNlevel_1(\F)) + 1$
	\end{enumerate}
	\cfout.
\end{athm}

\begin{aproof} 
Note that, e.g., Grohs et al.\ \cite[item (i) in Proposition~2.6]{GrohsHornungJentzen2019},  
\cref{item:sum:ANN:1} in \cref{lem:sum:ANN}, and 
\cref{padding_item2} in \cref{padding_lemma} 
prove \cref{non_rep2_item1}.
In addition, 
e.g., Grohs et al.\ \cite[item (v) in Proposition~2.6]{GrohsHornungJentzen2019},  
\cref{item:sum:ANN:3} in \cref{lem:sum:ANN}, and 
\cref{padding_item3} in \cref{padding_lemma} 
establish \cref{non_rep2_item2}.
It follows from, 
e.g., Grohs et al.\ \cite[item (v) in Proposition~2.6]{GrohsHornungJentzen2019},  
\cref{item:sum:ANN:4} in \cref{lem:sum:ANN}, and 
\cref{padding_item4} in \cref{padding_lemma} 
that for all $x\in\R$ it holds that 
\begin{equation}
\begin{split}
(\functionANN{a}(\F))(x) 
& = 
f_0 + \left(\functionANN{a}\left( \OSum{k=0}{K} \left( \scalar{ c_k }{ ( \compANN{ \ii_{1} }{ \A_{\alpha_k,\beta_k} } ) } \right) \right) \right)(x) \\
& = 
f_0 + \sum_{k=0}^K c_k ((\functionANN{a}(\ii_1))(\alpha_k x + \beta_k)) 
= f_0 + \sum_{k=0}^K c_k (a(\alpha_k x + \beta_k)) .
\end{split}
\end{equation}
This shows \cref{non_rep2_item3}. 
Moreover, note that \cref{non_rep2_item1} 
implies that 
\begin{equation}
\label{interpol_ANN_points:eq0}
\paramANN(\F)
=
(K+1) (1+1) + (K+1+1)
=
3K+4
= 3 (K+1) + 1
= 3(\dimANNlevel_1(\F)) + 1 .
\end{equation}
This proves \cref{non_rep2_item4}.
\end{aproof}

\cfclear
\begin{athm}{lemma}{non_rep3}
	Let $K\in\N$, 
	$f_0, f_1, \dots , f_K, \fx_0, \fx_1, \dots, \fx_K\in \R$
	satisfy 
	$\fx_0 < \fx_1 < \dots < \fx_K$, 
	let
	$\F \in \ANNs$ 
	satisfy
	\begin{equation}\label{interpol_ANN_points:ass1}
	\F 
	= 
	\compANN{ \A_{1,f_0} }{ \left( \OSum{k=0}{K} \left( \scalar{ \left(\tfrac{(f_{\min\{k+1,K\}}-f_k)}{(\fx_{\min\{k+1,K\}}-\fx_{\min\{k,K-1\}})} - \tfrac{(f_k - f_{\max\{k-1,0\}})}{(\fx_{\max\{k,1\}} - \fx_{\max\{k-1,0\}})}\right) }{ ( \compANN{ \ii_{1} }{ \A_{1,-\fx_k} } )\right) } \right) },
	\end{equation}
	and let $\relu\in C(\R,\R)$ satisfy for all $x\in\R$ that 
	$\relu(x)=\max\{x,0\}$ 
	\cfload.
	Then it holds that 
	$\functionANN{\relu}(\F) = \scrL_{\fx_0,\fx_1,\dots,\fx_K}^{f_0,f_1,\dots,f_K}$ 
	\cfout. 
\end{athm}

\begin{aproof}
Throughout this proof
let $ c_0, c_1, \dots, c_K  \in \R$ satisfy for all $k\in\{0,1,\allowbreak \dots, \allowbreak K\}$ that
\begin{equation}\label{interpol_ANN_points:setting1}
c_k = \frac{(f_{\min\{k+1,K\}}-f_k)}{(\fx_{\min\{k+1,K\}}-\fx_{\min\{k,K-1\}})} - \frac{(f_k - f_{\max\{k-1,0\}})}{(\fx_{\max\{k,1\}} - \fx_{\max\{k-1,0\}})} .
\end{equation}
Observe that 
\cref{non_rep2_item3} 
in \cref{non_rep2} 
ensures that for all $x\in\R$ it holds that
\begin{equation}
\label{interpol_ANN_points:eq1}
(\functionANN{\relu}(\F))(x)  =  f_0 + \SmallSum{k=0}{K} c_k \max\{x - \fx_k,0\}.
\end{equation}
This and the fact that $\forall\,k\in\{0,1,\dots,K\} \colon \fx_0 \le \fx_k$ assure that for all $x\in(-\infty,\fx_0]$ it holds that 
\begin{equation}\label{interpol_ANN_points:eq2}
(\functionANN{\relu}(\F))(x) = f_0 + 0 = f_0.
\end{equation}
Next we claim that for all
$k\in\{1,2,\dots,K\}$
it holds that
\begin{equation}
\label{interpol_ANN_points:eq1.1}
\SmallSum{n=0}{k-1} c_n = \tfrac{f_{k}-f_{k-1}}{\fx_{k}-\fx_{k-1}}.
\end{equation}
We now prove \eqref{interpol_ANN_points:eq1.1} by induction on $k\in\{1,2,\dots,K\}$.
For the base case $k=1$ observe that \cref{interpol_ANN_points:setting1} assures that 
$\sum_{n=0}^{0} c_n = c_0 = \tfrac{f_{1}-f_{0}}{\fx_{1}-\fx_{0}}$.
This proves \eqref{interpol_ANN_points:eq1.1} in the base case $k=1$.
For the induction step from $\{1,2,\ldots,K-1\} \ni (k-1) \induct k \in \{2,3,\ldots,K\}$ note that \cref{interpol_ANN_points:setting1} ensures that for all 
$k\in\{2,3,\dots,K\}$ 
with $\sum_{n=0}^{k-2} c_n = \tfrac{f_{k-1}-f_{k-2}}{\fx_{k-1}-\fx_{k-2}}$ it holds that
\begin{equation}
\SmallSum{n=0}{k-1} c_n 
= 
c_{k-1} + \SmallSum{n=0}{k-2} c_n 
=
\tfrac{f_{k}-f_{k-1}}{\fx_{k}-\fx_{k-1}}
- 
\tfrac{f_{k-1}-f_{k-2}}{\fx_{k-1}-\fx_{k-2}}
+ 
\tfrac{f_{k-1}-f_{k-2}}{\fx_{k-1}-\fx_{k-2}}
=
\tfrac{f_{k}-f_{k-1}}{\fx_{k}-\fx_{k-1}}.
\end{equation}
Induction thus proves \eqref{interpol_ANN_points:eq1.1}.
In addition, observe that 
\eqref{interpol_ANN_points:eq1}, \eqref{interpol_ANN_points:eq1.1}, and
the fact that $\forall\,k\in\{1,2,\dots,K\} \colon \allowbreak \fx_{k-1} < \fx_{k}$
show that for all $k\in\{1,2,\dots,K\}$, $x\in[\fx_{k-1},\fx_{k}]$ it holds that
\begin{equation}\label{interpol_ANN_points:eq3}
\begin{split}
& (\functionANN{\relu}(\F))(x) - (\functionANN{\relu}(\F))(\fx_{k-1}) 
= \SmallSum{n=0}{K} c_n \left(\max\{x-\fx_n,0\} - \max\{\fx_{k-1}-\fx_n,0\}\right)\\
& = \SmallSum{n=0}{k-1} c_n[(x-\fx_n)-(\fx_{k-1}-\fx_n)] = \SmallSum{n=0}{k-1} c_n(x-\fx_{k-1})
= \bigl(\tfrac{f_{k}-f_{k-1}}{\fx_{k}-\fx_{k-1}}\bigr)(x-\fx_{k-1}).
\end{split}
\end{equation}
Next we claim that for all $k\in\{1,2,\dots,K\}$, $x\in[\fx_{k-1},\fx_{k}]$ it holds that
\begin{equation}\label{interpol_ANN_points:eq4}
(\functionANN{\relu}(\F))(x) = f_{k-1} + \bigl(\tfrac{f_{k}-f_{k-1}}{\fx_{k}-\fx_{k-1}}\bigr)(x-\fx_{k-1}).
\end{equation}
We now prove \cref{interpol_ANN_points:eq4} by induction on $k\in\{1,2,\dots,K\}$.
For the base case $k=1$ observe that \cref{interpol_ANN_points:eq2} and \cref{interpol_ANN_points:eq3} demonstrate that for all $x\in[\fx_0,\fx_1]$ it holds that
\begin{equation}\label{interpol_ANN_points:eq5}
(\functionANN{\relu}(\F))(x) = (\functionANN{\relu}(\F))(\fx_0) + (\functionANN{\relu}(\F))(x) - (\functionANN{\relu}(\F))(\fx_0) = f_0 + \bigl(\tfrac{f_1-f_0}{\fx_1-\fx_0}\bigr)(x-\fx_0).
\end{equation}
This proves \cref{interpol_ANN_points:eq4} in the base case $k=1$.
For the induction step from $\{1,2,\ldots,K-1\} \ni (k-1) \induct k \in \{2,3,\ldots,K\}$ note that \cref{interpol_ANN_points:eq3} implies that for all $k\in\{2,3,\dots,K\}$, $x\in[\fx_{k-1},\fx_{k}]$ with $\forall\,y \in [\fx_{k-2},\fx_{k-1}] \colon (\functionANN{\relu}(\F))(y) = f_{k-2} + (\tfrac{f_{k-1} - f_{k-2}}{\fx_{k-1} - \fx_{k-2}})(y-\fx_{k-2})$ it holds that
\begin{align}\label{interpol_ANN_points:eq6}
& (\functionANN{\relu}(\F))(x) = (\functionANN{\relu}(\F))(\fx_{k-1}) + (\functionANN{\relu}(\F))(x) - (\functionANN{\relu}(\F))(\fx_{k-1}) \\
& = f_{k-2} + \bigl(\tfrac{f_{k-1} - f_{k-2}}{\fx_{k-1} - \fx_{k-2}})(\fx_{k-1} - \fx_{k-2}\bigr) + \bigl(\tfrac{f_{k} - f_{k-1}}{\fx_{k} - \fx_{k-1}}\bigr)(x - \fx_{k-1}) = f_{k-1} + \bigl(\tfrac{f_{k} - f_{k-1}}{\fx_{k} - \fx_{k-1}}\bigr)(x - \fx_{k-1}). \nonumber
\end{align}
Induction thus proves \cref{interpol_ANN_points:eq4}.
Furthermore, observe that 
\cref{interpol_ANN_points:setting1} and \cref{interpol_ANN_points:eq1.1} 
ensure that 
\begin{equation}
\label{TBD}
\SmallSum{n=0}{K} c_n 
= 
c_K + \SmallSum{n=0}{K-1} c_n 
= 
-\tfrac{f_{K}-f_{K-1}}{\fx_{K}-\fx_{K-1}} + \tfrac{f_{K}-f_{K-1}}{\fx_{K}-\fx_{K-1}} 
=
0.
\end{equation}
The fact that $\forall\,k\in\{0,1,\dots,K\} \colon \fx_k \le \fx_K$ and \cref{interpol_ANN_points:eq1} hence imply that for all $x\in[\fx_K,\infty)$ it holds that
\begin{equation}\label{interpol_ANN_points:eq7}
\begin{split}
(\functionANN{\relu}(\F))(x) - (\functionANN{\relu}(\F))(\fx_K) & = \left[\SmallSum{n=0}{K} c_n \left(\max\{x-\fx_n,0\} - \max\{\fx_K - \fx_n,0\}\right)\right]\\
& = \SmallSum{n=0}{K} c_n[(x-\fx_n)-(\fx_K-\fx_n)] = \SmallSum{n=0}{K} c_n(x-\fx_K) = 0.
\end{split}
\end{equation}
This and \cref{interpol_ANN_points:eq4} show that for all $x\in[\fx_K,\infty)$ it holds that
\begin{equation}\label{interpol_ANN_points:eq8}
(\functionANN{\relu}(\F))(x) = (\functionANN{\relu}(\F))(\fx_K) = f_{K-1} + \bigl(\tfrac{f_K - f_{K-1}}{\fx_K - \fx_{K-1}}\bigr)(\fx_K - \fx_{K-1}) = f_K.
\end{equation}
Combining this, \cref{interpol_ANN_points:eq2}, \cref{interpol_ANN_points:eq4}, and \cref{eq:lin_interp} establishes
that $\functionANN{\relu}(\F) = \scrL_{\fx_0,\fx_1,\dots,\fx_K}^{f_0,f_1,\dots,f_K}$.
\end{aproof}


\subsubsection{ANN approximations of one-dimensional functions}\label{sec:ANN_one}

\cfclear
\begin{athm}{lemma}{interpol_ANN_function1}
	Let $K\in\N$, $\beta \in (0,\infty)$,  
	$L, b_1 , \fx_0,\fx_1,\dots,\fx_K\in \R$, $b_2 \in (b_1,\infty)$
	satisfy
	for all $k\in\{0,1,\dots,K\}$ that
	$\fx_k = b_1 + \frac{k(b_2-b_1)}{K}$,
	let $\relu,a\in C(\R,\R)$ satisfy for all $x \in \R$ that $\relu(x)=\max\{x,0\}$ and $a(x)=\frac{1}{\beta} \ln(1+\exp(\beta x))$, 
	let
	$f \colon [b_1,b_2] \to \R$ satisfy for all $x,y\in [b_1,b_2]$ that  $\abs{f(x)-f(y)}\le L\abs{x-y}$,
	and let 
	$\interpolatingDNN \in \ANNs$ satisfy
	\begin{equation}
	\interpolatingDNN 
	= 
	\compANN{
		\AffineANN_{1,f(\fx_0)}
	}{
		\left({\bbigANNsum_{k=0}^K} \left(
		\scalarMultANN{
			\left(\tfrac{K(f(\fx_{\min\{k+1,K\}})-2f(\fx_k) + f(\fx_{\max\{k-1,0\}}))}{(b_1-b_2)}\right)
		}{
			(\compANN{\ii_{1} }{ \AffineANN_{1,-\fx_k}})
		}
		\right)\right)
	}
	\end{equation}
	\cfload.
	Then 
	\begin{enumerate}[label=(\roman *)]
		\item
		\label{interpol_ANN_function1:item3}
		it holds for all $x,y\in\R$ that $\abs{(\functionANN{\relu}(\interpolatingDNN))(x) - (\functionANN{\relu}(\interpolatingDNN))(y)} \le L\abs{x-y}$,
		\item
		\label{interpol_ANN_function1:item4}
		it holds that $\sup_{x\in[b_1,b_2]} \abs{ (\functionANN{\relu}(\interpolatingDNN))(x)-f(x) } \le 
		L(b_2-b_1)K^{-1}$, and 
		\item
		\label{interpol_ANN_function1:item6}
		it holds for all $x,y\in\R$ that $\abs{(\functionANN{a}(\interpolatingDNN))(x) - (\functionANN{a}(\interpolatingDNN))(y)} \le L\abs{x-y}$
	\end{enumerate}
	\cfout.
\end{athm}

\begin{aproof}
	Note that the fact that  
	$\forall \, k\in\{0,1,\dots,K\} \colon \fx_{\min\{k+1,K\}} - \fx_{\min\{k,K-1\}} = \fx_{\max\{k,1\}} - \fx_{\max\{k-1,0\}} = (b_2-b_1)K^{-1}$
	assures that for all
	$k\in\{0,1,\dots,K\}$ 
	it holds that 
	\begin{equation}
	\label{interpol_ANN_function1:eq1}
	\tfrac{(f(\fx_{\min\{k+1,K\}})-f(\fx_k))}{(\fx_{\min\{k+1,K\}}-\fx_{\min\{k,K-1\}})} - \tfrac{(f(\fx_k) - f(\fx_{\max\{k-1,0\}}))}{(\fx_{\max\{k,1\}} - \fx_{\max\{k-1,0\}})} 
	= 
	\tfrac{K(f(\fx_{\min\{k+1,K\}})-2f(\fx_k) + f(\fx_{\max\{k-1,0\}}))}{(b_2-b_1)}.
	\end{equation}
	This and 
	\cref{non_rep3} 
	demonstrate that 
	\begin{equation}
	\label{interpol_ANN_function1:eq2}
	\functionANN{\relu}(\interpolatingDNN)
	=
	\interpol{
		\fx_0,\fx_1,\dots,\fx_K
	}{
		f(\fx_0), f(\fx_1), \ldots, f(\fx_K)
	}.
	\end{equation}
	Combining this with the assumption that $\forall \, x,y\in [b_1,b_2] \colon \abs{f(x)-f(y)}\le L\abs{x-y}$ and \cref{lin_interp1_item1_new2} in \cref{lin_interp1_new2} establishes \cref{interpol_ANN_function1:item3}.
	Moreover, note that 
	\eqref{interpol_ANN_function1:eq2}, 
	the assumption that $\forall \, x,y\in [b_1,b_2] \colon \abs{f(x)-f(y)}\le L\abs{x-y}$,
	\cref{lin_interp1_item2_new2} in \cref{lin_interp1_new2}, and 
	the fact that  $\forall \, k\in\{1,2,\dots,K\} \colon \fx_{k} - \fx_{k-1} = (b_2-b_1)K^{-1}$ 
	demonstrate that for all 
	$x\in[b_1,b_2]$ 
	it holds that
	\begin{equation}
	\abs{ (\functionANN{\relu}(\interpolatingDNN))(x) - f(x) } \le L\left(\max_{k\in\{1,2,\dots,K\}} \abs{\fx_k - \fx_{k-1} } \right) = L(b_2-b_1)K^{-1}.
	\end{equation}
	This establishes \cref{interpol_ANN_function1:item4}. 
	Next, observe that \cref{non_rep2_item3} in \cref{non_rep2} 
	and \cref{interpol_ANN_function1:eq1} 
	imply for all $x \in \R$ that 
	\begin{equation}\label{interpol_ANN_function1:eq3}
	\begin{split}
	(\functionANN{a}(\interpolatingDNN))(x) 
	& = f(\fx_0) + \sum_{k=0}^K 
	\left( \tfrac{(f(\fx_{\min\{k+1,K\}})-f(\fx_k))}{(\fx_{\min\{k+1,K\}}-\fx_{\min\{k,K-1\}})} - \tfrac{(f(\fx_k) - f(\fx_{\max\{k-1,0\}}))}{(\fx_{\max\{k,1\}} - \fx_{\max\{k-1,0\}})}  \right)
	a(x-\fx_k) \\
	& = f(\fx_0) + \sum_{k=0}^{K-1} \frac{f(\fx_{k+1})-f(\fx_k)}{\fx_{k+1}-\fx_k} a(x-\fx_k) 
	- \sum_{k=1}^{K} \frac{f(\fx_{k})-f(\fx_{k-1})}{\fx_{k}-\fx_{k-1}} a(x-\fx_k) \\
	& = f(\fx_0) + \sum_{k=1}^{K} \frac{f(\fx_{k})-f(\fx_{k-1})}{\fx_{k}-\fx_{k-1}} \left( a(x-\fx_{k-1}) - a(x-\fx_k) \right) .
	\end{split}
	\end{equation}
	Note that $a$ is differentiable and it holds for all $x\in\R$ that $\frac{da(x)}{dx} = \frac{e^{\beta x}}{1+ e^{\beta x}}$. 
	It thus follows from \cref{interpol_ANN_function1:eq3} that $\functionANN{a}(\interpolatingDNN)$ is differentiable and it holds for all $x\in\R$ that 
	\begin{equation}
	\begin{split}
	\frac{d(\functionANN{a}(\interpolatingDNN))(x) }{dx}
	& = \sum_{k=1}^K \frac{f(\fx_{k})-f(\fx_{k-1})}{\fx_{k}-\fx_{k-1}} 
	\left( \frac{e^{\beta (x-\fx_{k-1})}}{1+ e^{\beta (x-\fx_{k-1})}} - \frac{e^{\beta (x-\fx_{k})}}{1+ e^{\beta (x-\fx_{k})}} \right) .
	\end{split}
	\end{equation}
	This, the triangle inequality, the assumption that $\forall \, x,y\in [b_1,b_2] \colon \abs{f(x)-f(y)}\le L\abs{x-y}$, and the fact that for all $k\in \{1,\ldots,K\}$, $x\in\R$ it holds that 
	\begin{equation}
	\frac{e^{\beta (x-\fx_{k-1})}}{1+ e^{\beta (x-\fx_{k-1})}} - \frac{e^{\beta (x-\fx_{k})}}{1+ e^{\beta (x-\fx_{k})}} 
	= \frac{ e^{\beta (x-\fx_{k-1})} - e^{\beta (x-\fx_{k})} }{\left( 1+ e^{\beta (x-\fx_{k-1})} \right) \left( 1+ e^{\beta (x-\fx_{k})} \right)}
	\ge 0
	\end{equation}
	yield for all $x \in \R$ that 
	\begin{equation}
	\begin{split}
	\abs[\bigg]{ \frac{d(\functionANN{a}(\interpolatingDNN))(x) }{dx} }
	& \le \sum_{k=1}^K \frac{\abs{f(\fx_{k})-f(\fx_{k-1})}}{\abs{\fx_{k}-\fx_{k-1}}} 
	\left( \frac{e^{\beta (x-\fx_{k-1})}}{1+ e^{\beta (x-\fx_{k-1})}} - \frac{e^{\beta (x-\fx_{k})}}{1+ e^{\beta (x-\fx_{k})}} \right) \\
	& \le L \left( \sum_{k=1}^K  \frac{e^{\beta (x-\fx_{k-1})}}{1+ e^{\beta (x-\fx_{k-1})}} - \sum_{k=1}^K  \frac{e^{\beta (x-\fx_{k})}}{1+ e^{\beta (x-\fx_{k})}} \right) \\
	& = L \left( \frac{e^{\beta (x-\fx_{0})}}{1+ e^{\beta (x-\fx_{0})}} - \frac{e^{\beta (x-\fx_{K})}}{1+ e^{\beta (x-\fx_{K})}} \right) 
	\le \frac{Le^{\beta (x-\fx_{0})}}{1+ e^{\beta (x-\fx_{0})}} 
	\le L .
	\end{split}
	\end{equation}
	Hence, it holds for all $x,y \in \R$ that $\abs{ (\functionANN{a}(\interpolatingDNN))(x) - (\functionANN{a}(\interpolatingDNN))(y) }
	\le L \abs{x-y}$.
	This establishes \cref{interpol_ANN_function1:item6}.
\end{aproof}

\cfclear
\begin{athm}{lemma}{interpol_ANN_ReLU}
	Let 
	$\varepsilon \in (0,1]$, $L\in[0,\infty)$, $q\in(1,\infty)$, 
	let $b\in [1,\infty)$ satisfy $\max\{1,2L\}=\varepsilon b^{q-1}$, 
	let $K\in\N \cap [\frac{2Lb}{\varepsilon},\frac{2Lb}{\varepsilon}+1]$, 
	let $f\colon \R\to\R$ satisfy for all $x,y\in\R$ that $\abs{f(x)-f(y)} \le L\abs{x-y}$, 
	let $\fx_0,\fx_1,\dots,\fx_K, c_0,c_1,\ldots,c_K \in\R$ satisfy for all $k\in\{0,1,\ldots,K\}$ that 
	$\fx_k = -b + \frac{2kb}{K}$ and 
	\begin{equation}
	c_k
	=
	\frac{K(f(\fx_{\min\{k+1,K\}})-2f(\fx_k) + f(\fx_{\max\{k-1,0\}}))}{2b} ,
	\end{equation}
	let $\interpolatingDNN\in\ANNs$ satisfy that 
	\begin{equation}
	\interpolatingDNN 
	= 
	\compANN{
		\AffineANN_{1,f(\fx_0)}
	}{
		\left({\bbigANNsum_{k=0}^K} \left(
		\scalarMultANN{
			c_k
		}{
			(\compANN{\ii_{1} }{ \AffineANN_{1,-\fx_k}})
		}
		\right)\right)
	},
	\end{equation}
	and let $\relu\in C(\R,\R)$ satisfy for all $x\in\R$ that $\relu(x)=\max\{x,0\}$
	\cfload. 
	Then 
	\begin{enumerate}[label=(\roman *)]
		\item 
		\label{interpol_ANN_ReLU:item4}
		it holds for all $x,y\in\R$ that $\abs{(\functionANN{\relu}(\interpolatingDNN))(x) - (\functionANN{\relu}(\interpolatingDNN))(y)} \le L\abs{x-y}$,
		\item 
		\label{interpol_ANN_ReLU:item5}
		it holds that $\sup_{x\in[-b,b]} \abs{(\functionANN{\relu}(\interpolatingDNN))(x)-f(x)} \le \frac{2Lb}{K} \le \varepsilon$,
		\item 
		\label{interpol_ANN_ReLU:item6}
		it holds for all $x\in\R$ that $\abs{(\functionANN{\relu}(\interpolatingDNN))(x)-f(x)} \le \varepsilon \max\{1,\abs{x}^q\}$, 
		\item
		\label{interpol_ANN_ReLU:item7}
		it holds that $\dimANNlevel_1(\interpolatingDNN) \le 2 (\max\{1,2L\})^{\nicefrac{q}{(q-1)}} \varepsilon^{-\nicefrac{q}{(q-1)}} + 1$, and
		\item 
		\label{interpol_ANN_ReLU:item8}
		it holds that 
		$\paramANN(\interpolatingDNN) = 3(\dimANNlevel_1(\interpolatingDNN)) + 1 \le 12 (\max\{1,2L\})^{\nicefrac{q}{(q-1)}} \varepsilon^{-\nicefrac{q}{(q-1)}} $
	\end{enumerate}
	\cfout.
\end{athm}

\begin{aproof}	
	Note that \cref{interpol_ANN_function1:item3}
	in \cref{interpol_ANN_function1} 
	yields 
	\cref{interpol_ANN_ReLU:item4}.  
	Next, observe that the fact that
	$K \in \N \cap [\frac{2Lb}{\varepsilon}, \frac{2Lb}{\varepsilon} + 1]$
	implies that $\frac{2Lb}{K} \le \varepsilon$.
	This and 
	\cref{interpol_ANN_function1:item4}
	in \cref{interpol_ANN_function1} 
	establish 
	\cref{interpol_ANN_ReLU:item5}. 	
	The triangle inequality, \cref{interpol_ANN_ReLU:item4}, the fact that $f(-b)=(\functionANN{\relu}(\interpolatingDNN))(-b)$, the fact that $f(b)=(\functionANN{\relu}(\interpolatingDNN))(b)$, and the fact that for all $x,y\in\R$ it holds that $\abs{f(x)-f(y)}\le L\abs{x-y}$ ensure that for all $x\in\R$ it holds that
	\begin{equation}\label{int_bd_pre1}
	\begin{split}
	\abs{(\functionANN{\relu}(\interpolatingDNN))(x)-f(x)}
	& \le \abs{(\functionANN{\relu}(\interpolatingDNN))(x) - f(b)} + \abs{f(b) -f(0)} + \abs{f(0)-f(x)} \\
	& = \abs{(\functionANN{\relu}(\interpolatingDNN))(x) - (\functionANN{\relu}(\interpolatingDNN))(b)} + \abs{f(b) -f(0)} + \abs{f(0)-f(x)} \\
	& \le L\abs{x-b} + L\abs{b} + L\abs{x} = L\pr{ \abs{x-b} + b + \abs{x} }
	\end{split}
	\end{equation}
	and
	\begin{equation}\label{int_bd_pre2}
	\begin{split}
	\abs{(\functionANN{\relu}(\interpolatingDNN))(x)-f(x)}
	& \le \abs{(\functionANN{\relu}(\interpolatingDNN))(x) - f(-b)} + \abs{f(-b) -f(0)} + \abs{f(0)-f(x)} \\
	& = \abs{(\functionANN{\relu}(\interpolatingDNN))(x) - (\functionANN{\relu}(\interpolatingDNN))(-b)} + \abs{f(-b) -f(0)} + \abs{f(0)-f(x)} \\
	& \le L\abs{x+b} + L\abs{b} + L\abs{x} = L\pr{ \abs{x+b} + b + \abs{x} } . 
	\end{split}
	\end{equation}
	It follows from \cref{int_bd_pre1} that for all $x\in(b,\infty)$ it holds that
	\begin{equation}\label{int_bd_pre3}
	\begin{split}
	\frac{ \abs{(\functionANN{\relu}(\interpolatingDNN))(x)-f(x)} }{ \max\{ 1 , \abs{x}^q \} }
	& \le \frac{ L \pr{ \abs{x-b} + b + \abs{x} } }{ \max\{ 1, \abs{x}^q \} }
	= \frac{ L \pr{ x-b + b + x } }{ \max\{ 1, \abs{x}^q \} } \\
	& = \frac{ 2L \abs{x} }{ \max\{ 1, \abs{x}^q \} }
	\le \frac{ \max\{1,2L\} }{ \abs{x}^{q-1} } 
	\le \frac{ \max\{1,2L\} }{ b^{q-1} } = \varepsilon .
	\end{split}
	\end{equation}
	Moreover, \cref{int_bd_pre2} demonstrates that for all $x \in (-\infty,-b)$ it holds that
	\begin{equation}
	\begin{split}
	\frac{ \abs{(\functionANN{\relu}(\interpolatingDNN))(x)-f(x)} }{ \max\{ 1 , \abs{x}^q \} }
	& \le \frac{ L \pr{ \abs{x+b} + b + \abs{x} } }{ \max\{ 1, \abs{x}^q \} }
	= \frac{ L \pr{ -(x+b) + b -x } }{ \max\{ 1, \abs{x}^q \} } \\
	& = \frac{ 2L \abs{x} }{ \max\{ 1, \abs{x}^q \} }
	\le \frac{ \max\{1,2L\} }{ \abs{x}^{q-1} } 
	\le \frac{ \max\{1,2L\} }{ b^{q-1} } = \varepsilon .
	\end{split}
	\end{equation}
	Combining this, \cref{int_bd_pre3}, and  \cref{interpol_ANN_ReLU:item5} shows that for all $x \in \R$ it holds that 
	$\abs{(\functionANN{\relu}(\interpolatingDNN))(x)-f(x)} \le \varepsilon \max\{1,\abs{x}^q\}$.
	This establishes \cref{interpol_ANN_ReLU:item6}.
	In addition, observe that 
	the fact that $\max\{1,2L\} = \varepsilon b^{q-1}$  
	and the fact that
	$K \le 1 + \frac{2Lb}{\varepsilon}$
	prove that
	\begin{equation}\label{interpol_ANN_function_implicit:eq3}
	K \le 1 + \frac{2Lb}{\varepsilon}
	\le 1 + \frac{\max\{1,2L\}b}{\varepsilon}
	= 1 + b^q 
	\le 2 b^q
	= 2 \pr*{ \frac{\max\{1,2L\}}{\varepsilon} }^{\!\!\nicefrac{q}{(q-1)}} .
	\end{equation}
	This and the fact that $\dimANNlevel_1(\interpolatingDNN)=K+1$ (cf.\ \cref{non_rep2_item1} in \cref{non_rep2}) 
	establish \cref{interpol_ANN_ReLU:item7}. 
	Moreover, observe that 
	\cref{non_rep2_item4} in \cref{non_rep2} 
	and \cref{interpol_ANN_ReLU:item7} 	
	ensure that
	\begin{align}\label{interpol_ANN_function_implicit:eq2a}
	\paramANN(\interpolatingDNN) 
	= 
	3(\dimANNlevel_1(\interpolatingDNN)) + 1
	\le 4(\dimANNlevel_1(\interpolatingDNN))
	& \le 8 (\max\{1,2L\})^{\nicefrac{q}{(q-1)}} \varepsilon^{-\nicefrac{q}{(q-1)}} + 4 \\
	& \le 12 (\max\{1,2L\})^{\nicefrac{q}{(q-1)}} \varepsilon^{-\nicefrac{q}{(q-1)}} . \nonumber
	\end{align}
	This establishes \cref{interpol_ANN_ReLU:item8}. 
\end{aproof}

\cfclear
\begin{athm}{corollary}{interpol_ANN_function_implicit_max}
	Let $\varepsilon \in (0,1]$, 
	$L \in [0,\infty)$, $q \in (1,\infty)$, 
	$\alpha \in [0,\infty)\backslash\{1\}$,  
	let
	$f \colon \R \to \R$ satisfy for all $x,y\in \R$ that  $\abs{ f(x)-f(y) } \le L \abs{ x-y }$, 
	and let $a\in C(\R,\R)$ satisfy for all $x\in\R$ that $a(x) = \max\{x,\alpha x\}$. 
	Then there exists $\G \in\ANNs$ such that
	\begin{enumerate}[label=(\roman *)]
		\item 
		\label{interpol_ANN_function_implicit_max:item1}
		it holds that $\functionANN{a}(\G) \in C(\R,\R)$, 
		\item
		\label{interpol_ANN_function_implicit_max:item2}
		it holds that 
		$\hiddenLength(\G) = 1$, 
		\item
		\label{interpol_ANN_function_implicit_max:item3}
		it holds for all $x,y\in\R$ that $\abs{(\functionANN{a}(\G))(x) - (\functionANN{a}(\G))(y)} \le L\abs{x-y}$,
		\item
		\label{interpol_ANN_function_implicit_max:item4}
		it holds for all $x\in\R$ that $\abs{(\functionANN{a}(\G))(x)-f(x)} \le \varepsilon \max\{1,\abs{x}^q\}$, 
		\item
		\label{interpol_ANN_function_implicit_max:item5}
		it holds that $\dimANNlevel_1(\G) \le 4 (\max\{1,2L\})^{\nicefrac{q}{(q-1)}} \varepsilon^{-\nicefrac{q}{(q-1)}} + 2$, and
		\item
		\label{interpol_ANN_function_implicit_max:item6}
		it holds that 
		$\paramANN(\G) = 3(\dimANNlevel_1(\G)) + 1 \le 24  (\max\{1,2L\})^{\nicefrac{q}{(q-1)}} \varepsilon^{-\nicefrac{q}{(q-1)}} $
	\end{enumerate}
	\cfout.
\end{athm}

\begin{aproof}
	Throughout this proof 
	let $b \in [1,\infty)$ satisfy $\max\{1,2L\} = \varepsilon b^{q-1}$, 
	let
	$K \in \N \cap [\frac{2L b}{\varepsilon}, \frac{2L b}{\varepsilon} + 1]$,  $\fx_0,\fx_1,\dots,\fx_K$, $c_0,c_1,\dots,c_K$, $h_0,h_1,\ldots,h_{2K+1}$, $\alpha_0,\alpha_1,\ldots,\alpha_{2K+1}$, $\beta_0,\beta_1,\ldots,\beta_{2K+1}\in\R$
	satisfy for all $k\in\{0,1,\dots,K\}$ that
	$\fx_k = -b + \frac{2kb}{K}$, 
	\begin{equation}
	c_k
	=
	\frac{K(f(\fx_{\min\{k+1,K\}})-2f(\fx_k) + f(\fx_{\max\{k-1,0\}}))}{2b} ,
	\end{equation} 
	$h_k = 
	\frac{c_k \abs{1-\alpha} \alpha}{(1-\alpha)(1-\alpha^2)}$, 
	$	h_{k+K+1} = 
	\frac{c_k \abs{1-\alpha}}{(1-\alpha)(1-\alpha^2)}$, 
	$\alpha_k = 
	\frac{-\abs{1-\alpha}}{1-\alpha}$, 
	$\alpha_{k+K+1} = 
	\frac{\abs{1-\alpha}}{1-\alpha}$, 
	$\beta_k = 
	\frac{\abs{1-\alpha} \fx_k}{1-\alpha}$,  
	and 
	$\beta_{k+K+1} = 
	\frac{-\abs{1-\alpha} \fx_{k}}{1-\alpha}$, 
	let $\interpolatingDNN\in\ANNs$ satisfy that
	\begin{equation}
	\interpolatingDNN 
	= 
	\compANN{
		\AffineANN_{1,f(\fx_0)}
	}{
		\left({\bbigANNsum_{k=0}^K} \left(
		\scalarMultANN{
			c_k
		}{
			(\compANN{\ii_{1} }{ \AffineANN_{1,-\fx_k}})
		}
		\right)\right), 
	}
	\end{equation}
	let  $\G\in\ANNs$ satisfy that
	\begin{equation}
	\G
	= 
	\compANN{
		\AffineANN_{1,f(\fx_0)}
	}{
		\left({\bbigANNsum_{k=0}^{2K+1}} \left(
		\scalarMultANN{
			h_k
		}{
			(\compANN{\ii_{1} }{ \AffineANN_{\alpha_k,\beta_k}})
		}
		\right)\right), 
	}
	\end{equation}
	and let $\relu \in C(\R,\R)$ satisfy for all $x\in\R$ that $\relu(x)=\max\{x,0\}$
	\cfload. 
	Note that \cref{non_rep2_item1,non_rep2_item2} in \cref{non_rep2} establish 
	\cref{interpol_ANN_function_implicit_max:item1,interpol_ANN_function_implicit_max:item2}. 
	Furthermore, observe that 
	the fact that for all $x \in \R$ it holds that $a(x)=\max\{x,\alpha x\}=\frac{x+\alpha x+\abs{x-\alpha x}}{2}$
	implies for all $x\in\R$ that 	
	\begin{equation}
	\begin{split}
	& \frac{\abs{1-\alpha}\alpha}{(1-\alpha)(1-\alpha^2)} 
	a\left( - \frac{\abs{1-\alpha}x}{1-\alpha} \right)
	+ \frac{\abs{1-\alpha}}{(1-\alpha)(1-\alpha^2)} 
	a\left( \frac{\abs{1-\alpha}x}{1-\alpha} \right) \\
	& = \frac{\abs{1-\alpha}\alpha}{(1-\alpha)(1-\alpha^2)} 
	\frac{ \frac{-\abs{1-\alpha}x}{1-\alpha} + \frac{-\alpha \abs{1-\alpha} x}{1-\alpha} + \abs[\big]{  \frac{-\abs{1-\alpha}x}{1-\alpha} - \frac{-\alpha\abs{1-\alpha}x}{1-\alpha}  } }{2} \\
	& \quad + \frac{\abs{1-\alpha}}{(1-\alpha)(1-\alpha^2)} 
	\frac{ \frac{\abs{1-\alpha}x}{1-\alpha} + \frac{\alpha \abs{1-\alpha} x}{1-\alpha} + \abs[\big]{  \frac{\abs{1-\alpha}x}{1-\alpha} - \frac{\alpha\abs{1-\alpha}x}{1-\alpha}  } }{2} \\
	& = \frac{ \abs{1-\alpha} \Big( -\alpha \abs{1-\alpha} x - \alpha^2 \abs{1-\alpha} x + \abs{1-\alpha} x + \alpha \abs{1-\alpha} x \Big) }{ 2(1-\alpha)(1-\alpha^2)(1-\alpha) } \\
	& \quad + \frac{ \abs{1-\alpha} \Big( \alpha \abs[\big]{ -\abs{1-\alpha} x + \alpha \abs{1-\alpha} x } + \abs[\big]{ \abs{1-\alpha} x - \alpha \abs{1-\alpha} x } \Big) }{ 2(1-\alpha)(1-\alpha^2)\abs{1-\alpha} } \\
	& = \frac{\abs{1-\alpha}^2 (1-\alpha^2) x}{2(1-\alpha)^2(1-\alpha^2)} 
	+ \frac{ \alpha \abs[\big]{  (-1+\alpha) \abs{1-\alpha} x } + \abs[\big]{ (1-\alpha) \abs{1-\alpha} x } }{ 2(1-\alpha)(1-\alpha^2) } \\
	& = \frac{x}{2} + \frac{ (\alpha +1) \abs{1-\alpha} \abs{1-\alpha} \abs{x}  }{ 2(1-\alpha)(1-\alpha)(1+\alpha) }
	= \frac{x+\abs{x}}{2} 
	= \max\{x,0\} 
	= \relu(x) .
	\end{split}
	\end{equation}
	Combining this and \cref{non_rep2_item3} in \cref{non_rep2} shows for all $x \in \R$ that 
	\begin{equation}
	\begin{split}
	(\functionANN{a}(\G))(x) 
	& = f(\fx_0) + \sum_{k=0}^K  
	\frac{c_k \abs{1-\alpha} \alpha}{(1-\alpha)(1-\alpha^2)} 
	\left( a\left( -\frac{\abs{1-\alpha} x}{1-\alpha} + \frac{\abs{1-\alpha} \fx_k}{1-\alpha} \right) \right) \\
	& \quad + \sum_{k=0}^K 
	\frac{c_k \abs{1-\alpha} }{(1-\alpha)(1-\alpha^2)}
	\left( a\left( \frac{\abs{1-\alpha} x}{1-\alpha} - \frac{\abs{1-\alpha} \fx_k}{1-\alpha} \right) \right) \\
	& = f(\fx_0) + \sum_{k=0}^K c_k \relu(x-\fx_k) 
	= 
	(\functionANN{\relu}(\interpolatingDNN))(x)
	.
	\end{split}
	\end{equation}
	Therefore, 
	\cref{interpol_ANN_ReLU:item4,interpol_ANN_ReLU:item6} 
	in \cref{interpol_ANN_ReLU}
	establish 
	\cref{interpol_ANN_function_implicit_max:item3,interpol_ANN_function_implicit_max:item4}.
	Note that 
	\cref{non_rep2_item1} in \cref{non_rep2} 
	shows that 
	$\dimANNlevel_1(\G) = 2(K+1) = 2 \dimANNlevel_1(\interpolatingDNN)$.  
	Therefore, \cref{interpol_ANN_ReLU:item7} in \cref{interpol_ANN_ReLU} implies \cref{interpol_ANN_function_implicit_max:item5}. 
	It follows from 
	\cref{non_rep2_item4} in \cref{non_rep2} 
	that 
	$\paramANN(\G) 
	= 3(\dimANNlevel_1(\G)) + 1 
	\leq 2 (3(\dimANNlevel_1(\interpolatingDNN)) + 1) = 2\paramANN(\interpolatingDNN)$.  
	This and \cref{interpol_ANN_ReLU:item8} in \cref{interpol_ANN_ReLU} ensure \cref{interpol_ANN_function_implicit_max:item6}. 
\end{aproof}

\cfclear
\begin{athm}{corollary}{interpol_ANN_function_implicit_softplus}
	Let $\varepsilon \in (0,1]$, 
	$L \in [0,\infty)$, $q \in (1,\infty)$, 
	let
	$f \colon \R \to \R$ satisfy for all $x,y\in \R$ that  $\abs{ f(x)-f(y) } \le L \abs{ x-y }$, 
	and let $a \in C(\R,\R)$ satisfy for all $x\in\R$ that $a(x) = \ln(1+\exp(x))$. 
	Then there exists $\G \in\ANNs$ such that
	\begin{enumerate}[label=(\roman *)]
		\item 
		\label{interpol_ANN_function_implicit_softplus:item1}
		it holds that $\functionANN{a}(\G) \in C(\R,\R)$, 
		\item
		\label{interpol_ANN_function_implicit_softplus:item2}
		it holds that 
		$\hiddenLength(\G) = 1$, 
		\item
		\label{interpol_ANN_function_implicit_softplus:item3}
		it holds for all $x,y\in\R$ that $\abs{(\functionANN{a}(\G))(x) - (\functionANN{a}(\G))(y)} \le L\abs{x-y}$,
		\item
		\label{interpol_ANN_function_implicit_softplus:item4}
		it holds for all $x\in\R$ that $\abs{(\functionANN{a}(\G))(x)-f(x)} \le 2 \varepsilon \max\{1,\abs{x}^q\}$, 
		\item
		\label{interpol_ANN_function_implicit_softplus:item5}
		it holds that $\dimANNlevel_1(\G) \le 2 (\max\{1,2L\})^{\nicefrac{q}{(q-1)}} \varepsilon^{-\nicefrac{q}{(q-1)}} + 1$, and 
		\item
		\label{interpol_ANN_function_implicit_softplus:item6}
		it holds that 
		$\paramANN(\G) = 3(\dimANNlevel_1(\G)) + 1 \le 12 (\max\{1,2L\})^{\nicefrac{q}{(q-1)}} \varepsilon^{-\nicefrac{q}{(q-1)}} $
	\end{enumerate}
	\cfout.
\end{athm}

\begin{aproof}
	Throughout this proof 
	let $b \in [1,\infty)$ satisfy $\max\{1,2L\} = \varepsilon b^{q-1}$, 
	let
	$K \in \N \cap [\frac{2L b}{\varepsilon}, \frac{2L b}{\varepsilon} + 1]$, 
	let $\fx_0,\fx_1,\dots,\fx_K,c_0,c_1,\dots,c_K\in\R$
	satisfy for all $k\in\{0,1,\dots,K\}$ that
	$\fx_k = -b + \frac{2kb}{K}$
	and
	\begin{equation}\label{eq:interpol_ANN_function_implicit_softplus:defck}
	c_k
	=
	\frac{K(f(\fx_{\min\{k+1,K\}})-2f(\fx_k) + f(\fx_{\max\{k-1,0\}}))}{2b} ,
	\end{equation} 
	let $\interpolatingDNN\in\ANNs$ satisfy that
	\begin{equation}
	\interpolatingDNN
	= 
	\compANN{
		\AffineANN_{1,f(\fx_0)}
	}{
		\left({\bbigANNsum_{k=0}^K} \left(
		\scalarMultANN{
			c_k
		}{
			(\compANN{\ii_{1} }{ \AffineANN_{1,-\fx_k}})
		}
		\right)\right) , 
	}
	\end{equation}
	let $\beta=\max\{2,2K^2L\ln(2)\varepsilon^{-1}\}$, let $\relu,\fa \in C(\R,\R)$ satisfy for all $x \in \R$ that 
	$\relu(x)=\max\{x,0\}$ and 
	\begin{equation}
	\fa(x) = \frac{1}{\beta} \ln\left( 1 + \exp\left( \beta x \right) \right) ,
	\end{equation}
	and let 
	$\G\in\ANNs$ satisfy that
	\begin{equation}
	\G  
	= 
	\compANN{
		\AffineANN_{1,f(\fx_0)}
	}{
		\left({\bbigANNsum_{k=0}^K} \left(
		\scalarMultANN{
			\frac{c_k}{\beta}
		}{
			(\compANN{\ii_{1} }{ \AffineANN_{\beta,-\beta\fx_k}})
		}
		\right)\right) 
	}
	\end{equation}
	\cfload. 
	Note that \cref{non_rep2_item1,non_rep2_item2} in \cref{non_rep2} establish 
	\cref{interpol_ANN_function_implicit_softplus:item1,interpol_ANN_function_implicit_softplus:item2}. 
	Moreover, 
	\cref{non_rep2_item3} in \cref{non_rep2} 
	implies 
	for all $x\in\R$ 
	that 
	\begin{equation}\label{eq:interpol_ANN_function_implicit_softplus:equalfctofANN}
	\begin{split}
	(\functionANN{a}(\G))(x) 
	& = f(\fx_0) +\sum_{k=0}^K  
	\frac{c_k}{\beta}
	\left( a\left( \beta x - \beta \fx_k \right) \right) 
	= f(\fx_0) +\sum_{k=0}^K  
	\frac{c_k}{\beta} \ln\left( 1 + \exp(\beta (x - \fx_k)) \right) \\
	& = f(\fx_0) +\sum_{k=0}^K  
	c_k \left( \fa \left( x - \fx_k \right) \right) 
	= (\functionANN{\fa}(\interpolatingDNN))(x)
	.
	\end{split}
	\end{equation}
	This and 
	\cref{interpol_ANN_function1:item6} in \cref{interpol_ANN_function1} 
	prove \cref{interpol_ANN_function_implicit_softplus:item3}. 
	Observe that the triangle inequality and  
	\cref{interpol_ANN_ReLU:item6} in \cref{interpol_ANN_ReLU}
	show 
	for all $x \in \R$ that 
	\begin{equation}\label{eq:interpol_ANN_function_implicit_softplus:diff1}
	\begin{split}
	\abs{ (\functionANN{a}(\G))(x)  - f(x) } 
	& \le \abs{ (\functionANN{a}(\G))(x) - (\functionANN{\relu}(\interpolatingDNN))(x) } + \abs{ (\functionANN{\relu}(\interpolatingDNN))(x) - f(x) } \\
	& \le  \abs{ (\functionANN{a}(\G))(x) - (\functionANN{\relu}(\interpolatingDNN))(x) } + \varepsilon \max\{1,\abs{x}^q\} .
	\end{split}
	\end{equation}
	Furthermore, 	
	\cref{eq:interpol_ANN_function_implicit_softplus:equalfctofANN}, 
	\cref{non_rep2_item3} in \cref{non_rep2},  
	and the triangle inequality 
	imply 
	for all $x\in\R$ 
	that 
	\begin{equation}\label{eq:interpol_ANN_function_implicit_softplus:diff2}
	\begin{split}
	\abs{ (\functionANN{a}(\G))(x) - (\functionANN{\relu}(\interpolatingDNN))(x) } 
	& = \abs{ (\functionANN{\fa}(\interpolatingDNN))(x) - (\functionANN{\relu}(\interpolatingDNN))(x) } \\
	& \le
	\sum_{k=0}^K \abs{ c_k } \, \abs{ \fa(x-\fx_k) - \relu(x-\fx_k) } .
	\end{split}
	\end{equation}
	Note that it holds for all $x \in [0,\infty)$ that 
	$\abs{\frac{1}{\beta}\ln(1+\exp(\beta x))-x} = \abs{\frac{1}{\beta}\ln(\frac{1+\exp(\beta x)}{\exp(\beta x)})} \le \frac{1}{\beta} \ln(2)$ 
	and that it holds for all $x\in(-\infty,0)$ that 
	$\abs{\frac{1}{\beta}\ln(1+\exp(\beta x))} \le \frac{1}{\beta} \ln(2)$. 
	This and \cref{eq:interpol_ANN_function_implicit_softplus:diff2} 
	yield for all $x\in\R$ that 
	\begin{equation}\label{eq:interpol_ANN_function_implicit_softplus:diff3}
	\abs{ (\functionANN{a}(\G))(x) - (\functionANN{\relu}(\interpolatingDNN))(x) } 
	\le \frac{\ln(2)}{\beta} \sum_{k=0}^K \abs{ c_k } .
	\end{equation}
	Moreover,  \cref{eq:interpol_ANN_function_implicit_softplus:defck}, the triangle inequality,  
	and the assumption that $\forall x,y \in \R\colon \abs{ f(x)-f(y) } \le L \abs{ x-y }$ 
	imply 
	for all $k\in\{0,1,\ldots,K\}$ that 
	\begin{equation}
	\begin{split}
	\abs{c_k} 
	& \le \tfrac{K}{2b} \left( \abs{ f(\fx_{\min\{k+1,K\}})-f(\fx_k) } + \abs{ f(\fx_{\max\{k-1,0\}}) - f(\fx_k) } \right) \\
	& \le \tfrac{KL}{2b} \left( \abs{ \fx_{\min\{k+1,K\}}-\fx_k } + \abs{ \fx_{\max\{k-1,0\}} - \fx_k } \right) 
	\le 2KL .
	\end{split}
	\end{equation}
	This, \cref{eq:interpol_ANN_function_implicit_softplus:diff3}, and the fact that 
	$\beta\ge 2K^2L\ln(2)\varepsilon^{-1}$ 
	show for all $x \in \R$ that 
	\begin{equation}
	\abs{ (\functionANN{a}(\G))(x) - (\functionANN{\relu}(\interpolatingDNN))(x) } 
	\le \tfrac{2\ln(2) K^2 L}{\beta} 
	\le \varepsilon .
	\end{equation}
	Combining this and \cref{eq:interpol_ANN_function_implicit_softplus:diff1} 
	proves 
	\cref{interpol_ANN_function_implicit_softplus:item4}.
	Moreover, observe that 
	\cref{non_rep2_item1,non_rep2_item4} in \cref{non_rep2} 
	show that $\dimANNlevel_1(\G)=K+1=\dimANNlevel_1(\interpolatingDNN)$ 
	and 
	$\paramANN(\G) = 3(\dimANNlevel_1(\G))+1 
	= 3(\dimANNlevel_1(\interpolatingDNN))+1 
	= \paramANN(\interpolatingDNN)$.
	Therefore, 
	\cref{interpol_ANN_ReLU:item7,interpol_ANN_ReLU:item8} in \cref{interpol_ANN_ReLU} 
	establish 
	\cref{interpol_ANN_function_implicit_softplus:item5,interpol_ANN_function_implicit_softplus:item6}.
\end{aproof}

\subsection{ANN approximation results with specific activation functions}\label{sec:4_2}

\cfclear
\begin{athm}{corollary}{cor:almostfinal}
	Let 
	$\gamma,T,\constantAssumpMainThm,\fc \in(0,\infty)$,  
	$r\in\N$, 
	$p\in\N\backslash\{1\}$, 
	$\fq\in [2,\infty)$, 
	let $f\colon \R\to\R$ be Lipschitz continuous, 
	for every $d\in\N$ let $u_d \in C^{1,2}([0,T]\times \R^d,\R)$ 
	satisfy for all $t \in [0,T]$, $x\in\R^d$ that 
	\begin{equation}\label{final_cor_eq1}
	\textstyle{(\tfrac{\partial}{\partial  t}u_d)(t,x) 
	+ 
	\fc (\Delta_x u_d)(t,x) 
	+ 
	f\pr[\big]{ u_d(t,x) } 
	= 
	0}, 
	\end{equation}
	let 
	$\indexAct\in\{0,1\}$, 
	$\alpha \in [0,\infty)\backslash\{1\}$,   
	$\fa_0,\fa_1\in C(\R,\R)$ satisfy for all $x\in\R$ that $\fa_0(x) = \max\{x,\alpha x\}$ and $\fa_1(x)=\ln(1+\exp(x))$,  
	for every $d\in\N$ let $\mu_d\colon \cB(\R^d)\to[0,1]$ be a probability measure
	with 
	\begin{equation}
		\textstyle{\int_{\R^d} \norm{y}^{p^2\fq} \mu_{d}(\dxx y) \le \constantAssumpMainThm d^{r p^2\fq}},
	\end{equation}	
	and assume for all $d\in\N$, $\varepsilon \in (0,1]$ 
	that there exists 
	$\ANNassumpMainThm \in \ANNs$ such that for all $t \in [0,T]$, $x\in\R^d$ it holds that 
	\begin{equation}\label{eq:assG1}
	\textstyle{\functionANN{\fa_\indexAct}(\ANNassumpMainThm) \in  C(\R^d,\R), \qquad 
		\paramANN(\ANNassumpMainThm) \le \constantAssumpMainThm d^p \varepsilon^{-\gamma},} \qquad \text{and}
	\end{equation}
	\begin{equation}\label{eq:assG2}
	\textstyle{
		\varepsilon \abs{ u_d(t,x) } + \abs{ u_d(T,x) - (\functionANN{\fa_\indexAct}(\ANNassumpMainThm))(x) } 
		\le \varepsilon \constantAssumpMainThm d^p (1 + \norm{x})^p } 
	\end{equation}
	\cfload.
	Then there exists 
	$c\in \R$ 
	such that for all 
	$d\in\N$, 
	$\varepsilon\in(0,1]$ 
	there exists $\ANNresultMainThm\in\ANNs$
	such that 
	\begin{equation}
	\functionANN{\fa_\indexAct}(\ANNresultMainThm) \in C(\R^{d},\R), \qquad 
	\paramANN(\ANNresultMainThm) \le c d^c \varepsilon^{-c}, 
	\qquad \text{and}
	\end{equation}
	\begin{equation}\label{eq:cor:approx}	
		\textstyle{
		\sup_{q\in (0,\fq]}
		\bigl[ \int_{\R^d} \abs[]{ u_d(0,x) - (\functionANN{\fa_\indexAct}(\ANNresultMainThm))(x) }^q \mu_d(\dxx x)  \bigr]^{\nicefrac{1}{q}} 
		\le 
		\varepsilon }
	\dpp
	\end{equation}
\end{athm}

\begin{aproof}
	Throughout this proof 
	let $\fsc\colon \R\to \R$ satisfy for all $z\in\R$ that 
	$\fsc(z)=(2\fc)^{-1} f(z)$ 
	and for every $d\in\N$ let 
	$\usc_d\colon [0,2\fc T]\times\R^d \to \R$ 
	satisfy for all 
	$x\in\R^d$, $s\in[0,2\fc T]$ that 
	\begin{equation}
		\usc_d(s,x)=u_d\big((2\fc)^{-1}s,x\big) .
	\end{equation}		
	\Nobs that  
	\cref{eq:assG1}
	and~\cref{eq:assG2} 
	assure that 	
	there exist $(\F_{d,\varepsilon})_{(d,\varepsilon) \in \N\times(0,1]} \subseteq \ANNs$ which satisfy 	
	for all $d\in\N$, $\varepsilon\in(0,1]$, $t\in[0,T]$, $x\in\R^d$ that 
	$\functionANN{\fa_\indexAct}(\F_{d,\varepsilon}) \in  C(\R^d,\R)$, 
	$\paramANN(\F_{d,\varepsilon}) \le \constantAssumpMainThm d^p \varepsilon^{-\gamma}$, and 
	\begin{equation}\label{eq:1202}
		\varepsilon \abs{ u_d(t,x) } + \abs{ u_d(T,x) - (\functionANN{\fa_\indexAct}(\F_{d,\varepsilon}))(x) } 
		\le \varepsilon \constantAssumpMainThm d^p (1 + \norm{x})^p .
	\end{equation}
	\Nobs that \cref{eq:1202} 
	proves that for all $d\in\N$, $t\in[0,T]$, $x\in\R^d$ it holds that 
	\begin{equation}
	\abs{u_d(t,x)} \le \constantAssumpMainThm d^p (1 + \norm{x})^p .
	\end{equation}
	This assures that for all $d\in\N$ it holds that $u_d$ is at most polynomially growing.
	Therefore, we obtain that for all $d\in\N$ it holds that $\usc_d$ is at most 
	polynomially growing. 	
	\Moreover it holds for all $d\in\N$ that 
	$\usc_d \in C^{1,2}([0,2\fc T]\times \R^d,\R)$. 
	In addition, \cref{final_cor_eq1} implies for all $d\in\N$, $x\in\R^d$, $s \in [0,2\fc T]$ that 
	\begin{equation}
	(\tfrac{\partial}{\partial  t}\usc_d)(s,x) 
	+ 
	\tfrac{1}{2} (\Delta_x \usc_d)(s,x) 
	+ 
	\fsc\pr[\big]{ \usc_d(s,x) } 
	= 
	0 .
	\end{equation} 
	\Nobs that the assumption that $f$ is Lipschitz continuous establishes that there exists $L\in[0,\infty)$ 
	such that for all $w,z\in\R$ it holds that 
	$\abs{ f(z)-f(w) } \le L \abs{ z-w }$.  
	This yields for all $w,z\in\R$ that 
	$\abs{ \fsc(z)-\fsc(w) } \le (2\fc)^{-1} L \abs{ z-w }$. 
	\Moreover[Next] \cref{eq:1202} and the triangle inequality ensure for all $d\in\N$, $x\in\R^d$, $\varepsilon \in (0,1]$ that
	\begin{align}
	& \varepsilon \abs{(\functionANN{\fa_\indexAct}(\F_{d,\varepsilon}))(x)}
	+ \abs{\usc_d(2\fc T,x)-(\functionANN{\fa_\indexAct}(\F_{d,\varepsilon}))(x)}  \nonumber \\
	& = \varepsilon \abs{(\functionANN{\fa_\indexAct}(\F_{d,\varepsilon}))(x)}
	+ \abs{u_d(T,x)-(\functionANN{\fa_\indexAct}(\F_{d,\varepsilon}))(x)}  \nonumber 
	\le \varepsilon \abs{u_d(T,x) }
	+ 2 \abs{u_d(T,x)-(\functionANN{\fa_\indexAct}(\F_{d,\varepsilon}))(x)} \nonumber \\
	& \le 2 \varepsilon \constantAssumpMainThm d^p (1 + \norm{x})^p .
	\end{align}	
	\Moreover 
	the fact  
	that for all $d\in\N$, $\varepsilon\in(0,1]$ it holds that $\paramANN(\F_{d,\varepsilon}) \le \constantAssumpMainThm d^p \varepsilon^{-\gamma}$ and \cref{lem:dims_and_params} establish that for all $d\in\N$, $\varepsilon\in(0,1]$ it holds that 
	$\varepsilon^{\gamma}\lengthANN(\F_{d,\varepsilon})+\varepsilon^{\gamma}\normmm{\dims(\F_{d,\varepsilon})}\le \constantAssumpMainThm d^p$. 	
	\Moreover 
	it follows from \cref{interpol_ANN_function_implicit_max,interpol_ANN_function_implicit_softplus} 
	that there exists $(\F_{0,\varepsilon})_{\varepsilon\in (0,1]} \subseteq \ANNs$ such that for all $\varepsilon\in (0,1]$ the following properties hold true: 
	\begin{enumerate}[label=(\Roman *)]
		\item 
		\label{interpol_ANN_function_implicit:pfitem1}
		it holds that $\functionANN{\fa_\indexAct}(\F_{0,\varepsilon}) \in C(\R,\R)$, 
		\item
		\label{interpol_ANN_function_implicit:pfitem3}
		it holds for all $x,y\in\R$ that $\abs{(\functionANN{\fa_\indexAct}(\F_{0,\varepsilon}))(x) - (\functionANN{\fa_\indexAct}(\F_{0,\varepsilon}))(y)} \le \frac{L}{2\fc}\abs{x-y}$,
		\item
		\label{interpol_ANN_function_implicit:pfitem4}
		it holds for all $x\in\R$ that $\abs{(\functionANN{\fa_\indexAct}(\F_{0,\varepsilon}))(x)-\fsc(x)} \le 2 \varepsilon \max\{1,\abs{x}^2\}$, 
		\item
		\label{interpol_ANN_function_implicit:pfitem6}
		it holds that 
		$\paramANN(\F_{0,\varepsilon}) \le 24  (\max\{1,\fc^{-1}L\})^{2} \varepsilon^{-2} $ .
	\end{enumerate}
	Note that 
	\cref{interpol_ANN_function_implicit:pfitem6} and \cref{lem:dims_and_params} 
	prove  
	for all $\varepsilon \in (0,1]$ that 
	$\varepsilon^2 \lengthANN(F_{0,\varepsilon}) + \varepsilon^2 \normmm{\dims(\F_{0,\varepsilon})} \le 24  (\max\{1,\fc^{-1}L\})^{2} \varepsilon^{-2}$. 
	Furthermore, observe that \cref{interpol_ANN_function_implicit:pfitem3,interpol_ANN_function_implicit:pfitem4} 
	and the triangle inequality 
	imply for all $\varepsilon\in(0,1]$, $x \in \R$ that 
	\begin{align}
	\abs{ \left( \functionANN{\fa_\indexAct}(\F_{0,\varepsilon}) \right)(x) }
	& \le 
	\abs{(\functionANN{\fa_\indexAct}(\F_{0,\varepsilon}))(x) - (\functionANN{\fa_\indexAct}(\F_{0,\varepsilon}))(0)}
	+ \abs{(\functionANN{\fa_\indexAct}(\F_{0,\varepsilon}))(0)-\fsc(0)}
	+ \abs{\fsc(0)} \nonumber \\
	& \le \tfrac{L}{2\fc} \abs{x} + 2 \varepsilon + \abs{\fsc(0)}
	\le \left(2+(2\fc)^{-1}(L+\abs{f(0)})\right) (1+\abs{x})^2 .
	\end{align}
	This, \cref{interpol_ANN_function_implicit:pfitem4}, the fact that $\forall x \in \R\colon 1+\abs{x}^2 \le (1+\abs{x})^2$, and the assumption that $p\ge 2$ demonstrate for all $\varepsilon\in (0,1]$, $x\in\R$ that 
	\begin{align}
	& \varepsilon \abs{(\functionANN{\fa_\indexAct}(\F_{0,\varepsilon}))(x)} + 
	\abs{\fsc(x) - (\functionANN{\fa_\indexAct}(\F_{0,\varepsilon}))(x)} \nonumber \\
	& \le \varepsilon \left( \left(2+(2\fc)^{-1}(L+\abs{f(0)})\right) (1+\abs{x})^2 + 2\max\{1,\abs{x}^2\} \right)\nonumber \\
	& \le \varepsilon (4+(2\fc)^{-1}(L+\abs{f(0)})) (1+\abs{x})^p .
	\end{align}
	Moreover, 
	\cref{item:lem:Relu:dims,item:lem:ReluLeaky:dims,item:lem:ReluLeaky:real}
	in \cref{lem:ReluLeaky:identity} 
	and 
	\cref{item:lem:SoftPlus:real} 
	in \cref{lem:SoftPlus:identity}
	ensure that there exists 
	$\fJ \in \ANNs$ such that $\hiddenLength(\fJ)=1$ and $\functionANN{\fa_\indexAct}(\fJ)=\operatorname{id}_\R$. 
	\cref{theorem:final} 
	(applied with $\constantAssumpMainThm \with \max\{2\constantAssumpMainThm,24(\max\{1,\fc^{-1} L\})^2,4+(2\fc)^{-1}(L+\abs{f(0)})\}$,  $(u_d)_{d\in\N} \with (\usc_d)_{d\in\N}$, $L \with (2\fc)^{-1} L$, $\alpha_0 \with 2$, $\beta_0 \with 2$, $\alpha_1 \with \gamma$, $\beta_1 \with \gamma$, $T \with 2\fc T$, $a\with \fa_\indexAct$, $f_0 \with \fsc$, $(f_d)_{d\in\N} \with (\R^d \ni x \mapsto u_d(T,x)  \in \R)_{d\in\N}$, $(\nu_d)_{d\in\N}\with(\mu_d)_{d\in\N}$ in the notation of \cref{theorem:final}) 
	establishes 
	that there exist 
	$(\U_{d,\varepsilon})_{(d,\varepsilon)\in \N \times (0,1]} \subseteq \ANNs$ 
	and
	$c\in(0,\infty)$
	which satisfy for all 
	$d\in\N$, 
	$\varepsilon\in(0,1]$ 
	that 
	$\functionANN{\fa_\indexAct}(\U_{d,\varepsilon}) \in C(\R^{d},\R)$, 
	$\paramANN(\U_{d,\varepsilon}) \le c d^c \varepsilon^{-c}$, 
	and 	
	\begin{equation}\label{eq:cor:proof:approxV}
	\left( \int_{\R^d} \abs[\big]{ \usc_d(0,x) - (\functionANN{\fa_\indexAct}(\U_{d,\varepsilon}))(x) }^\fq \mu_d(\dxx x) \right)^{\!\!\nicefrac{1}{\fq}} 
	\le 
	\varepsilon 
	\dpp
	\end{equation} 
	The fact that for all $d \in \N$, $x \in \R^d$ it holds that $\usc_d(0,x)=u_d(0,x)$ and 
	Jensen's inequality hence imply \cref{eq:cor:approx}. 
\end{aproof}

\cfclear
\begin{athm}{corollary}{cor:final}
	Let 
	$T,\constantAssumpMainThm,\fc,p \in(0,\infty)$, 
	$a\in \R$, $b\in(a,\infty)$, 
	let $f\colon \R\to\R$ be Lipschitz continuous, 
	for every $d\in\N$ let $u_d \in C^{1,2}([0,T]\times \R^d,\R)$ 
	satisfy for all $t \in [0,T]$, $x\in\R^d$ that 
	\begin{equation}\label{final_cor_pde_reversed}
	\textstyle{(\tfrac{\partial}{\partial  t}u_d)(t,x) 
	= 
	\fc (\Delta_x u_d)(t,x) 
	+ 
	f\pr[\big]{ u_d(t,x) } 
	, }
	\end{equation}
	let 
	$\indexAct\in\{0,1\}$, 
	$\alpha \in [0,\infty)\backslash\{1\}$,   
	$\fa_0,\fa_1\in C(\R,\R)$ satisfy for all $x\in\R$ that $\fa_0(x) = \max\{x,\alpha x\}$ and $\fa_1(x)=\ln(1+\exp(x))$, 
	and assume for all $d\in\N$, $\varepsilon \in (0,1]$ 
	that there exists 
	$\ANNassumpMainThm \in \ANNs$ such that for all $t \in [0,T]$, $x\in\R^d$ it holds that 
	\begin{equation}
	\textstyle{
		\functionANN{\fa_\indexAct}(\ANNassumpMainThm) \in  C(\R^d,\R), \qquad 
		\paramANN(\ANNassumpMainThm) \le \constantAssumpMainThm d^{\constantAssumpMainThm} \varepsilon^{-\constantAssumpMainThm}, } \qquad \text{and}
	\end{equation}
	\begin{equation}\label{eq:0931}
	\textstyle{
		\varepsilon \abs{ u_d(t,x) } + \abs{ u_d(0,x) - (\functionANN{\fa_\indexAct}(\ANNassumpMainThm))(x) } 
		\le \varepsilon \constantAssumpMainThm d^{\constantAssumpMainThm} (1 + \norm{x})^{\constantAssumpMainThm}  }
	\end{equation}
	\cfload.
	Then there exists 
	$c\in \R$ 
	such that for all 
	$d\in\N$, 
	$\varepsilon\in(0,1]$ 
	there exists $\ANNresultMainThm\in\ANNs$
	such that 	
	\begin{equation}
	\textstyle{\functionANN{\fa_\indexAct}(\ANNresultMainThm) \in C(\R^{d},\R), \qquad
	\paramANN(\ANNresultMainThm) \le c d^c \varepsilon^{-c}}, 
	\qquad \text{and}
	\end{equation}
	\begin{equation}\label{eq:cor:approxLp}
		\textstyle{
		\bigl[ \frac{1}{(b-a)^d} \int_{[a,b]^d} \abs{ u_d(T,x) - (\functionANN{\fa_\indexAct}(\ANNresultMainThm))(x) }^p \dx x \bigr]^{\nicefrac{1}{p}} 
	\le 
	\varepsilon }
	\dpp
	\end{equation}
\end{athm}

\begin{aproof}
	Let $\fp = \inf\{k\in\N\colon k\ge \max\{\constantAssumpMainThm,2\}\}$ and $\fq = \max\{p,2\}$. 
	For every $d\in\N$ let $v_d\colon [0,T]\times\R^d\to\R$ satisfy for all $x\in\R^d$, $t\in[0,T]$ that $v_d(t,x)=u_d(T-t,x)$. 
	\Nobs that \cref{final_cor_pde_reversed} shows that 
	it holds for all $d\in\N$, $x\in\R^d$, $t\in[0,T]$ that 
	\begin{equation}
	(\tfrac{\partial}{\partial  t}v_d)(t,x) 
	= -  
	\fc (\Delta_x v_d)(t,x) 
	- 
	f\pr[\big]{ v_d(t,x) } 
	.
	\end{equation}
	\Nobs that \cref{eq:0931} implies for all $d\in\N$, $\varepsilon\in(0,1]$, $t\in[0,T]$, $x\in\R^d$ that 
	\begin{equation}
	\varepsilon \abs{v_d(t,x)} + \abs{ v_d(T,x) - (\functionANN{\fa_\indexAct}(\ANNassumpMainThm))(x) } 
	\le \varepsilon \constantAssumpMainThm d^{\constantAssumpMainThm} (1 + \norm{x})^{\constantAssumpMainThm} .
	\end{equation}
	For all $d\in\N$ let $\mu_d\colon \cB(\R^d)\to[0,1]$ be the uniform distribution on $[a,b]^d$ and note that there exists $K\in[1,\infty)$ such that for all $d\in\N$ it holds that 
	$\int_{\R^d} \norm{y}^{\fp^2\fq} \mu_{d}(\dxx y) \le K^{\fp^2\fq
	} d^{\fp^2\fq}$.
	\cref{cor:almostfinal} (applied with 
	$p\with \fp$, 
	$\fq\with \fq$,
	$\gamma\with \constantAssumpMainThm$,
	$r\with 1$,
	$\constantAssumpMainThm \with \max\{\constantAssumpMainThm,K^{\fp^2\fq}\}$, 
	$(u_d)_{d\in\N}\with(v_d)_{d\in\N}$  
	in the notation of \cref{cor:almostfinal}) 
	establishes 
	that there 
	exists 
	$c\in \R$ 
	such that for all 
	$d\in\N$, 
	$\varepsilon\in(0,1]$ 
	there exists $\ANNresultMainThm\in\ANNs$ 
	such that  
	$\functionANN{\fa_\indexAct}(\ANNresultMainThm) \in C(\R^{d},\R)$, 
	$\paramANN(\ANNresultMainThm) \le c d^c \varepsilon^{-c}$, 
	and 	
	\begin{equation}
	\sup_{q\in(0,\fq]}\left( \int_{\R^d} \abs[\big]{ v_d(0,x) - (\functionANN{\fa_\indexAct}(\ANNresultMainThm))(x) }^q \mu_d(\dxx x) \right)^{\!\!\nicefrac{1}{q}} 
	\le 
	\varepsilon 
	\dpp
	\end{equation} 
	This, the definition of $(\mu_d)_{d\in\N}$, the fact that $p\in (0,\fq]$, and the fact that for all $d\in\N$, $x\in\R^d$ it holds that $v_d(0,x)=u_d(T,x)$ prove \cref{eq:cor:approxLp}. 
\end{aproof}


\section*{Acknowledgments}

This work has been partially funded 
by the Deutsche Forschungsgemeinschaft (DFG, German Research Foundation) under Germany's Excellence Strategy EXC 2044-390685587, Mathematics M{\"u}nster: Dynamics-Geometry-Structure.
In addition, this work has been partially 
funded by the Deutsche Forschungsgemeinschaft (DFG, German Research Foundation) in the frame of the priority programme SPP 2298 ``Theoretical Foundations of Deep Learning'' -- Project no.\ 464123384 and Project no.\ 464101154. 
This work has also been supported by the Ministry of Culture and Science NRW as part of the Lamarr Fellow Network.
Furthermore, funding 
by the National Science Foundation (NSF 1903450)
is gratefully acknowledged.
Moreover, the second author gratefully acknowledges the support of 
the startup fund project of the Shenzhen Research Institute of Big Data under grant No.\ T00120220001.

\bibliographystyle{acm}
\bibliography{bibfile}

\end{document}